\documentclass{amsart} 
\usepackage[pdftex]{graphicx}  
\usepackage{mathtools,amssymb} 
\usepackage{amsthm} 
\usepackage{mathrsfs} 
\usepackage{amscd} 
\usepackage{tikz} 
\usepackage{hyperref} 
\usepackage{cleveref} 
\usepackage{enumitem} 
\usepackage{multicol} 
\usepackage{ascmac} 
\usepackage{cases} 
\usepackage{bm} 
\usepackage[mark=o]{dynkin-diagrams} 
\usepackage{longtable} 
\usepackage{colortbl,xcolor} 
\usepackage{rotating} 
\usepackage{here} 
\usepackage{ulem} 
\usepackage{lscape}
\usepackage{circledsteps} 
\usepackage{caption}
\usetikzlibrary{nfold}

\DeclareMathOperator{\ch}{char} 
\newcommand{\Z}{\mathbb{Z}} 

\renewcommand{\O}{\mathcal{O}} 
 
\DeclareMathOperator{\ord}{ord}

\DeclareMathOperator{\Spec}{Spec}

\newcommand{\A}{\mathbb{A}} 
\renewcommand{\P}{\mathbb{P}}
\DeclareMathOperator{\Pic}{Pic}

\DeclareMathOperator{\NS}{NS}
\DeclareMathOperator{\Supp}{Supp}

\DeclareMathOperator{\Triv}{Triv}

\theoremstyle{definition}
\newtheorem{def.}{Definition}[section] 
\theoremstyle{definition}
\newtheorem{prop.}[def.]{Proposition}
\theoremstyle{definition}
\newtheorem{th.}[def.]{Theorem}
\theoremstyle{definition}
\newtheorem{cor.}[def.]{Corollary} 
\theoremstyle{definition}
\newtheorem{rem.}[def.]{Remark}
\theoremstyle{definition}
\newtheorem{ex.}[def.]{Example}
\theoremstyle{definition}
\newtheorem{lem.}[def.]{Lemma}
\theoremstyle{definition}
\newtheorem{fact.}[def.]{Fact}
\theoremstyle{definition}
\newtheorem{cond.}[def.]{Condition}  
\theoremstyle{definition}

\address{DEPARTMENT OF MATHEMATICS, FACULTY OF SCIENCE AND TECHNOLOGY, 
TOKYO UNIVERSITY OF SCIENCE, 2641 Yamazaki, Noda, Chiba, 278-8510, Japan}  
\email{ito\_hiroyuki@rs.tus.ac.jp} 

\address{DEPARTMENT OF MATHEMATICS, FACULTY OF SCIENCE AND TECHNOLOGY, 
TOKYO UNIVERSITY OF SCIENCE, 2641 Yamazaki, Noda, Chiba, 278-8510, Japan}  
\email{shota140204@gmail.com} 
\email{6124514@ed.tus.ac.jp} 

\date{August 24, 2025} 

\thanks{This work was supported by JSPS KAKENHI Grant number JP22K03254.}  

\begin{document} 

\tikzset{cross/.style={preaction={-,draw=white,line width=6pt}}} 

\title{Unirational quasi-hyperelliptic surfaces 
in characteristic five}    

\author{Hiroyuki Ito \and Shota Takayashiki}

\begin{abstract} 
As a generalization of a quasi-elliptic surface,  
there is a quasi-hyperelliptic surface, a nonsingular projective surface which has a fibration structure whose general fiber is a quasi-hyperelliptic curve ($=$ singular hyperelliptic curve 
with one cuspidal singular point) of genus $(p-1)/2$ 
in characteristic $p \ge 5$. 
In this paper, we consider unirational quasi-hyperelliptic surfaces in characteristic $5$, and classify singular fibers and give a formula for the arithmetic genus and the self-intersection number of the canonical divisor.   
As a corollary, we classify rational quasi-hyperelliptic surfaces by determining the combinations of singular fibers, the defining equations and sections.    
\end{abstract} 

\maketitle

\section{Introduction}  
In the study of algebraic varieties in positive characteristic,
the existence of quasi-fibrations, 
which are objects specific to algebraic varieties in positive characteristic,
is a peculiar feature due to the failure of Bertini’s theorem.
The existence of such objects may be dismissed as pathological phenomena, 
such as exceptional varieties or phenomena in positive characteristic,
but the existence of varieties with such a fibration structure in which singular varieties appear as general fibers is extremely important in the study of algebraic geometry in positive characteristic.

A typical example of a variety with a quasi-fibration structure is a quasi-elliptic surface,
which played an important role in the classification of algebraic surfaces in positive characteristic deeply studied by Bombieri-Mumford \cite{BM77,BM76}.  
Furthermore, Schr\"oer \cite{Schrooer08} studied fibrations in which general fibers are singular,
that is,
varieties with a quasi-fibrational structure, and obtained interesting results,
especially regarding singular points that appear in general fibers.

The generic fiber of such a fibration is a regular variety over the function field of the base variety,
but it becomes a singular variety over a purely inseparable extension field. 
Therefore,
when a general fiber is one-dimensional,
the genus of the curve over the function field of the base variety is different from the genus over an algebraically closed field.
According to Tate \cite{T52}, 
the genus of the curve over an algebraically closed field is less than the genus of a general fiber,
and the difference is divisible by $\frac{p-1}{2}$. 
This shows that even if the fiber is one-dimensional, 
quasi-fibrations of various genera exist in all positive characteristics.

When the arithmetic genus of a general fiber is one,
which is a special case of a quasi-fibration,
called a quasi-elliptic fibration,
there has been much research on it both geometrically and number theoretically.
(There has been much research on quasi-elliptic surfaces, 
starting with Bombieri-Mumford \cite{BM77,BM76}, also see Miyanishi-Ito \cite{MI21}.)

In this paper, 
we deal with quasi-fibrations defined on nonsingular projective algebraic surfaces, 
whose general fiber is a singular rational curve of arithmetic genus two (so the characteristic of the ground field is $5$).   
Let us call such surfaces, quasi-hyperelliptic surfaces.  
More generally, we will be back to the case where a general fiber is a singular rational curve of arithmetic genus $g=\frac{p-1}{2}$ in a forthcoming paper.

When one take a enough purely inseparable base extension of the general fiber of such a surface,
a one-place singular point appears on the curve, 
and the affine curve excluding this singular point becomes isomorphic to an affine line after the field extension,
so it is called the $\A^{1}$-form. 
The $\A^{1}$-form of genus $\geq 1$ has been studied deeply by Kambayashi, Miyanishi and Takeuchi \cite{KMT74} 
and Kambayashi and Miyanishi \cite{KM77}.     
In addition, Miyanishi \cite{M78} (see also Miyanishi-Ito \cite{MI21}) has obtained partial results on quasi-hyperelliptic fibrations in characteristic 5, 
although with technical restrictions.

The first main theorem (Theorem \ref{main.th.1}) is a local result for the quasi-hyperelliptic fibration,
which gives the classification of the configurations of singular fibers that appear
in a quasi-hyperelliptic fibration of genus 2 by using a defining equation.  

\begin{th.}[{Theorem \ref{main.th.1}}]   
  Let $f : X \to \Spec(k[[t]])$ be a local quasi-hyperelliptic fibration of genus 2 with a section. 
  Then $X$ is defined by $y^{2} = x^{5} + t^{e}u$, where 
  $e \in \{ 1,2,3,4,6,7,8,9 \}$ and $u$ is a unit of $k[[t]]$. 
  Furthermore, the fiber is classified into $8$ types 
  corresponding to the value $e$. 
  \begin{table}[H] 
    \centering 
    \scalebox{0.93}{ 
    \begin{tabular}{|c|c|c|c|c|c|c|c|c|} 
    \hline 
    $e$ &  1 & 2 & 3 & 4 & 6 & 7 & 8 & 9   
    \\ \hline  
    Type & $C(1,0)$ & $C(5,1)$ & $C(9,2)$ & $C(3,3)$ & 
    $C(11,5)$ & $C(4,6)$ & $C(9,7)$ & $C(13,8)$
    \\ \hline   
    \end{tabular}   
    } 
  \end{table}  
  The configurations and topological Euler numbers of these $C(-,n)$'s are given as in Table \ref{table:configurations of fibers}.
\end{th.}   

The second main theorem is a global result when the base curve is a projective line,
or equivalently, 
the surface is unirational. 
We remove the technical restrictions of Miyanishi's result and give a formula for the arithmetic genus and the self-intersection number of the canonical divisor from the information on the combination of singular fibers.  

\begin{th.}[{Theorem \ref{main.th.2}}] 
  Let $f : X \to \P^{1}$ be a relatively minimal unirational quasi-hyperelliptic fibration 
  of genus $2$ in characteristic $5$, which is birationally equivalent to 
  an affine hypersurface defined by $y^{2} = x^{5} + \varphi(t)$  
  satisfying Condition \ref{assumption}.   
  Put $d = \deg_{t}\varphi(t)$ and $m = [d/10]$, 
  where $[d/10]$ is the greatest integer less than or equal to $d/10$. 
  Then, we have    
  \begin{align*} 
    p_{a}(X) &= 4m+2 -(N_{3}+N_{5}+2N_{6}+2N_{7}+2N_{8}), \\ 
    (K_{X}^{2}) &= 8m-(2N_{3}+2N_{5}+3N_{6}+4N_{7}+4N_{8}),   
  \end{align*}
  where $N_{i}$ is the number of reducible fibers of $f$ of type $C(-,i)$.    
\end{th.}

Furthermore, 
as an application of these results,
when the surface is rational, 
that is, 
rational surfaces with a quasi-hyperelliptic fibration structure of genus 2,
we classify the combinations of singular fibers, 
defining equations, and sections of the fibrations, 
as well as the Mordell-Weil groups of the generalized Jacobian variety of the generic fiber.

The method used to prove both the local and global results 
is the double coverings and canonical resolution,
deelply studied by M.\ Artin in his doctoral thesis, 
and the results are obtained by calculating precise branch loci and coverings. 
In addition, 
Shioda's Theory of Mordell-Weil lattices for higher genus fibrations 
makes it possible to apply two main theorems to specific surfaces. 

The structure of the entire paper is as follows. 
After this introduction, we recall the definition of a quasi-hyperelliptic fibration 
and recall also the theory of Mordell-Weil lattices for higher genus fibrations. 
In Section \ref{sec:Main Th}, we state two main theorems, local and global ones. 
After preparing the setting for the proofs 
in Seciton \ref{sec:On double covering}, 
we give proofs of two main theorems 
in Sections \ref{sec:pf.of.main.th.1} and 
\ref{sec:pf.of.main.th.2} each. 
Finally, we apply these results to the rational surfaces and 
give the classification theorem in Section \ref{sec:rat QE}.  

\section{A quasi-hyperelliptic fibration} \label{sec:q.h.e. fib} 

Throughout this paper, we work over an algebraically closed field $k$ of characteristic $p > 0$ 
unless otherwise specified.   
We recall the definition of a hyperelliptic field.  

\begin{def.} 
  A field $K$ is {\itshape hyperelliptic over $k$} if 
  it satisfies the following conditions;  
  \begin{enumerate}[label=(\roman*),topsep=1pt, labelsep=3pt]
    \item $K$ is a one variable function field over $k$ and its $k$-genus is greater than or equal to $2$. 
    Here, the ($k$-)genus of an algebraic function field $K/{k}$ is defined to be the dimension of the first cohomology group of the normalized curve corresponding to $K$.   
    \item $K$ contains a subfield $K_{0}$ which contains $k$ and its $k$-genus is $0$.  
    \item The degree of the extension $K/{K_{0}}$ is equal to $2$.
  \end{enumerate} 
\end{def.} 

\begin{def.}
  Let $K$ be a field. 
  A $K$-form $X$ of the affine line $\A_{K}^{1}$ is a {\itshape hyperelliptic $K$-form of genus $g$}
  if its function field $k(X)$ is hyperelliptic of genus $g$ over $k$. 
  Here, a $K$-form of $\A_{K}^{1}$ is, by definition,   
  a scheme $X$ over $K$ such that $X \times_{K} \overline{K}$ is $\overline{K}$-isomorphic to $\A_{\overline{K}}^{1}$ for an algebraic closure $\overline{K}$ of $K$.  
\end{def.}

\begin{lem.}[{\cite[Theorem I.4.3]{MI21}}] \label{equation of a hyperelliptic form} 
  Let $K$ be a non-perfect field of characteristic $p \ge 3$. 
  Let $X$ be a hyperelliptic $K$-form of $\A_{K}^{1}$ of $k$-genus 
  $g \ge 2$. 
  Assume that $X$ has a $K$-rational point.  
  Then $X$ is $K$-biraitionally equivalent to an affine plane curve 
  of the type 
  \[ 
  y^{2} = x^{p^{m}} + a , \ \ a \in K - K^{p}, 
  \]  
  for which $g = \frac{1}{2}(p^{m}-1)$.  
  Conversely, the complete $K$-normal model $C$ of 
  every such plane curve has a unique singular point 
  $P$ and the smooth part $C_{\mathrm{sm}} = C \setminus \{ P \}$ of $C$  
  is a $K$-form of $\A_{K}^{1}$ of $K$-genus $\frac{1}{2}(p^{m}-1)$. 
\end{lem.}

\begin{def.}   
  Let $f : X \to C$ be a fibration over $k$.   
  Assume that $\ch(k) = p \ge 5$.  
  A fibration $f : X \to C$ is a 
  {\itshape quasi-hyperelliptic fibration of genus $g=(p-1)/2$} 
  if the generic fiber is a hyperelliptic $K$-form of $\A_{K}^{1}$ 
  with arithmetic genus $g = (p-1)/2$, 
  where $K$ is the function field of $C$. 
  We sometimes call such $X$ a {\itshape quasi-hyperelliptic surface}.        
\end{def.} 

\begin{rem.} 
  We remark on the term ``quasi-hyperelliptic surface''.   
  This term is used in Bombieri and Mumford's papers 
  \cite{BM77,BM76} in different meaning. 
  They uses this term for 
  a surface $X$ of Kodaira dimension zero with $b_{1}(X)=2$, $b_{2}(X)=2$ and $\chi(\O_{X})=0$ 
  whose Albanese map defines a quasi-elliptic fibration, 
  where $b_{i}(X)$ is the $i$-th Betti number of $X$.    
\end{rem.} 

The next lemma follows from 
Lemma \ref{equation of a hyperelliptic form}. 

\begin{lem.}[{\cite[Lemma II.3.1.3]{MI21}}] \label{standard form of QE} 
  Assume that the characteristic of the base field $k$ is 
  greater than $3$. 
  Let $f : X \to C$ be a quasi-hyperelliptic fibration of 
  genus $g = (p-1)/2$. 
  Let $X_{\eta}$ be the generic fiber of $f$. 
  Let $K$ be the function field of $C$. 
  Then we have : 
  \begin{enumerate}[label=(\roman*)] 
    \item The generic fiber $X_{\eta}$ 
    has a unique singular point $P_{\infty}$ which is a one-place point. 
    Furthermore, $X_{\eta} \setminus \{ P_{\infty} \}$ is a $K$-form 
    of the affine line $\A_{K}^{1}$. 
    Hence almost all fibers of $f$ are irreducible 
    rational curves wtih unique one-place singular point, 
    i.e., cuspidal singular point. 
    \item assume that $f$ has a section.  
    Then, $X_{\eta}$ is birationally equivalnt to an affine plane 
    curve defined by 
    \[ 
      y^{2} = x^{p} + a \ \ \text{with} \ \ a \in K - K^{p}. 
    \]
  \end{enumerate} 
  The closure of $\{ P_{\infty} \}$ in the surface $X$ is called {\itshape moving cusp} of $f$  
  just like for a quasi-elliptic fibration 
  and we denote it by $\Gamma$  
\end{lem.} 

Recall that a $n$-dimensional variety $X$ defined over $k$ is unirational if 
there exists a dominant rational map from $\P_{k}^{n}$ to $X$.  
For a uniraional surface, the following property is well-known. 

\begin{lem.} \label{irregularity}
Let $X$ be a unirational surface. 
Then the irregularity of $X$ is equal to zero, where 
the irregularity of $X$ is defined to be the dimension of 
the Albanese variety $\mathrm{Alb}(X)$ of $X$.    
\end{lem.}

\begin{proof}
By the definition of the unirationality, there exists 
a dominant rational map from a projective plane $\P^{2}$ to a surface $X$.  
Thus, we have a dominant morphism from a rational surface $Y$ to $X$. 
This morphism induces a surjective homomorphism 
$\mathrm{Alb}(Y) \to \mathrm{Alb}(X)$ and 
$\mathrm{Alb}(Y)$ is a point 
since $Y$ is rational. 
Therefore, we have $q(X) = \dim \mathrm{Alb}(X) = 0$. 
\end{proof}  

For a quasi-hyperelliptic surface, 
we have the following characterization of unirationality. 

\begin{lem.}
  Let $f : X \to C$ be a quasi-hyperelliptic fibration of genus $(p-1)/{2}$ 
  in characteristic $p \ge 5$ with a section.  
  Then, the following conditions are equivalent. 
  \begin{enumerate}[label=(\roman*)]
    \item $X$ is unirational. \label{lem:characterization of unirat 1} 
    \item $C$ is $k$-rational, that is $C \cong \P_{k}^{1}$. \label{lem:characterization of unirat 2} 
  \end{enumerate} 
\end{lem.}

\begin{proof} 
\ref{lem:characterization of unirat 1} $\implies$ \ref{lem:characterization of unirat 2}: 
This implication follows from L\"{u}roth's theorem. 

\ref{lem:characterization of unirat 2} $\implies$ \ref{lem:characterization of unirat 1}:  
Assume that $C$ is $k$-rational. 
Then, the function field $K$ of $C$ is isomorphic to $k(t)$ and 
$X_{K} \otimes_{K} K^{1/{p}}$ is rational over $K^{1/{p}}$  
by Lemma \ref{standard form of QE}, where $X_{K}$ is the generic fiber of $f$.   
Hence $K^{1/{p}}(X_{K})$ is isomorphic to the polynomial ring over $K^{1/{p}}$, 
that is, $K^{1/{p}}(X_{K}) \cong K^{1/{p}}(u)$.   
Since $K^{1/{p}} = k(t^{1/{p}})$, we know that $k(X) \hookrightarrow k(t^{1/{p}},u)$. 
Hence $X$ is unirational.    
\end{proof}

\section{The Mordell-Weil groups} \label{sec:MWG} 

In this section, we consider 
a unirational quasi-hyperelliptic fibration $f : X \to \P^{1}$ 
of genus $g = (p-1)/{2}$ in characteristic $p \ge 5$. 
Let $X_{\eta}$ be the generic fiber of $f$ with the unique singular point $P_{\infty}$
and $K$ be the function field of $\P^{1}$. 

\begin{def.} 
  Let $J_{X_{\eta}}$ be the generalized Jacobian of $X_{\eta}$. 
  We define the {\itshape Mordell-Weil group} of $f$ to be 
  the group of $K$-rational points of $J_{X_{\eta}}$, and 
  denote it by $J_{X_{\eta}}(K)$.     
\end{def.} 

Then, there exists a closed immersion 
\[ 
  \iota : X_{\eta} \setminus \{ P_{\infty} \} 
  \hookrightarrow 
  J_{X_{\eta}} = \Pic_{X_{\eta}/{K}}^{0}  
\] 
and the image $\iota(X)$ generates $J_{X_{\eta}}$ as a $K$-group scheme (\cite[Theorem I.2.4.3]{MI21}). 

Note that this closed immersion $\iota$ induces an injection of $K$-rational points 
\[
  \iota(K) :  (X_{\eta} \setminus \{ P_{\infty} \})(K) \to J_{X_{\eta}}(K). 
\]

As (quasi-)elliptic fibrations, we use the following notations. 
\begin{itemize}
  \item $R = \{ v \in \P^{1}(k) \mid F_{v} \coloneqq f^{-1}(v) \ \text{is reducible} \}$, 
  \item For a $v \in R$, we write 
  \[ F_{v} = R_{v,0} + \sum_{i=1}^{m_{v}-1}\mu_{v,i}R_{v,i}, \]
  where $R_{v,0}$ is the irreducible component meeting the zero section $(O)$, 
  $m_{v}$ is the number of irreducible components of $F_{v}$,  
  and $\mu_{v,i}$ is the multiplicity of the component $R_{v,i}$. 
  \item $\Triv(X)$ is the sublattice of $\NS(X)$ generated by 
  the zero section $(O)$, a general fiber $F$, and 
  the irreducible components $R_{v,i}$ for $v \in R$ and $1 \le i \le m_{v}-1$. 
  The sublattice $\Triv(X)$ is called the {\itshape trivial lattice} of $f$.  
\end{itemize} 

In the case that $X$ is a rational surface, 
the following theorem holds. 

\begin{th.}[{\cite[Theorem II.4.1.1]{MI21}}] \label{th:MWG of rat q-HE}  
Let $f : X \to \P^{1}$ be a rational quasi-hyperelliptic fibration with a section.  
Then the following assertions hold. 
\begin{enumerate}[label=(\roman*)]
  \item The N\'{e}ron-Severi group $\NS(X)$ of $X$ is a free abelian group of 
  finite rank endowed with a non-degenerate intersection pairing. 
  \item As an abelian group, the Mordell-Weil group 
  $J_{X_{\eta}}(K)$ is isomorphic to $\NS(X)/{\Triv(X)}$.   
  \item The Mordell-Weil group $J_{X_{\eta}}(K)$ is a finite-dimensional vector space 
  over $\mathbb{F}_{p}$ of dimension $r$, 
  and $|J_{X_{\eta}}(K)|^{2} = p^{2r} = \det(\Triv(X))/{\det(\NS(X))}$.   
\end{enumerate} 
We call this $r$ the {\itshape torsion-rank} of the Mordell-Weil group.  
\end{th.} 

\begin{cor.} \label{th:str of MWG}
  For a rational quasi-hyperelliptic fibration $f : X \to \P^{1}$ with a section, 
  we have 
  \[
    J_{X_{\eta}}(K) \cong (\Z/{p\Z})^{\oplus r}  
  \]
  for some $r \in \Z_{\ge 0}$.  
\end{cor.} 

Finally, we mention the height pairing on the Mordell-Weil group. 
As similar fashion in \cite{Shioda-Schutt19}, \cite{Sioda},   
we can define the height pairing $\langle \ , \ \rangle$ on $J_{X_{\eta}}(K)$. 
We omit the formula for general elements of $J_{X_{\eta}}(K)$, 
but we remark that, for elements which are the images of 
the generic fiber $(X_{\eta} \setminus \{ P_{\infty }\})$ under the closed immersion $\iota$,  
say $\iota(P)$ and $\iota(Q)$, we have   
\begin{equation} \label{explcit height pairing} 
  \langle \iota(P), \iota(Q) \rangle = -((P) - (O) \cdot (Q) - (O))-\sum_{v \in R}\mathrm{contr}_{v}(P,Q),  
\end{equation} 
\noindent where 
$(P)$ and $(Q)$ are the corresponding sections and 
$\mathrm{contr}_{v}$ is the local contribution at $v \in R$ defined by   
\[
  \mathrm{contr}_{v}(P,Q) 
  = \biggl(
    \bigl((P) \cdot R_{v,1}\bigr), \dots , \bigl((P) \cdot R_{v,m_{v}-1}\bigr) 
    \biggr) 
    (-A_{v}^{-1})
  \begin{pmatrix}
    ((Q) \cdot R_{v,1}) \\ 
    \vdots \\
    ((Q) \cdot R_{v,m_{v}-1})
  \end{pmatrix}, 
\]
where $A_{v}$ is the gram matrix of a reducible fiber $f^{-1}(v)$. 
Note that $\mathrm{contr}_{v}(P,Q) = 0$ 
if either $(P)$ or $(Q)$ intersects the identity component of the reducible fiber $f^{-1}(v)$.   
Furthermore, since every element of $J_{X_{\eta}}(K)$ is a torsion element, 
we have $\langle \iota(P),\iota(P) \rangle = 0$ for $P \in (X_{\eta} \setminus \{ P_{\infty} \})$ 
(cf.\ \cite{Sioda}).  

\section{The main theorems} \label{sec:Main Th}   

From now on, we consider the case $p=5$,  
that is, a quasi-hyperelliptic fibratiton of genus $2$. 
Let $k$ be an algebraically closed field of characteristic $5$,  
and $f : X \to C$ be a quasi-hyperelliptic fibration of 
genus $2$ on a unirational surface $X$ with a section. 
By the lemma \ref{standard form of QE}, 
we already know that $X$ is birationally equivalent to a hypersurface 
in $\A_{k}^{3}=\Spec(k[t,x,y])$ defined by $y^{2} = x^{5} + \varphi(t)$ 
with $\varphi(t) \in k[t]$. 
From a birational view point, 
we can assume that $\varphi(t)$ satisfies the following conditions. 

\begin{cond.} \label{assumption} 
  \begin{enumerate}[label= (\roman*)] 
    \item $\varphi(t)$ has no terms of degree multiple of $5$. 
    In particular, $\varphi(t)$ has no constant terms.  
    \item For every root $\alpha$ of $\varphi'(t) = 
    \frac{d\varphi}{dt} = 0$, 
    \[
    e_{\alpha} \coloneqq v_{\alpha}(\varphi(t)-\varphi(\alpha)) 
    \in \{ 2,3,4,5,6,7,8,9 \}, 
    \] 
    where, $v_{\alpha}$ is a $(t-\alpha)$-adic valuation of $k[t]$  
    with $v_{\alpha}(t-\alpha) = 1$.  
  \end{enumerate} 
\end{cond.} 

Indeed, 
if $\varphi(t)$ contains the term of degree multiple of $5$, 
a birational transformation of the type 
\[ 
  (t,x,y) \mapsto (t,x+\phi(t),y) 
\]
makes $\varphi(t)$ into a desired form with a suitable polynomial $\phi(t)$.   
Furthermore, if there exists a root $\alpha$ 
of the $\varphi'(t)$ whose $e_{\alpha}$ is greater than $9$, 
then we can drop $e_{\alpha}$ by a birational transformation of the type 
\[
  (t,x,y) \mapsto 
  \left(
    t,\frac{x}{(t-\alpha)^{2}},\frac{y}{(t-\alpha)^{5}}
  \right).
\]

Let us state the first main theorem which 
classifies non-multiple fibers of quasi-hyperelliptic fibrations.  

\begin{th.} \label{main.th.1} 
  Let $f : X \to \Spec(k[[t]])$ be a local fibration  
  defined by $y^{2} = x^{5} + t^{e}u$, where 
  $e \in \{ 1,2,3,4,6,7,8,9 \}$ and $u$ is a unit of $k[[t]]$.   
  Then, the fiber is classified into $8$ types
  corresponding to the value $e$. 
  The configurations and topological Euler numbers are given as in Table \ref{table:configurations of fibers}. 
\end{th.} 

\newpage 

\begin{table}[H] 
  \centering
  \scalebox{0.7}{
  \begin{tabular}{|c|c|c|c|} 
  \hline
  $e$ & 
  type & 
  configuration & 
  \begin{tabular}{c} 
    topological \\ Euler number 
  \end{tabular}
  \\ \hline 
  1 & 
  $C(1,0)$ &  
  \begin{tikzpicture}
  \draw[very thick] (-2,1) to [out=0,in=90]  (0,0) 
  to [out=90, in=-180] (2,1);
  \node  (cusp explain) at (1,0.2) {$(2,5)\text{-cusp}$};
  \node (selfint) at (2.2,1) {$\Circled{0}$}; 
  \draw[dashed] (0,1) -- (0,-0.5); 
  \end{tikzpicture}
  & $2$ 
  \\ \hline 

  2 & 
  $C(5,1)$ & 
  \begin{tikzpicture} 
    \draw[very thick] (-1.5,0) -- (1.5,0);
    \draw (-1,-1) -- (1,1); 
    \draw (-1,1) -- (1,-1); 
    \draw (0.5,0.7) -- (2.5,0.7); 
    \draw (0.5,-0.7) -- (2.5,-0.7); 
    \node (r0selfint) at (-1.8,0) {$\Circled{-4}$}; 
    \node (r1coeff) at (-0.65,0.9) {$2$}; 
    \node (r3coeff) at (-0.65,-0.9) {$2$}; 
    \draw[dashed] (0,1) -- (0,-1); 
  \end{tikzpicture} 
  & $6$ 
  \\ \hline 
  
  3 & 
  $C(9,2)$ & 
  \begin{tikzpicture} 
    \draw (-2,0) -- (2,0); 
    \draw (-1,0.2) -- (-1,-1.2); 
    \draw (-0.8,-1) -- (-2.2,-1);
    \draw[very thick] (-2,-0.8) -- (-2,-2.2); 
    \draw (0,0.2) -- (0,-1.8);   
    \draw (1,0.2) -- (1,-1.2); 
    \draw (0.8,-1) -- (2.2,-1);
    \draw (2,-0.8) -- (2,-2.2);  
    \draw (1.8,-2) -- (3.2,-2);
    \node (r0selfint) at (-2,-2.5) {$\Circled{-4}$}; 
    \node (r1coeff) at (-1.6,-0.8) {$4$}; 
    \node (r2coeff) at (-0.8,-0.5) {$7$}; 
    \node (r3coeff) at (-1.9,0.2) {$10$}; 
    \node (r4coeff) at (0.2,-0.5) {$5$}; 
    \node (r5coeff) at (1.2,-0.5) {$8$}; 
    \node (r6coeff) at (1.6,-0.8) {$6$}; 
    \node (r7coeff) at (2.2,-1.5) {$4$}; 
    \node (r8coeff) at (2.7,-1.8) {$2$}; 
    \draw[dashed] (-1,-1.6) -- (1,-1.6); 
  \end{tikzpicture} 
  & $10$ 
  \\ \hline 
  
  4 & 
  $C(3,3)$ & 
  \begin{tikzpicture} 
    \draw[very thick] (0,-1) -- (0,1);  
    \draw[domain=-1:1] plot(\x,{pow(\x,2)}); 
    \draw[domain=-1:1] plot(\x,-{pow(\x,2)});
    \node (r1selfint) at (1,1.2) {$\Circled{-3}$}; 
    \node (r2selfint) at (1,-1.2) {$\Circled{-3}$}; 
    \draw[dashed] (-1,0)--(1,0); 
  \end{tikzpicture}
  & $4$  
  \\ \hline
  
  6 & 
  $C(11,5)$ & 
  \begin{tikzpicture} 
    \draw (-2,0) -- (4,0);  
    \draw (-1,0.2) -- (-1,-1.2);
    \draw[very thick] (-0.8,-1) -- (-2.2,-1); 
    \draw (1,0.2) -- (1,-1.2); 
    \draw (0.8,-1) -- (2.2,-1);
    \draw (2,-0.8) -- (2,-2.2);  
    \draw (1.8,-2) -- (3.2,-2);
    \draw (3,0.2) -- (3,-1.2); 
    \draw (2.8,-1) -- (4.2,-1);
    \draw (4,-0.8) -- (4,-2.2);  
    \draw (3.8,-2) -- (5.2,-2); 
    \node (r1selfint) at (-1,-1.5) {$\Circled{-3}$}; 
    \node (r1coeff) at (-0.8,-0.5) {$2$}; 
    \node (r2coeff) at (-1.8,0.2) {$5$};
    \node (r3coeff) at (1.2,-0.5) {$4$}; 
    \node (r4coeff) at (1.6,-0.8) {$3$}; 
    \node (r5coeff) at (2.2,-1.5) {$2$}; 
    \node (r6coeff) at (2.7,-1.8) {}; 
    \node (r7coeff) at (3.2,-0.5) {$4$}; 
    \node (r8coeff) at (3.6,-0.8) {$3$}; 
    \node (r9coeff) at (4.2,-1.5) {$2$}; 
    \node (r10coeff) at (4.7,-1.8) {};
    \draw[dashed] (0,0.2) -- (0,-1.2); 
  \end{tikzpicture}
  & $12$  
  \\  \hline

  7 & 
  $C(4,6)$ &  
    \begin{tikzpicture} 
      \draw (-2,0) -- (1,0);   
      \draw[domain=-1:1] plot(\x,{pow(\x,2)}); 
      \draw (-1.5,0.2) -- (-1.5,-1.2); 
      \draw[very thick] (-1.3,-1) -- (-2.7,-1); 
      \node (r3selfint) at (1,1.2) {$\Circled{-3}$}; 
      \node (r1coeff) at (-1.3,-0.5) {$2$}; 
      \node (r2coeff) at (-1.8,0.2) {$3$}; 
      \node (r3coeff) at (-0.7,0.9) {$2$}; 
      \draw[dashed] (0,1) -- (0,-1); 
    \end{tikzpicture}
  & $5$  
    \\ \hline 

  8 & 
  $C(9,7)$ &  
   \begin{tikzpicture} 
    \draw (-2,0) -- (4,0);  
    \draw (-1,0.2) -- (-1,-1.2);
    \draw (-0.8,-1) -- (-2.2,-1); 
    \draw (-2,-0.8) -- (-2,-2.2); 
    \draw[very thick] (-1.8,-2) -- (-3.2,-2); 
    \draw (1,0.2) -- (1,-1.2); 
    \draw (0.8,-1) -- (2.2,-1);
    \draw (3.5,0.2) -- (3.5,-1.2); 
    \draw (3.3,-1) -- (4.7,-1);
    \node (r6selfint) at (2.5,-1) {$\Circled{-3}$}; 
    \node (r8selfint) at (5,-1) {$\Circled{-3}$}; 
    \node (r1coeff) at (-1.8,-1.5) {$2$}; 
    \node (r2coeff) at (-1.5,-0.8) {$3$};
    \node (r3coeff) at (-0.8,-0.5) {$4$}; 
    \node (r4coeff) at (-1.8,0.2) {$5$}; 
    \node (r5coeff) at (1.2,-0.5) {$3$};  
    \node (r7coeff) at (3.7,-0.5) {$3$}; 
    \draw[dashed] (0,0.2) -- (0,-1.2); 
    \end{tikzpicture} 
  & $10$ 
    \\ \hline

  9 & 
  $C(13,8)$ &  
    \begin{tikzpicture} 
      \draw (-2,0) -- (4,0);  
      \draw (-1,0.2) -- (-1,-1.2);
      \draw (-0.8,-1) -- (-2.2,-1); 
      \draw (-2,-0.8) -- (-2,-2.2); 
      \draw (-1.8,-2) -- (-3.2,-2); 
      \draw (-3,-1.8) -- (-3,-3.2); 
      \draw (-2.8,-3) -- (-4.2,-3); 
      \draw (-4,-2.8) -- (-4,-4.2); 
      \draw (-3.8,-4) -- (-5.2,-4); 
      \draw[very thick] (-5,-3.8) -- (-5,-4.8);
      \draw (1,0.2) -- (1,-1.8);
      \draw (2,0.2) -- (2,-1.2); 
      \draw (1.8,-1) -- (3.2,-1);
      \node (r12selfint) at (3.5,-1) {$\Circled{-3}$}; 
      \node (r1coeff) at (-4.5,-3.8) {$2$}; 
      \node (r2coeff) at (-3.8,-3.5) {$3$}; 
      \node (r3coeff) at (-3.5,-2.8) {$4$};
      \node (r4coeff) at (-2.8,-2.5) {$5$}; 
      \node (r5coeff) at (-2.5,-1.8) {$6$}; 
      \node (r6coeff) at (-1.8,-1.5) {$7$}; 
      \node (r7coeff) at (-1.5,-0.8) {$8$};
      \node (r8coeff) at (-0.8,-0.5) {$9$}; 
      \node (r9coeff) at (-1.8,0.2) {$10$}; 
      \node (r10coeff) at (1.2,-0.5) {$5$};
      \node (r11coeff) at (2.2,-0.5) {$6$}; 
      \node (r12coeff) at (2.7,-0.8) {$2$};
      \draw[dashed] (0,-1.6) -- (2,-1.6);
    \end{tikzpicture} 
  & $14$ 
    \\ \hline 
  \end{tabular}
  }
  \caption{Configurations of non-multiple fibers} \label{table:configurations of fibers}   
\end{table} 

\noindent In Table 1, 
a solid line stands for an irreducible component of a fiber  
while a dashed line stands for a moving cusp, which is not an irreducible component of the fiber.  
All of the irreducible components are rational curves of genus $0$.  
A circled number beside an irreducible component represents 
its self-intersection number, but we omit $-2$ for $(-2)$-curves. 
We note that the self-intersection number of 
the moving cusp is determined by the global 
combination of reducible fibers 
(cf.\ Corollary \ref{selfint of the moving cusp}), 
so it is not necessary to be equal to $-2$ 
in contrast to an (quasi-)elliptic fibration.  
A non-circled number stands for the multiplicity of each irreducible component,  
but we do not write $1$ for a component with multiplicity $1$. 
Each thick component is the identity component, that is, 
the unique component meeting the zero section. 
Note that a fiber of type $C(1,0)$ is a general fiber. 

\begin{rem.}
  In \cite{NU}, Namikawa and Ueno classified the type of singular fibers 
  for genus $2$ fibration in characteristic $0$. 
  Here is the table for the corresponding notations (Table \ref{table:Comparison table}).   
  
  \begin{table}[H]  
    \centering 
    \scalebox{0.93}{
    \begin{tabular}{|c|c|c|c|c|c|c|c|c|} 
    \hline 
    type & $C(1,0)$ & $C(5,1)$ & $C(9,2)$ & $C(3,3)$ & 
    $C(11,5)$ & $C(4,6)$ & $C(9,7)$ & $C(13,8)$
    \\ \hline   
    \cite{NU} &  VIII-1 & IX-1 & VIII-2 & IX-2 & IX-3 & VIII-3 & IX-4& VIII-4  
    \\ \hline  
    \end{tabular}  
    }
    \caption{Corresponding notations} \label{table:Comparison table} 
  \end{table}
\end{rem.}

The second main theorem gives a formula for the arithmetic genus $p_{a}$ and $(K_{X}^{2})$ by the data of $\varphi(t)$, where $K_{X}$ is the canonical divisor of a quasi-hyperelliptic surface $X$. 

\begin{th.} \label{main.th.2} 
  Let $f : X \to \P^{1}$ be a relatively minimal unirational quasi-hyperelliptic fibration 
  of genus $2$ in characteristic $5$, which is birationally equivalent to 
  an affine hypersurface defined by $y^{2} = x^{5} + \varphi(t)$  
  satisfying Condition \ref{assumption}. 
  Put $d = \deg_{t}\varphi(t)$ and $m = [d/10]$, 
  where $[d/10]$ is the greatest integer less than or equal to $d/10$. 
  Then, we have    
  \begin{align*} 
    p_{a}(X) &= 4m+2 -(N_{3}+N_{5}+2N_{6}+2N_{7}+2N_{8}), \\ 
    (K_{X}^{2}) &= 8m-(2N_{3}+2N_{5}+3N_{6}+4N_{7}+4N_{8}),  
  \end{align*}  
  where $N_{i}$ is the number of reducible fibers of $f$ of type $C(-,i)$.    
\end{th.}

In \cite[Theorem III.3.1]{M78} 
(or \cite[Theorem II.3.9]{MI21}), 
Miyanishi computed $(K_{X}^{2})$ under the assumption that 
fibrations admit reducible fibers of type $C(5,1)$ or $C(9,2)$ only.  
Our theorem is the natural generalization of Miyanishi's result. 

As an application, we classify rational quasi-hyperelliptic surfaces in characteristic $5$ 
(see Theorem \ref{main.th.3}).  

\section{Setting for the proofs} \label{sec:On double covering} 

In this section, we set a stage toward proofs of Theorem 
\ref{main.th.1} and Theorem \ref{main.th.2},   
and introduce some important lemmata on a double covering.  
First, consider an affine hypersurface defined 
by $y^{2} =x^{5} + \varphi(t)$ with 
$\varphi(t) \in k[t]$. 
We assume that $\varphi(t)$ satisfies Condition \ref{assumption}.  
We denote by $\mathscr{K}$ the function field of 
the above affine hypersurface.  
Clearly, this affine hypersurface is a double covering 
of the affine plane $\A_{tx}^{2} \coloneqq \Spec(k[t,x])$ 
with the branch locus $\{ (t,x) \mid x^{5}+\varphi(t)=0 \}$.

Embed the affine plane $\A_{tx}^{2}$ into 
the Hirzebruch surface 
$\Sigma_{0} \coloneqq \P_{k}^{1} \times \P_{k}^{1}$ naturally.    
We denote by $L_{\infty}$ the line defined $t=\infty$ and 
by $M_{\infty}$ the line defined by $x=\infty$. 
Define a curve $C$ in $\Sigma_{0}$ by 
\[
  C \coloneqq \{ (t,x) \in \A_{tx}^{2} \mid x^{5}+\varphi(t)=0 \} 
  \cup \{ (t=\infty,x=\infty) \}.  
\] 

Putting $\tau =1/t$ and $\xi = 1/x$, the defining equation of $C$ at 
$(t=\infty,x=\infty)$ is given by 
\[ 
  \tau^{d}+\xi^{5}\psi(\tau)=0, 
\]
where $\psi(\tau)=\tau^{d}\varphi(1/{\tau})$ with $\psi(0) \neq 0$. 

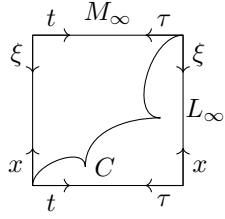
\begin{figure}[H] 
  \centering
  \begin{tikzpicture}
    \draw (0,0) -- (2,0) -- (2,2) -- (0,2) -- (0,0);
    \draw[->] (0,0) to node[below] {$t$} (0.5,0);
    \draw[->] (0,0) to node[left] {$x$} (0,0.5);
    \draw[->] (2,0) to node[below] {$\tau$} (1.5,0);
    \draw[->] (2,0) to node[right] {$x$} (2,0.5);
    \draw[->] (0,2) to node[above] {$t$} (0.5,2);
    \draw[->] (0,2) to node[left] {$\xi$} (0,1.5);
    \draw[->] (2,2) to node[above] {$\tau$} (1.5,2);
    \draw[->] (2,2) to node[right] {$\xi$} (2,1.5); 
    \node (L) at (2.3,1) {$L_{\infty}$}; 
    \node (M) at (1,2.3) {$M_{\infty}$}; 
    \draw (0,0) 
    to [out=90,in=90] (0.7,0.25)node[right] {$C$} 
    to [out=90,in=180] (1.7,0.9)
    to [out=180,in=180] (2,2); 
  \end{tikzpicture}
\end{figure}

Let $S_{0}$ be the normalization of $\Sigma_{0}$ in $\mathscr{K}$ and 
$\rho_{0} : S_{0} \to \Sigma_{0}$ be 
the corresponding normalization morphism. 
Regarding a polynomial $x^{5}+\varphi(t)$ 
as an element of the function field $k(\Sigma_{0})$ of $\Sigma_{0}$, 
we define the principal divisor $A$ on $\Sigma_{0}$ as 
\[ A \coloneqq (x^{5}+\varphi(t)). \] 
Then, $\rho_{0}$ is a double covering of $\Sigma_{0}$ 
with the branch locus $A$.  
Note that we can write the divisor $A$ as 
\[
  A = C - 5M_{\infty} - dL_{\infty}   
\]
because 
$x^{5}+\varphi(t) = (\tau^{d}+\xi^{5}\psi(\tau))\xi^{-5}\tau^{-d}$.

Let $\overline{\sigma} : \overline{\Sigma} \to \Sigma_{0}$ be the shortest 
succession of blowing-ups of $\Sigma_{0}$ at the singular points 
of $C$ and their infinitely near points such that 
the proper transform $\overline{C} \coloneqq \overline{\sigma}'(C)$ 
becomes nonsingular. 
Let $\overline{L_{\infty}}$ and $\overline{M_{\infty}}$ be 
the proper transforms of $L_{\infty}$ and $M_{\infty}$ 
by $\overline{\sigma}$.  
Then, we write the divisor $\overline{\sigma}^{*}A$ as 
\[
  \overline{\sigma}^{*}A = \overline{B} - 2Z,  
\]
where $\overline{B}$ is a reduced divisor, that is, coefficients of 
all irreducible components of $\overline{B}$ are $1$. 
If the set of the irreducible components of $\overline{B}$ and 
that of the connected components of $\overline{B}$ do not coincide, 
we need to blow-up more.  
Let $\sigma : \Sigma \to \overline{\Sigma}$ be the shortest 
succession of blowing-ups of $\overline{\Sigma}$ at the singular points 
of $\overline{B}$ and their infinitely near points 
such that the set of the irreducible components of $B$ and 
that of the connected components of $B$ coincide when we write the divisor 
$(\sigma \circ \overline{\sigma})^{*}A$ as the above form  
\[ 
  (\sigma \circ \overline{\sigma})^{*}A = B -2Z.  
\] 
If the irreducible components of $\overline{B}$ and 
the connected components of $\overline{B}$ coincide on $\overline{\Sigma}$,  
we regard $\sigma$ as the identity map $\mathrm{id}_{\overline{\Sigma}}$.

Let $\overline{S}$ (resp.\ $S$) be the normalization of 
$\overline{\Sigma}$ (resp.\ $\Sigma$) in $\mathscr{K}$ 
and $\overline{\rho}$ (resp.\ $\rho$) be the corresponding 
normalization morphism. 
We define a morphism $\pi : S \to \overline{\Sigma}$ as the composite 
$\pi \coloneqq \sigma \circ \rho$.    
Then, it is known that $S$ is nonsingular and the following 
diagram is commutative (cf.\ \cite{Artin60}, \cite{M78}):  

\begin{center} 
  \begin{tikzpicture}
    \node (a) at (0,0) {$\Sigma_{0}$}; 
    \node (b) at (2,0) {$\overline{\Sigma}$}; 
    \node (c) at (4,0) {$\Sigma$}; 
    \node (d) at (0,2) {$S_{0}$}; 
    \node (e) at (2,2) {$\overline{S}$}; 
    \node (f) at (4,2) {$S$}; 
    \draw[->] (b) to node[above] {$\overline{\sigma}$} (a); 
    \draw[->] (c) to node[above] {$\sigma$} (b); 
    \draw[->] (e) to (d); 
    \draw[->] (f) to (e);  
    \draw[->] (d) to node[left] {$\rho_{0}$} (a); 
    \draw[->] (e) to node[left] {$\overline{\rho}$} (b); 
    \draw[->] (f) to node[left] {$\rho$} (c); 
  \end{tikzpicture}
\end{center} 

Furthermore, we can compute some numerical invariants of $S$ 
by using infomation of $\overline{\Sigma}$ and $\overline{Z}$.    

\begin{lem.}[{\cite[Chapter II]{Artin60}}]  \label{lem:Artin} 
  We continue to use the above notations.    
  Assume that $\overline{B}$ has only negligible singularities, 
  that is, by definition, every point of $\overline{B}$ is one of the following types: 
  \begin{enumerate}[label=(\alph*)]
    \item a simple point of $\overline{B}$, 
    \item a double point of $\overline{B}$, 
    \item a triple point of $\overline{B}$ 
    with at most double point (not necessarily ordinary) infinitely near.  
  \end{enumerate}   
  Then, the following assertions hold.   
  \begin{enumerate}[label=(\roman*)]
    \item $K_{S}$ is linearly equivalent to $\pi^{*}(K_{\overline{\Sigma}}+\overline{Z})$. 
    \item $p_{a}(S) = 2p_{a}(\Sigma_{0})+p_{a}(\overline{Z}) = p_{a}(\overline{Z})$. 
    \item $(K_{S}^{2}) = 2((K_{\overline{\Sigma}}+\overline{Z})^{2})$. 
  \end{enumerate}  
\end{lem.} 

\begin{rem.} 
In our setting, the branch locus $\overline{B}$ has only negligible singularities, 
so the assumption is satisfied and we can use Artin's result.   
\end{rem.}

\begin{lem.} 
Let $p_{1} : \Sigma_{0} \to \P^{1}$ be the first projection, i.e., $p_{1}(t,x)=t$. 
Then, the morphism 
\[ 
  f \coloneqq p_{1} \circ \overline{\sigma} \circ \sigma \circ \rho : S \to \P^{1} 
\] 
defines a quasi-hyperelliptic fibration. 
\end{lem.} 

\begin{proof}
  We denote by $l_{\alpha}$ the line on $\Sigma_{0}$ defined by $t=\alpha$. 
  Let $L = \{ l_{\alpha }\}$ be a linear system on $\Sigma_{0}$ and 
  $l$ be a general member of $L$. 
  Then, $\rho^{-1}((\sigma \circ \overline{\sigma})^{*}(l))$ has 
  the unique singular point $(x=-\varphi(\alpha)^{1/5},y=0)$, 
  so $f$ is a quasi-hyperelliptic fibration.   
\end{proof}  

The following lemma is useful 
to figure out the configurations and self-intersection numbers 
after taking a double covering. 

\begin{lem.}[{\cite[Lemma III.2.7.1]{M78}}] \label{double covering}
  With the above notations, the following assertions hold. 
  \begin{enumerate}[label=(\roman*)]
    \item For divisors $D_{1}, D_{2}$ on $\Sigma$,  
    we have 
    $(\rho^{*}(D_{1}) \cdot \rho^{*}(D_{2})) = 2(D_{1} \cdot D_{2})$. 
    \item For an irreducible component $D$ of $B$, 
    we have $\rho^{*}(D) = 2\Delta$, where $\Delta$ is a nonsingular curve.  
    Furthermore, if $D \cong \P^{1}$, $\Delta \cong \P^{1}$.  
    \item Let $D$ be a curve on $\Sigma$ 
    such that $D \cong \P^{1}$ and $D \nsubseteq \Supp(B)$.   
    \begin{enumerate}
      \item If $D \cap \Supp(B) = \emptyset$, 
      then $\rho^{*}(D) = D_{1} + D_{2}$, $D_{i} \cong \P^{1}$, 
      $D_{1} \cap D_{2} = \emptyset$, and $(D_{i}^{2}) = (D^{2})$  
      for $i=1,2$.  
      \item If $D$ intersects exactly two components of $B$ transversally, 
      say $B_{1}, B_{2}$, and $D \cap B_{1} \neq D \cap B_{2}$,  
      then $\rho^{*}(D)$ is irreducible and $\rho^{*}(D) \cong \P^{1}$.  
    \end{enumerate}
  \end{enumerate} 
\end{lem.}

\section{Proof of Theorem \ref{main.th.1}}  \label{sec:pf.of.main.th.1}

We use the same notations as in the previous section. 
This section is not only the proof of Theorem \ref{main.th.1} but also the preparation for the proof of Theorem \ref{main.th.2},  
so, put $P = (t=0,x=0)$ and consider the equation $x^{5}+t^{e_{P}}=0$ 
for $e_{P} \in \{ 2,3,4,5,6,7,8,9 \}$ due to Condition  \ref{assumption}. 
Furthermore, in local calculations, we can also assume that $e_{P} \neq 5$. 
Indeed, if $e_{\alpha} = 5$ and $\varphi(t) = (t-\alpha)^{e_{\alpha}}\varphi_{\alpha}(t)$ 
with $\varphi_{\alpha}(\alpha) \neq 0$,    
a birational transformation of the type 
\[ 
(t,x,y) \mapsto (t,x+\gamma(t-\alpha),y),
\]
where $\gamma^{5} = \varphi_{\alpha}(\alpha)$, enables us 
to assume $e_{\alpha} \ge 6$.

We introduce some notations here for later use.     
\begin{itemize} 
\item Let $E_{P}$ be the contribution to the effective divisor 
$\overline{\sigma}^{*}C - \overline{C}$ and we write  
$E_{P} = E_{P}^{(1)}+2E_{P}^{(2)}$, 
where $E_{P}^{(i)} \ge 0$ for $i=1,2$ and $E_{P}^{(1)}$ is a reduced divisor
which is possibly to be zero.   
\item Let $E_{\mathrm{can},P}$ be the contribution to the effective divisor 
$K_{\overline{\Sigma}} - \overline{\sigma}^{*}K_{\Sigma_{0}}$.  
\item  $\mu_{P} \coloneqq 
\frac{1}{2}(E_{P}^{(2)} 
\cdot E_{P}^{(2)}-E_{\mathrm{can},P})$.  
\item $\nu_{P} \coloneqq ((E_{P}^{(2)}-E_{\mathrm{can},P})^{2})$.  
\end{itemize} 

\subsubsection*{Case \texorpdfstring{$e_{P}=2$}{e=2}} 

We shall start the equation $x^{5}+t^{2}=0$.   
After twice blowing-ups, 
we have the configuration on $\overline{\Sigma}$ 
(Figure \ref{sec:pf.of.main.th.1}-1).     

\begin{figure}[H]  
  \centering 
  \begin{tikzpicture}
    \draw[domain=-1:1] plot(\x,{pow(\x,2)})node[right]{$\overline{C}$}; 
    \draw[dashed] (-1,0)node[below]{$E_{2}$}--(2,0); 
    \draw[dashed] (1.5,0.5)-- (1.5,-1)node[left]{$E_{1}$}; 
    \draw[dashed] (0,1)--(0,-0.5)node[right]{$E_{0}$}; 
    \node (a) at (3,1) {$(\overline{C} \cdot E_{2})=2$}; 
    \node (E0selfint) at (0,1.3) {$\Circled{-2}$}; 
    \node (E1selfint) at (1.5,-1.5) {$\Circled{-2}$}; 
    \node (E2selfint) at (-1.5,0) {$\Circled{-1}$}; 
  \end{tikzpicture} 
  \caption*{Figure \ref{sec:pf.of.main.th.1}-1} 
\end{figure} 

\noindent In Figure \ref{sec:pf.of.main.th.1}-1, the rational curve $E_{0}$ is not an exceptional curve, 
but the proper transform of $x$-axis on $\Sigma_{0}$.  
The rational curves $E_{1}$ and $E_{2}$ are exceptional curves and 
indexed numbers represent the order of appearance.  
Circled numbers close to lines are the self-intersection numbers. 

Then, we have $E_{P} = 2E_{1}+4E_{2}$, 
$E_{P}^{(1)}=0$, $E_{P}^{(2)}=E_{1}+2E_{2}$, 
$E_{\mathrm{can},P}=E_{1}+2E_{2}$, 
$\mu_{P}=0$ and $\nu_{P}=0$. 
The solid lines (including $\overline{C}$) are contained in $\overline{B}$, while 
the broken lines are not contained in $\overline{B}$.

As one can see,  
we do not need to blow-up more, 
so we take a double cover (cf.\ Section \ref{sec:On double covering}). 
By Lemma \ref{double covering}, we get the configuration on $S$ as in Figure \ref{sec:pf.of.main.th.1}-2, 
which is $R_{0}+2(R_{1}+R_{3})+(R_{2}+R_{4})$ as a divisor. 

\begin{figure}[H]
  \centering 
  \begin{tikzpicture} 
    \draw (-1.5,0) -- (1.5,0)node[right]{$\widetilde{E_{0}} = R_{0}$};
    \draw (-1,-1)node[below]{$\widetilde{E_{2}'} = R_{3}$} -- (1,1); 
    \draw (-1,1)node[above]{$\widetilde{E_{2}} = R_{1}$} -- (1,-1); 
    \draw (0.5,0.7) -- (2.5,0.7)node[above]{$\widetilde{E_{1}} = R_{2}$};  
    \draw (0.5,-0.7) -- (2.5,-0.7)node[below]{$\widetilde{E_{1}'} = R_{4}$};  
    \node (r0selfint) at (-1.8,0) {$\Circled{-4}$}; 
    \draw[dashed] (0,1) -- (0,-1)node[below]{$\Gamma$}; 
  \end{tikzpicture} 
  \caption*{Figure \ref{sec:pf.of.main.th.1}-2}  
\end{figure}

\noindent In Figure \ref{sec:pf.of.main.th.1}-2, 
The dashed line $\Gamma$ is the moving cusp and  
a circled number represents the self-intersection number 
of each curve except for $\Gamma$.  
We omit the self-intersection numbers for $(-2)$-curves to simplify the figure. Note that the self-intersection number $(\Gamma^{2})$ cannot be determined locally, so it is unknown at present. 
We will compute it 
from the global data of a fibration 
(see Corollary \ref{selfint of the moving cusp}). 
The component $R_{0} = \widetilde{E_{0}}$ is the unique component meeting the zero section.

For other $e_{P}$, we will describe configurations as similar fashion.  
Therefore, we omit detailed calculations.    

\subsubsection*{Case \texorpdfstring{$e_{P}=3$}{e=3}}

Starting the equation $x^{5}+t^{3}=0$, we have 
the configuration on $\overline{\Sigma}$ as in Figure \ref{sec:pf.of.main.th.1}-3.    
\begin{figure}[H] 
  \centering 
  \begin{tikzpicture}
    \draw[domain=-1:1] plot(\x,{pow(\x,2)})node[right]{$\overline{C}$}; 
    \draw (-2,0)--(1,0)node[below]{$E_{2}$}; 
    \draw (0,1)-- (0,-0.5)node[left]{$E_{1}$}; 
    \draw[dashed] (-1.5,0.5)--(-1.5,-1)node[left]{$E_{0}$}; 
    \node (E0selfint) at (-1.5,0.8) {$\Circled{-2}$}; 
    \node (E1selfint) at (0,1.3) {$\Circled{-2}$}; 
    \node (E2selfint) at (1.5,0) {$\Circled{-1}$}; 
    \node (mult) at (3,1) {$(\overline{C} \cdot E_{2})=2$};
  \end{tikzpicture} 
  \caption*{Figure \ref{sec:pf.of.main.th.1}-3} 
\end{figure} 
\noindent Then, we have $E_{P} = 3E_{1}+5E_{2}$, 
$E_{P}^{(1)}=E_{1}+E_{2}$, $E_{P}^{(2)}=E_{1}+2E_{2}$, 
$E_{\mathrm{can},P}=E_{1}+2E_{2}$, 
$\mu_{P}=0$, and $\nu_{P}=0$.

As one can see, 
we need to blow-up several times at the point $\overline{C} \cap E_{1} \cap E_{2}$.  
Since three irreducible components $\overline{C}$, $E_{1}$ and $E_{2}$ 
are contained in $\overline{B}$, 
appearing exceptional curves from this point 
is also components of $\overline{B}$. 
Thus, we have the configuration on $\Sigma$ as in Figure \ref{sec:pf.of.main.th.1}-4.   

\begin{figure}[H] 
  \centering 
  \begin{tikzpicture} 
    \draw (-2,0) -- (2,0); 
    \draw[dashed] (-1,0.2) -- (-1,-1.2); 
    \draw (-0.8,-1) -- (-2.2,-1)node[left]{$E_{2}$}; 
    \draw[dashed] (-2,-0.8) -- (-2,-2.2)node[left]{$E_{0}$};  
    \draw[dashed] (0,0.2) -- (0,-1.8);   
    \draw[dashed] (1,0.2) -- (1,-1.2); 
    \draw (0.8,-1) -- (2.2,-1);
    \draw[dashed] (2,-0.8) -- (2,-2.2);  
    \draw (1.8,-2) -- (3.2,-2)node[below]{$E_{1}$};
    \node (E1selfint) at (3.7,-2) {$\Circled{-4}$};
    \node (E0selfint) at (-2,-2.7) {$\Circled{-4}$}; 
    \node (E3selfint) at (2.5,-1) {$\Circled{-4}$}; 
    \node (E4selfint) at (2.3,0) {$\Circled{-4}$}; 
    \node (left-tate) at (-1,0.5) {$\Circled{-1}$};
    \node (right-tate) at (1,0.5) {$\Circled{-1}$};
    \node (E2selfint) at (-3.1,-1) {$\Circled{-4}$};
    \node (chuo-tate) at (0,0.5) {$\Circled{-1}$};
    \node (right-tate-shita) at (2,-2.5) {$\Circled{-1}$};
    \draw (-1,-1.6) -- (1,-1.6)node[below]{$\sigma'\overline{C}$};  
  \end{tikzpicture} 
  \caption*{Figure \ref{sec:pf.of.main.th.1}-4} 
\end{figure}

Therefore, we have the configuration on $S$ as in Figure \ref{sec:pf.of.main.th.1}-5, 
which is $R_{0}+4R_{1}+7R_{2}+10R_{3}+5R_{4}+8R_{5}+6R_{6}+4R_{7}+2R_{8}$ as a divisor. 

\begin{figure}[H] 
  \centering 
  \begin{tikzpicture} 
    \draw (-2,0)node[left]{$R_{3}$} -- (2,0); 
    \draw (-1,0.2)node[above]{$R_{2}$} -- (-1,-1.2); 
    \draw (-0.8,-1) -- (-2.2,-1)node[left]{$R_{1} = \widetilde{E_{2}}$};
    \draw (-2,-0.8) -- (-2,-2.2)node[left]{$R_{0} = \widetilde{E_{0}}$};  
    \draw (0,0.2)node[above]{$R_{4}$} -- (0,-1.8);   
    \draw (1,0.2)node[above]{$R_{5}$} -- (1,-1.2); 
    \draw (0.8,-1) -- (2.2,-1)node[right]{$R_{6}$};
    \draw (2,-0.8)node[above]{$R_{7}$} -- (2,-2.2);  
    \draw (1.8,-2) -- (3.2,-2)node[right]{$R_{8} = \widetilde{E_{1}}$};
    \node (r0selfint) at (-1.7,-2.5) {$\Circled{-4}$}; 
    \draw[dashed] (-1,-1.6) -- (1,-1.6)node[below]{$\Gamma$};  
  \end{tikzpicture} 
  \caption*{Figure \ref{sec:pf.of.main.th.1}-5} 
\end{figure}

\subsubsection*{Case \texorpdfstring{$e_{P}=4$}{e=4}}

Starting the equation $x^{5}+t^{4}=0$, we have 
the configuration on $\overline{\Sigma}$ as in Figure \ref{sec:pf.of.main.th.1}-6. 
\begin{figure}[H] 
  \centering 
  \begin{tikzpicture}
    \draw[domain=-1:1] plot(\x,{pow(\x,4)})node[right]{$\overline{C}$}; 
    \draw[dashed] (-1,0)--(1,0)node[above]{$E_{1}$};  
    \draw[dashed] (0,1)--(0,-0.5)node[right]{$E_{0}$};  
    \node (E0selfint) at (0,1.3) {$\Circled{-1}$};  
    \node (E1selfint) at (1.5,0) {$\Circled{-1}$};  
    \node (mult) at (3,1) {$(\overline{C} \cdot E_{1})=4$}; 
  \end{tikzpicture} 
  \caption*{Figure \ref{sec:pf.of.main.th.1}-6} 
\end{figure}

\noindent Then, we have $E_{P} = 4E_{1}$,  
$E_{P}^{(1)}=0$, $E_{P}^{(2)}=2E_{1}$, 
$E_{\mathrm{can},P}=E_{1}$, 
$\mu_{P}=-1$, and $\nu_{P}=-1$.  

As one can see, we do not have to blow-up more. 
Thus, the configuration on $S$ is as in Figure \ref{sec:pf.of.main.th.1}-7, 
which is $R_{0}+R_{1}+R_{2}$ as a divisor. 

\begin{figure}[H]
  \centering
  \begin{tikzpicture} 
    \draw (0,-1)node[below]{$R_{0} = \widetilde{E_{0}}$} -- (0,1);  
    \draw[domain=-1:1] plot(\x,{pow(\x,2)}); 
    \draw[domain=-1:1] plot(\x,-{pow(\x,2)});
    \node (r1selfint) at (1,1.2) {$\Circled{-3}$}; 
    \node (r2selfint) at (1,-1.2) {$\Circled{-3}$}; 
    \draw[dashed] (-1,0)--(1,0)node[below]{$\Gamma$}; 
    \node (r1name) at (-1.2,1.3) {$R_{1}$}; 
    \node (r2name) at (-1.2,-1.3) {$R_{2}$}; 
  \end{tikzpicture} 
  \caption*{Figure \ref{sec:pf.of.main.th.1}-7}  
\end{figure}

Moreover, we give the expression of $\pi^{*}(E_{\mathrm{can},P}-E_{P}^{(2)})$ here 
for later use. 
Since $\pi^{*}(E_{\mathrm{can},P}-E_{P}^{(2)})=-E_{1}$, we have  
\begin{align*}
  \pi^{*}(E_{\mathrm{can},P}-E_{P}^{(2)})=-R_{1}-R_{2}.  
\end{align*}

\subsubsection*{Case \texorpdfstring{$e_{P}=6$}{e=6}}

Starting the equation $x^{5}+t^{6}=0$, we have 
the configuration on $\overline{\Sigma}$ as in Figure \ref{sec:pf.of.main.th.1}-8.   
\begin{figure}[H] 
  \centering 
  \begin{tikzpicture}
    \draw[domain=-1:1] plot(\x,{pow(\x,5)})node[right]{$\overline{C}$}; 
    \draw (-2,0)--(1,0)node[below]{$E_{1}$};   
    \draw[dashed] (-1.5,0.5)--(-1.5,-1)node[left]{$E_{0}$}; 
    \node (E0selfint) at (-1.5,0.8) {$\Circled{-1}$}; 
    \node (E1selfint) at (1.5,0) {$\Circled{-1}$}; 
    \node (mult) at (3,1) {$(\overline{C} \cdot E_{1})=5$};  
  \end{tikzpicture} 
  \caption*{Figure \ref{sec:pf.of.main.th.1}-8} 
\end{figure} 

\noindent Then, we have $E_{P} = 5E_{1}$, 
$E_{P}^{(1)}=E_{1}$, $E_{P}^{(2)}=2E_{1}$, 
$E_{\mathrm{can},P}=E_{1}$, 
$\mu_{P}=-1$, and $\nu_{P}=-1$.    

After some blowing-ups and taking a double cover, 
we have the configuration on $S$ as in Figure \ref{sec:pf.of.main.th.1}-9, 
which is $R_{0}+2R_{1}+5R_{2}+4(R_{3}+R_{7})+(3R_{4}+R_{8})+2(R_{5}+R_{9})+(R_{6}+R_{10})$ as a divisor. 

\begin{figure}[H] 
  \centering
  \begin{tikzpicture} 
    \draw (-2,0)node[left]{$R_{2}$} -- (4,0);  
    \draw (-1,0.2)node[above]{$R_{1} = \widetilde{E_{1}}$} -- (-1,-1.2);
    \draw (-0.8,-1) -- (-2.2,-1)node[left]{$R_{0} = \widetilde{E_{0}}$}; 
    \draw (1,0.2)node[above]{$R_{3}$} -- (1,-1.2); 
    \draw (0.8,-1)node[left]{$R_{4}$} -- (2.2,-1);
    \draw (2,-0.8)node[above]{$R_{5}$} -- (2,-2.2);  
    \draw (1.8,-2) -- (3.2,-2)node[below]{$R_{6}$};
    \draw (3,0.2)node[above]{$R_{7}$} -- (3,-1.2); 
    \draw (2.8,-1) -- (4.2,-1)node[right]{$R_{8}$};
    \draw (4,-0.8)node[above]{$R_{9}$} -- (4,-2.2);  
    \draw (3.8,-2) -- (5.2,-2)node[below]{$R_{10}$}; 
    \draw[dashed] (0,0.2) -- (0,-1.2)node[below]{$\Gamma$}; 
    \node (selfint) at (-1,-1.5) {$\Circled{-3}$}; 
  \end{tikzpicture} 
  \caption*{Figure \ref{sec:pf.of.main.th.1}-9} 
\end{figure} 

Moreover, we give the expression of $\pi^{*}(E_{\mathrm{can},P}-E_{P}^{(2)})$ here 
for later use. 
Since $\pi^{*}(E_{\mathrm{can},P}-E_{P}^{(2)})=-E_{1}$, we have  
\begin{align*}
  &\pi^{*}(E_{\mathrm{can},P}-E_{P}^{(2)}) \\
  &=-2R_{1}-5R_{2}-4(R_{3}+R_{7})-(3R_{4}+R_{8})-2(R_{5}+R_{9})
  -(R_{6}+R_{10}).   
\end{align*}

\subsubsection*{Case \texorpdfstring{$e_{P}=7$}{e=7}}

Starting the equation $x^{5}+t^{7}=0$, we have 
the configuration on $\overline{\Sigma}$ as in Figure \ref{sec:pf.of.main.th.1}-10.   
\begin{figure}[H] 
  \centering 
  \begin{tikzpicture}
    \draw[domain=-1:1] plot(\x,{pow(\x,2)})node[right]{$\overline{C}$}; 
    \draw[dashed] (-1,0)node[left]{$E_{3}$}--(2,0); 
    \draw (1.5,0.5)-- (1.5,-1)node[right]{$E_{2}$}; 
    \draw (0,1)--(0,-1.5)node[right]{$E_{1}$}; 
    \draw[dashed] (0.5,-1)--(-1,-1)node[below]{$E_{0}$}; 
    \node (E0selfint) at (-1.5,-1) {$\Circled{-1}$}; 
    \node (E1selfint) at (0,1.3) {$\Circled{-3}$}; 
    \node (E2selfint) at (1.5,-1.5) {$\Circled{-2}$};
    \node (E3selfint) at (2.5,0) {$\Circled{-1}$};
    \node (mult) at (3,1) {$(\overline{C} \cdot E_{3})=2$}; 
  \end{tikzpicture} 
  \caption*{Figure \ref{sec:pf.of.main.th.1}-10}  
\end{figure}

\noindent Then, we have $E_{P} = 5E_{1}+7E_{2}+14E_{3}$, 
$E_{P}^{(1)}=E_{1}+E_{2}$, $E_{P}^{(2)}=2E_{1}+3E_{2}+7E_{3}$, 
$E_{\mathrm{can},P}=E_{1}+2E_{2}+4E_{3}$, 
$\mu_{P}=-2$, and $\nu_{P}=-2$.   

After blowing-up once, taking a double cover and contracting a $(-1)$-curve contained in the fiber,   
we have the configuration as in Figure \ref{sec:pf.of.main.th.1}-11, which is  $R_{0}+2R_{1}+3R_{2}+2R_{3}$ as a divisor. 

\begin{figure}[H] 
  \centering
  \begin{tikzpicture} 
    \draw (-2,0)node[left]{$R_{2}$} -- (1,0);   
    \draw[domain=-1:1] plot(\x,{pow(\x,2)})node[right]{$R_{3} = \widetilde{E_{3}}$}; 
    \draw (-1.5,0.2) -- (-1.5,-1.2)node[below]{$R_{1} = \widetilde{E_{1}}$}; 
    \draw (-1.3,-1) -- (-2.7,-1)node[left]{$R_{0} = \widetilde{E_{0}}$}; 
    \node (e3selfint) at (0.8,1.4) {$\Circled{-3}$};  
    \node (e2selfint) at (2.5,0.5) {$\Circled{-1}$}; 
    \draw[dashed] (0,1) -- (0,-1)node[right]{$\Gamma$}; 
  \end{tikzpicture} 
  \caption*{Figure \ref{sec:pf.of.main.th.1}-11} 
\end{figure} 

Moreover, we give the expression of $\pi^{*}(E_{\mathrm{can},P}-E_{P}^{(2)})$ here 
for later use. 
Since $\pi^{*}(E_{\mathrm{can},P}-E_{P}^{(2)})=-E_{1}-E_{2}-3E_{3}$, we have  
\begin{align*}
  &\pi^{*}(E_{\mathrm{can},P}-E_{P}^{(2)}) 
  =-2R_{1}-4R_{2}-3R_{3}.  
\end{align*}

\subsubsection*{Case \texorpdfstring{$e_{P}=8$}{e=8}}

Starting the equation $x^{5}+t^{8}=0$, we have 
the configuration on $\overline{\Sigma}$ as in Figure \ref{sec:pf.of.main.th.1}-12.  
\begin{figure}[H] 
  \centering 
  \begin{tikzpicture}
    \draw[domain=-1:1] plot(\x,{pow(\x,2)})node[right]{$\overline{C}$}; 
    \draw (-2,0)--(1,0)node[above]{$E_{3}$}; 
    \draw (-1.5,0.5)-- (-1.5,-1.2)node[below]{$E_{1}$}; 
    \draw[dashed] (0,1)--(0,-1)node[right]{$E_{2}$}; 
    \draw[dashed] (-1.3,-1)--(-2.5,-1)node[below]{$E_{0}$}; 
    \node (E0selfint) at (-3,-1) {$\Circled{-1}$}; 
    \node (E1selfint) at (0,1.3) {$\Circled{-3}$}; 
    \node (E2selfint) at (-1.5,0.8) {$\Circled{-2}$};
    \node (E3selfint) at (1.5,0) {$\Circled{-1}$};
    \node (mult) at (3,1) {$(\overline{C} \cdot E_{3})=2$}; 
  \end{tikzpicture} 
  \caption*{Figure \ref{sec:pf.of.main.th.1}-12}  
\end{figure} 

\noindent Then, we have $E_{P} = 5E_{1}+8E_{2}+15E_{3}$,  
$E_{P}^{(1)}=E_{1}+E_{3}$, $E_{P}^{(2)}=2E_{1}+4E_{2}+7E_{3}$, 
$E_{\mathrm{can},P}=E_{1}+2E_{2}+4E_{3}$, 
$\mu_{P}=-2$, and $\nu_{P}=-2$.   

After blowing-ups several times and taking a double cover, 
we have the configuration on $S$ as in Figure \ref{sec:pf.of.main.th.1}-13, which is 
$R_{0}+2R_{1}+3R_{2}+4R_{3}+5R_{4}+3(R_{5}+R_{7})+(R_{6}+R_{8})$ as a divisor. 

\begin{figure}[H] 
  \centering
  \begin{tikzpicture} 
    \draw (-2,0) -- (4,0);  
    \draw (-1,0.2) -- (-1,-1.2)node[below]{$R_{3} = \widetilde{E_{3}}$}; 
    \draw (-0.8,-1) -- (-2.2,-1)node[left]{$R_{2}$}; 
    \draw (-2,-0.8) -- (-2,-2.2)node[below]{$R_{1} = \widetilde{E_{1}}$};  
    \draw (-1.8,-2) -- (-3.2,-2)node[left]{$R_{0} = \widetilde{E_{0}}$};  
    \draw (1,0.2) -- (1,-1.2)node[below]{$R_{5}$}; 
    \draw (0.8,-1) -- (2.2,-1)node[below=3.5pt]{$R_{6} = \widetilde{E_{2}}$};
    \draw (3.5,0.2) -- (3.5,-1.2)node[below]{$R_{7}$}; 
    \draw (3.3,-1) -- (4.7,-1)node[below=3.5pt]{$R_{8} = \widetilde{E_{2}'}$}; 
    \node (r6selfint) at (2.7,-1) {$\Circled{-3}$}; 
    \node (r8selfint) at (5.2,-1) {$\Circled{-3}$}; 
    \draw[dashed] (0,0.2) -- (0,-1.2)node[right]{$\Gamma$}; 
    \end{tikzpicture} 
    \caption*{Figure \ref{sec:pf.of.main.th.1}-13}  
\end{figure} 

Moreover, we give the expression of $\pi^{*}(E_{\mathrm{can},P}-E_{P}^{(2)})$ here 
for later use. 
Since $E_{\mathrm{can},P}-E_{P}^{(2)}=-E_{1}-2E_{2}-3E_{3}$, we have 
\begin{align*}
  &\pi^{*}(E_{\mathrm{can},P}-E_{P}^{(2)}) \\
  &=-2R_{1}-4R_{2}-6R_{3}-8R_{4}-5(R_{5}+R_{7})-2(R_{6}+R_{8}).   
\end{align*}

\subsubsection*{Case \texorpdfstring{$e_{P}=9$}{e=9}} 

Starting the equation $x^{5}+t^{9}=0$, we have 
the configuration on $\overline{\Sigma}$ as in Figure \ref{sec:pf.of.main.th.1}-14.  
\begin{figure}[H] 
  \centering 
  \begin{tikzpicture}
    \draw[domain=-1:1] plot(\x,{pow(\x,2)})node[right]{$\overline{C}$}; 
    \draw (-1,0)--(1,0)node[below]{$E_{2}$}; 
    \draw (0,1)-- (0,-1.5)node[left]{$E_{1}$}; 
    \draw[dashed] (-1,-1)--(2,-1)node[below]{$E_{0}$}; 
    \node (E0selfint) at (2.5,-1) {$\Circled{-1}$}; 
    \node (E1selfint) at (0,1.3) {$\Circled{-2}$}; 
    \node (E2selfint) at (1.5,0) {$\Circled{-1}$}; 
    \node (mult) at (3,1) {$(\overline{C} \cdot E_{2})=4$}; 
  \end{tikzpicture} 
  \caption*{Figure \ref{sec:pf.of.main.th.1}-14} 
\end{figure} 

\noindent Then, we have $E = 5E_{1}+9E_{2}$, 
$E_{P}^{(1)}=E_{1}+E_{2}$, $E_{P}^{(2)}=2E_{1}+4E_{2}$, 
$E_{\mathrm{can},P}=E_{1}+2E_{2}$, 
$\mu_{P}=-2$, $\nu_{P}=-2$. 

After blowing-ups several times and taking a double covering, 
we have the configuration as in Figure \ref{sec:pf.of.main.th.1}-15, which is 
$R_{0}+2R_{1}+3R_{2}+4R_{3}+5R_{4}+6R_{5}+7R_{6}+8R_{7}+9R_{8}+10R_{9}+5R_{10}+6R_{11}+2R_{12}$ 
as a divisor. 

\begin{figure}[H] 
  \centering
  \scalebox{0.8}{
  \begin{tikzpicture} 
    \draw (-2,0)node[left]{$R_{9}$} -- (4,0);  
    \draw (-1,0.2) -- (-1,-1.2)node[below]{$R_{8}$};
    \draw (-0.8,-1) -- (-2.2,-1)node[left]{$R_{7}$}; 
    \draw (-2,-0.8) -- (-2,-2.2)node[below]{$R_{6}$}; 
    \draw (-1.8,-2) -- (-3.2,-2)node[left]{$R_{5}$}; 
    \draw (-3,-1.8) -- (-3,-3.2)node[below]{$R_{4}$}; 
    \draw (-2.8,-3) -- (-4.2,-3)node[left]{$R_{3}$}; 
    \draw (-4,-2.8) -- (-4,-4.2)node[below]{$R_{2}$}; 
    \draw (-3.8,-4) -- (-5.2,-4)node[left]{$R_{1} = \widetilde{E_{1}}$}; 
    \draw (-5,-3.8) -- (-5,-4.8)node[left]{$R_{0} = \widetilde{E_{0}}$}; 
    \draw (1,0.2)node[above]{$R_{10}$} -- (1,-1.8);
    \draw (2,0.2)node[above]{$R_{11}$} -- (2,-1.2); 
    \draw (1.8,-1) -- (3.2,-1)node[below=3.5pt]{$R_{12} = \widetilde{E_{2}}$};
    \node (r12selfint) at (3.7,-1) {$\Circled{-3}$}; 
    \draw[dashed] (0,-1.6) -- (2,-1.6)node[below]{$\Gamma$}; 
  \end{tikzpicture}
  } 
  \caption*{Figure \ref{sec:pf.of.main.th.1}-15}  
\end{figure} 

Moreover, we give the expression of $\pi^{*}(E_{\mathrm{can},P}-E_{P}^{(2)})$ here 
for later use. 
Since $E_{\mathrm{can},P}-E_{P}^{(2)}=-E_{1}-2E_{2}$, we have 
\begin{multline*}
  \pi^{*}(E_{\mathrm{can},P}-E_{P}^{(2)}) 
  =-2R_{1}-2R_{2}-6R_{3}-8R_{4}-10R_{5}-12R_{6} \\
    -14R_{7}-16R_{8}-18R_{9}-9R_{10}-11R_{11}-4R_{12}.
\end{multline*} 

This completes the classification of reducible fibers. 
Since all of the irreducible components of reducible fibers are rational curves of genus $0$,  
it is easy to compute the topological Euler number of each reducible fiber 
from the configuration.   
$\qed$ 

\medskip 

Towards the next section, we put 
\[ 
E_{\mathrm{fin}} \coloneqq \sum_{P}E_{P}, \ 
E_{\mathrm{can}} \coloneqq \sum_{P}E_{\mathrm{can},P}, \   
\mu = \sum_{P}\mu_{P}, \ \nu = \sum_{P}\nu_{P}, 
\]
where $P$ runs through all singular points of $C$ in $\A_{tx}^{2}$.  

For $P,Q \in X_{\eta} \setminus \{ P_{\infty} \}$, 
Table \ref{table:contribution} represents local contributions 
of each reducible fiber for the height pairing (cf.\ Section \ref{sec:MWG}).    

\begin{table}[H] 
  \centering 
  \scalebox{0.9}{
  \begin{tabular}{|c|c|c|c|c|c|c|c|} 
  \hline 
  type & 
  $C(5,1)$ & $C(9,2)$ & $C(3,3)$ & 
  $C(11,5)$ & $C(4,6)$ & $C(9,7)$ & 
  $C(13,8)$ 
  \\ \hline   
  $\mathrm{contr}_{v}(P)$ & $\frac{4}{5}$ & - & $\frac{3}{5}$ & $\frac{7}{5}$ & - & $\frac{6}{5}$ & -  
  \\ \hline 
  $\mathrm{contr}_{v}(P,Q)$ & $\frac{1}{5}$ & - & $\frac{2}{5}$ & $\frac{3}{5}$ & - & $\frac{4}{5}$ & -
  \\ \hline
  \end{tabular}     
  }
  \caption{Local contributions for the height pairing} \label{table:contribution}  
\end{table} 

\section{Proof of Theorem \ref{main.th.2}}  \label{sec:pf.of.main.th.2}

In this section, we prove Theorem \ref{main.th.2}. 
In the same setting as in Section \ref{sec:On double covering}, 
we look at the singular point $P_{\infty} = (t=\infty, x=\infty)$ of the curve at the infinity. 
At the point $P_{\infty}$, the curve $C$ is defined by 
\[
  \tau^{d}+\xi^{5}\psi(\tau) = 0, 
\]
where $\psi(\tau) = \tau^{d}\varphi(1/\tau)$ with $\psi(0) \neq 0$. 
First, we divide into four cases according to the degree $d$, 
case $d=5n+1$, case $d=5n+2$, case $5n+3$ and case $5n+4$ with $n \in \Z_{>0}$.    
Then, we divide into two subcases for each case.  
For cases $n=0$, we need to modify indices in the proof a little bit. 
However, we omit the detail because most of computations are same as 
in the case $n>0$ (cf.\ \cite{M78}). 

\subsubsection*{Case \texorpdfstring{$1$}{1} \texorpdfstring{$d=5n+1$ $(n>0)$}{d=5n+1}}  

We have to blow-up $n$ times to obatain the nonsingular curve $\overline{C}$ (Figure \ref{sec:pf.of.main.th.2}-1).  

\begin{figure}[H]
  \centering
  \begin{tikzpicture}
    \draw (0.5,1) to [out=225,in=135] (0.25,0.25) 
    to [out=315,in=225] (1.5,0.5)node[below]{$\overline{C}$}; 
    \draw (1,1)node[right]{$\overline{M_{\infty}}$}--(0,0); 
    \draw (0.5,0)node[below]{$E_{n}$}--(-0.5,1); 
    \draw (0,1)--(-1,0)node[below]{$E_{n-1}$};   
    \node (tenten) at (-1.5,0.5) {$\cdots$}; 
    \draw (-2,1)--(-3,0)node[below]{$E_{2}$};  
    \draw (-2.5,0)--(-3.5,1)node[above]{$E_{1}$}; 
    \draw (-3,1)--(-4,0)node[below]{$\overline{L_{\infty}}$};  
    \node (Lselfint) at (-4,-0.8) {$\Circled{-1}$}; 
    \node (E1selfint) at (-3.5,1.7) {$\Circled{-2}$}; 
    \node (E2selfint) at (-3,-0.7) {$\Circled{-2}$}; 
    \node (En-1selfint) at (-1,-0.7) {$\Circled{-2}$}; 
    \node (Enselfint) at (0.5,-0.7) {$\Circled{-1}$}; 
    \node (Mselfint) at (1.5,1.5) {$\Circled{-n}$}; 
    \node (intersection) at (3,0) {$(\overline{C} \cdot E_{n})=5$};
  \end{tikzpicture} 
  \caption*{Figure \ref{sec:pf.of.main.th.2}-1}  
\end{figure}
Then, we have 
\begin{align*} 
  \overline{\sigma}^{*}A 
  &= \overline{\sigma}^{*}C - 5\overline{\sigma}^{*}M_{\infty} 
  -d\overline{\sigma}^{*}L_{\infty} \\
  &= (\overline{C}+5(E_{1}+2E_{2}+ \dots + nE_{n})+E_{\mathrm{fin}}) \\ 
  & -5(\overline{M_{\infty}}+E_{1}+2E_{2}+ \dots +nE_{n})  
  -d(\overline{L_{\infty}}+E_{1}+E_{2}+ \dots +E_{n}) \\ 
  &= \overline{C}-5\overline{M_{\infty}}-d\overline{L_{\infty}}
  -d(E_{1}+E_{2}+ \dots +E_{n}) + E_{\mathrm{fin}},  
\end{align*}
and 
\begin{align*}
  K_{\overline{\Sigma}} 
  &\sim -2\overline{\sigma}^{*}M_{\infty} 
  -2\overline{\sigma}^{*}L_{\infty} +(E_{1}+2E_{2}+ \dots +nE_{n}) 
  +E_{\mathrm{can}} \\ 
  &= -2\overline{M_{\infty}}-2\overline{L_{\infty}}
  -(3E_{1}+4E_{2}+ \dots + (n+2)E_{n}) + E_{\mathrm{can}}, 
\end{align*} 
\noindent where $E_{\mathrm{fin}}$ and $E_{\mathrm{can}}$ are introduced 
in the end of the previous section.  
In order to consider the branch locus, 
we need to separate into two subcases.  

\medskip 

\paragraph*{Subcase 1-1 \texorpdfstring{$d=10m+1$ $(n=2m)$}{d=10m+1}}\quad\\ 

We have 
\begin{align*}
  \overline{B} 
  &= \overline{C}+\overline{M_{\infty}}+\overline{L_{\infty}}
  +E_{1}+ \dots + E_{2m} + E_{\mathrm{fin}}^{(1)}, \\
  \overline{Z} 
  &= 3\overline{M_{\infty}}+(5m+1)\overline{L_{\infty}}
  +(5m+1)(E_{1}+ \dots + E_{2m})-E_{\mathrm{fin}}^{(2)}, \\
  K_{\overline{\Sigma}}+\overline{Z}
  &\sim \overline{M_{\infty}}+(5m-1)\overline{L_{\infty}} \\
  &+(5m-2)E_{1}+(5m-3)E_{2}+ \dots + (3m-1)E_{2m} 
  +(E_{\mathrm{can}}-E_{\mathrm{fin}}^{(2)}), \\
  p_{a}(\overline{Z}) 
  &= \frac{1}{2}(\overline{Z} \cdot K_{\overline{\Sigma}}+\overline{Z})+1 
  = 4m+\mu, \\
  (K_{\overline{\Sigma}}+\overline{Z})^{2} 
  &= 2m-2+\nu, 
\end{align*} 
\noindent where $\mu$ and $\nu$ are introduced in the end of the previous section. 
Reflecting the data about the branch locus, 
we have the configuration as in Figure \ref{sec:pf.of.main.th.2}-1-1A,  
\begin{figure}[H]
  \centering
  \begin{tikzpicture}
    \draw (0.5,1) to [out=225,in=135] (0.25,0.25) 
    to [out=315,in=225] (1.5,0.5)node[below]{$\overline{C}$}; 
    \draw (1,1)node[right]{$\overline{M_{\infty}}$}--(0,0); 
    \draw (0.5,0)node[below]{$E_{2m}$}--(-0.5,1); 
    \draw (0,1)--(-1,0)node[below]{$E_{2m-1}$};   
    \node (tenten) at (-1.5,0.5) {$\cdots$}; 
    \draw (-2,1)--(-3,0)node[below]{$E_{2}$};  
    \draw (-2.5,0)--(-3.5,1)node[above]{$E_{1}$}; 
    \draw (-3,1)--(-4,0)node[below]{$\overline{L_{\infty}}$};  
    \node (Lselfint) at (-4,-0.8) {$\Circled{-1}$}; 
    \node (E1selfint) at (-3.5,1.7) {$\Circled{-2}$}; 
    \node (E2selfint) at (-3,-0.7) {$\Circled{-2}$}; 
    \node (En-1selfint) at (-1,-0.7) {$\Circled{-2}$}; 
    \node (Enselfint) at (0.5,-0.7) {$\Circled{-1}$}; 
    \node (Mselfint) at (1.5,1.5) {$\Circled{-2m}$}; 
    \node (intersection) at (3,0) {$(\overline{C} \cdot E_{2m})=5$};
  \end{tikzpicture} 
  \caption*{Figure \ref{sec:pf.of.main.th.2}-1-1A}   
\end{figure} 
\noindent where the solid lines are contained in $\overline{B}$, while 
the broken lines are not contained in $\overline{B}$ 
as in the previous subsection. 
After some blowing-ups, 
we have the dual graph on $S$ by taking a double covering of $\Sigma$ as in Figure \ref{sec:pf.of.main.th.2}-1-1B,     
\begin{figure}[H] 
  \centering 
  \scalebox{0.8}{
  \begin{tikzpicture}
    \node[draw,shape=circle,inner sep=2pt] (a0) at (0,0) {};
    \node[draw,shape=circle,inner sep=2pt,fill=black] (a1) at (1,0) {};
    \node[draw,shape=circle,inner sep=2pt] (a2) at (2,0) {};
    \node[draw,shape=circle,inner sep=2pt,fill=black] (a3) at (3,0) {};
    \node (a4) at (4,0) {$\cdots$};
    \node[draw,shape=circle,inner sep=2pt,fill=black] (a5) at (5,0) {};
    \node[draw,shape=circle,inner sep=2pt] (a6) at (6,0) {};
    \node[draw,shape=circle,inner sep=2pt,fill=black] (a7) at (7,0) {};
    \node[draw,shape=circle,inner sep=2pt] (a8) at (8,0) {};
    \node[draw,shape=circle,inner sep=2pt,fill=black] (a9) at (9,0) {};
    \node[draw,shape=circle,inner sep=2pt] (a10) at (10,0) {};
    \node[draw,shape=circle,inner sep=2pt,fill=black] (a11) at (11,0) {};
    \node[draw,shape=circle,inner sep=2pt] (l5) at (11,-1) {};
    \node[draw,shape=circle,inner sep=2pt] (b) at (9,-1) {}; 
    \node (takedots) at (11,-2) {$\vdots$}; 
    \node[draw,shape=circle,inner sep=2pt] (l11) at (11,-3) {}; 
    \node[draw,shape=circle,inner sep=2pt,fill=black] (m) at (11,-4) {}; 
    \draw (a0)--(a1)--(a2)--(a3)--(3.5,0);
    \draw (4.5,0)--(a5)--(a6)--(a7)--(a8)--(a9)--(a10)--(a11)--(l5);
    \draw (a9)--(b);
    \draw (l5)--(11,-1.5); 
    \draw (11,-2.5)--(l11); 
    \draw (l11)--(m); 
      \node (a0name) at (0,-0.5) {$\widetilde{L_{\infty}}$}; 
      \node (a1name) at (1,-0.5) {$\widetilde{L_{\infty}'}$}; 
      \node (a2name) at (2,-0.5) {$\widetilde{E_{1}}$};
      \node (a3name) at (3,-0.5) {$\widetilde{E_{1}'}$};
      \node (a5name) at (5,-0.5) {$\widetilde{E_{2m-1}}$};
      \node (a6name) at (6,-0.5) {$\widetilde{E_{2m-1}'}$};
      \node (a7name) at (7,-0.5) {$\widetilde{E_{2m}}$};
      \node (a8name) at (8,-0.5) {$L_{1}$}; 
      \node (a9name) at (9.5,-0.5) {$L_{2}$}; 
      \node (a10name) at (10.3,-0.5) {$L_{3}$};
      \node (a11name) at (11.5,0) {$L_{4}$}; 
      \node (a12name) at (11.5,-1) {$L_{5}$}; 
      \node (bname) at (9.5,-1.2) {$L_{2}'$}; 
      \node (l11name) at (11.5,-3) {$L_{11}$}; 
      \node (mname) at (11.5,-4) {$\widetilde{M_{\infty}}$};  
      \node (l0selfint) at (0,0.5) {$-1$};
      \node (a1selfint) at (1,0.5) {$-2$};
      \node (e1selfint) at (2,0.5) {$-2$};
      \node (a3selfint) at (3,0.5) {$-2$};
      \node (e2elfint) at (5,0.5) {$-2$};
      \node (a6selfint) at (6,0.5) {$-2$};
      \node (a7selfint) at (7,0.5) {$-4$};
      \node (a8selfint) at (8,0.5) {$-2$};
      \node (a9selfint) at (9,0.5) {$-2$};
      \node (a10selfint) at (10,0.5) {$-2$};
      \node (a11selfint) at (11,0.5) {$-2$};
      \node (bselfint) at (8.5,-1.2) {$-2$};
      \node (l11selfint) at (10.5,-3) {$-2$};
      \node (mselfint) at (10,-4) {$-(m+1)$}; 
  \end{tikzpicture}
  } 
  \caption*{Figure \ref{sec:pf.of.main.th.2}-1-1B}  
\end{figure}  
\noindent where 
$\pi^{*}(\overline{L_{\infty}})=2\widetilde{L_{\infty}}+\widetilde{L_{\infty}'}$,  
$\pi^{*}(E_{1})=\widetilde{L_{\infty}'}+2\widetilde{E_{1}}+\widetilde{E_{1}'}$, 
$\pi^{*}(E_{i})=\widetilde{E_{i-1}'}+2\widetilde{E_{i}}+\widetilde{E_{i}'}$ for $2 \le i \le 2m-1$, 
$\pi^{*}(E_{2m})=\widetilde{E_{2m-1}'}+2\widetilde{E_{2m}}
+6L_{1}+5L_{2}'+\sum_{i=2}^{11}(12-i)L_{i}$ 
and 
$\pi^{*}(\overline{M_{\infty}})=2\widetilde{M_{\infty}}+L_{1}+L_{2}'+\sum_{i=2}^{11}L_{i}$. 

By contracting $\widetilde{L_{\infty}}$, $\widetilde{E_{1}}$, $\dots$, 
$\widetilde{E_{2m-1}}$, $\widetilde{E_{2m-1}'}$,  
and 
other $(-1)$-curves contained in reducible fibers in the affine part,        
we have the relatively minimal model $X$. 
By relabeling suitably, we have the dual graph on $X$ as in Figure \ref{sec:pf.of.main.th.2}-1-1C.    
\begin{figure}[H]
  \centering
  \scalebox{0.8}{
\begin{tikzpicture} 
  \node[draw,shape=circle] (0) at (13,0) {};
  \node[draw,shape=circle] (a0) at (12,0) {};
  \node[draw,shape=circle] (a1) at (11,0) {};
  \node[draw,shape=circle] (a2) at (10,0) {};
  \node[draw,shape=circle] (a3) at (9,0) {};
  \node[draw,shape=circle] (a4) at (8,0) {};
  \node[draw,shape=circle] (a5) at (7,0) {};
  \node[draw,shape=circle] (a6) at (6,0) {};
  \node[draw,shape=circle] (a7) at (5,0) {};
  \node[draw,shape=circle] (a8) at (4,0) {};
  \node[draw,shape=circle] (a9) at (3,0) {};
  \node[draw,shape=circle] (a10) at (3,-1) {};
  \node[draw,shape=circle] (a11) at (2,0) {};
  \node[draw,shape=circle] (a12) at (1,0) {};
  \draw (0)--(a0)--(a1)--(a2)--(a3)--(a4)--(a5)--(a6)
  --(a7)--(a8)--(a9)--(a11)--(a12);
  \draw (a9)--(a10);
  \node (0name) at (13,-0.5) {$(O)$};
  \node (a0name) at (12,-0.5) {$A_{0}$};
  \node (a1name) at (11,-0.5) {$A_{1}$};
  \node (a2name) at (10,-0.5) {$A_{2}$};
  \node (a3name) at (9,-0.5) {$A_{3}$};
  \node (a4name) at (8,-0.5) {$A_{4}$};
  \node (a5name) at (7,-0.5) {$A_{5}$};
  \node (a6name) at (6,-0.5) {$A_{6}$};
  \node (a7name) at (5,-0.5) {$A_{7}$};
  \node (a8name) at (4,-0.5) {$A_{8}$};
  \node (a9name) at (3.3,-0.5) {$A_{9}$};
  \node (a10name) at (3.3,-1.5) {$A_{10}$};
  \node (a11name) at (2,-0.5) {$A_{11}$};
  \node (a12name) at (1,-0.5) {$A_{12}$};
  \node (a12selfint) at (1,0.5) {$-3$}; 
  \node (a11selfint) at (1,0.5) {};
  \node (a10selfint) at (2.5,-1) {};
  \node (a9selfint) at (2,0.5) {};
  \node (a8selfint) at (3,0.5) {};
  \node (a7selfint) at (4,0.5) {};
  \node (a6selfint) at (5,0.5) {};
  \node (a5selfint) at (6,0.5) {};
  \node (a4selfint) at (7,0.5) {};
  \node (a3selfint) at (8,0.5) {};
  \node (a2selfint) at (9,0.5) {};
  \node (a1selfint) at (10,0.5) {};
  \node (Zerosecselfint) at (13,0.4) {$-(m+1)$}; 
  \end{tikzpicture}
  } 
  \caption*{Figure \ref{sec:pf.of.main.th.2}-1-1C}  
\end{figure}
Then, we obtain an explicit expression of 
the canonical divisor of $X$ (cf.\ Lemma \ref{lem:Artin}). 
\[
  K_{X} \sim 
  2(O)+(3m-1)F+2\sum_{i=0}^{9}A_{i}+A_{10}+A_{11}
  +\pi^{*}(E_{\mathrm{can}}-E_{\mathrm{fin}}^{(2)}),   
\]
where $F$ is a general fiber.

\medskip 

\paragraph*{Subcase 1-2 \texorpdfstring{$d=10m+6$ $(n=2m+1)$}{d=10m+6}}\quad\\ 

We have 
\begin{align*}
  \overline{B} 
  &= \overline{C}+\overline{M_{\infty}}+E_{\mathrm{fin}}^{(1)}, \\
  \overline{Z} 
  &= 3\overline{M_{\infty}}+(5m+3)\overline{L_{\infty}}
  +(5m+3)(E_{1}+ \dots + E_{2m+1})-E_{\mathrm{fin}}^{(2)}, \\
  K_{\overline{\Sigma}}+\overline{Z}
  &\sim \overline{M_{\infty}}+(5m-1)\overline{L_{\infty}} \\
  &+5mE_{1}+(5m+1)E_{2}+ \dots + 3mE_{2m+1} 
  +(E_{\mathrm{can}}-E_{\mathrm{fin}}^{(2)}), \\
  p_{a}(\overline{Z}) 
  &= 4m+1+\mu, \\
  (K_{\overline{\Sigma}}+\overline{Z})^{2} 
  &= 2m-2+\nu. 
\end{align*} 

Reflecting the data about the branch locus, 
we have the configuration as in Figure \ref{sec:pf.of.main.th.2}-1-2A,      
\begin{figure}[H]
  \centering
  \begin{tikzpicture}
    \draw (0.5,1) to [out=225,in=135] (0.25,0.25) 
    to [out=315,in=225] (1.5,0.5)node[below]{$\overline{C}$}; 
    \draw (1,1)node[right]{$\overline{M_{\infty}}$}--(0,0); 
    \draw[dashed] (0.5,0)node[below]{$E_{2m+1}$}--(-0.5,1); 
    \draw[dashed] (0,1)--(-1,0)node[below]{$E_{2m}$};   
    \node (tenten) at (-1.5,0.5) {$\cdots$}; 
    \draw[dashed] (-2,1)--(-3,0)node[below]{$E_{2}$};  
    \draw[dashed] (-2.5,0)--(-3.5,1)node[above]{$E_{1}$}; 
    \draw[dashed] (-3,1)--(-4,0)node[below]{$\overline{L_{\infty}}$};  
    \node (Lselfint) at (-4,-0.8) {$\Circled{-1}$}; 
    \node (E1selfint) at (-3.5,1.7) {$\Circled{-2}$}; 
    \node (E2selfint) at (-3,-0.7) {$\Circled{-2}$}; 
    \node (En-1selfint) at (-1,-0.7) {$\Circled{-2}$}; 
    \node (Enselfint) at (0.5,-0.7) {$\Circled{-1}$}; 
    \node (Mselfint) at (1.5,1.5) {$\Circled{-(2m+1)}$}; 
    \node (intersection) at (3,0) {$(\overline{C} \cdot E_{2m+1})=5$};
  \end{tikzpicture} 
  \caption*{Figure \ref{sec:pf.of.main.th.2}-1-2A}  
\end{figure} 
\noindent where the solid lines are contained in $\overline{B}$, while 
the broken lines are not contained in $\overline{B}$ 
as in the previous subsection.  
After some blowing-ups, 
we have the dual graph on $S$ 
by taking a double covering of $\Sigma$ as in Figure \ref{sec:pf.of.main.th.2}-1-2B,     
\begin{figure}[H] 
  \centering
  \scalebox{0.8}{
  \begin{tikzpicture}
    \node[draw,shape=circle,inner sep=2pt] (linf) at (3,1) {};
    \node[draw,shape=circle,inner sep=2pt] (llinf) at (3,-1) {};
    \node[draw,shape=circle,inner sep=2pt] (e1) at (4,1) {};
    \node[draw,shape=circle,inner sep=2pt] (ee1) at (4,-1) {}; 
    \node (ueyoko) at (5,1) {$\cdots$}; 
    \node (shitayoko) at (5,-1) {$\cdots$}; 
    \node[draw,shape=circle,inner sep=2pt] (e2m) at (6,1) {};
    \node[draw,shape=circle,inner sep=2pt] (ee2m) at (6,-1) {};
    \node[draw,shape=circle,inner sep=2pt] (e2m+1) at (7,1) {};
    \node[draw,shape=circle,inner sep=2pt] (ee2m+1) at (7,-1) {};
    \node[draw,shape=circle,inner sep=2pt] (l) at (8,0) {};
    \node[draw,shape=circle,inner sep=2pt,fill=black] (m) at (9,0) {};   
    \draw (linf)--(e1)--($(e1)+(0.5,0)$); 
    \draw (llinf)--(ee1)--($(ee1)+(0.5,0)$);  
    \draw ($(e2m)+(-0.5,0)$)--(e2m);  
    \draw ($(ee2m)+(-0.5,0)$)--(ee2m);   
    \draw (e2m)--(e2m+1)--(l); 
    \draw (ee2m)--(ee2m+1)--(l);  
    \draw (l)--(m);   
    \draw[double distance=2pt,nfold=2] (e2m+1)--(ee2m+1); 
    \node (linfname) at ($(linf)+(0,-0.5)$) {$\widetilde{L_{\infty}}$}; 
    \node (llinfname) at ($(llinf)+(0,-0.5)$) {$\widetilde{L_{\infty}'}$}; 
    \node (e1name) at ($(e1)+(0,-0.5)$) {$\widetilde{E_{1}}$};  
    \node (ee1name) at ($(ee1)+(0,-0.5)$) {$\widetilde{E_{1}'}$}; 
    \node (e2mname) at ($(e2m)+(0,-0.5)$) {$\widetilde{E_{2m}}$};  
    \node (ee2mname) at ($(ee2m)+(0,-0.5)$) {$\widetilde{E_{2m}'}$}; 
    \node (e2m+1name) at ($(e2m+1)+(1,0.5)$) {$\widetilde{E_{2m+1}}$};  
    \node (ee2m+1name) at ($(ee2m+1)+(1,-0.5)$) {$\widetilde{E_{2m+1}'}$};  
    \node (lname) at ($(l)+(0,-0.5)$) {$L$};
    \node (mname) at ($(m)+(0,-0.5)$) {$\widetilde{M_{\infty}}$};   
    \node (linfselfint) at ($(linf)+(0,0.5)$) {$-1$}; 
    \node (llinfselfint) at ($(llinf)+(0,0.5)$) {$-1$}; 
    \node (e1selfint) at ($(e1)+(0,0.5)$) {$-2$};  
    \node (ee1selfint) at ($(ee1)+(0,0.5)$) {$-2$}; 
    \node (e2m+1selfint) at ($(e2m+1)+(0,0.5)$) {$-4$};  
    \node (ee2m+1selfint) at ($(ee2m+1)+(0,-0.5)$) {$-4$}; 
    \node (e2mselfint) at ($(e2m)+(0,0.5)$) {$-2$};  
    \node (ee2mselfint) at ($(ee2m)+(0,0.5)$) {$-2$};  
    \node (lselfint) at ($(l)+(0,0.5)$) {$-2$}; 
    \node (mselfint) at ($(m)+(0,0.5)$) {$-(m+1)$};   
  \end{tikzpicture}  
  }
  \caption*{Figure \ref{sec:pf.of.main.th.2}-1-2B}  
\end{figure} 
\noindent where 
$\pi^{*}(\overline{L_{\infty}})=\widetilde{L_{\infty}}+\widetilde{L_{\infty}'}$, 
$\pi^{*}(E_{i})=\widetilde{E_{i}}+\widetilde{E_{i}'}$ for $1 \le i \le 2m$, 
$\pi^{*}(E_{2m+1})=\widetilde{E_{2m+1}}+\widetilde{E_{2m+1}'}+L$  
and 
$\pi^{*}(\overline{M_{\infty}})=2\widetilde{M_{\infty}}+L$. 
Note that components $\widetilde{E_{2m+1}}, \widetilde{E_{2m+1}'}$ and $L$ meet 
in one point with $(\widetilde{E_{2m+1}} \cdot \widetilde{E_{2m+1}'})=2$ and 
$(\widetilde{E_{2m+1}} \cdot L)=(\widetilde{E_{2m+1}'} \cdot L)=1$. 

By contracting 
$\widetilde{L_{\infty}}$, $\widetilde{E_{1}}$, $\dots$, $\widetilde{E_{2m}}$ 
and 
$\widetilde{L_{\infty}'}$, $\widetilde{E_{1}'}$, $\dots$, $\widetilde{E_{2m}'}$ 
and 
other $(-1)$-curves contained in reducible fibers in the affine part,        
we have the relatively minimal model $X$. 
By relabeling suitably, we have the dual graph on $X$ as in Figure \ref{sec:pf.of.main.th.2}-1-2C.     
\begin{figure}[H]
  \centering
  \scalebox{0.8}{
    \begin{tikzpicture} 
      \node[draw,shape=circle] (0) at (2.5,0) {};
      \node[draw,shape=circle] (a0) at (1.5,0) {};
      \node[draw,shape=circle] (a1) at (0.5,0.5) {};
      \node[draw,shape=circle] (a2) at (0.5,-0.5) {};
      \draw (a1)--(a0)--(a2);
      \draw (0)--(a0); 
      \draw[double distance=4pt] (a1)--(a2); 
      \node (0name) at (2.5,-0.5) {$(O)$};
      \node (a0name) at (1.5,-0.5) {$A_{0}$};
      \node (a1name) at (0,0.5) {$A_{1}$};
      \node (a2name) at (0,-0.5) {$A_{2}$};
      \node (0selfint) at (2.5,0.5) {$-(m+1)$};
      \node (a0selfint) at (1.5,0.3) {};
      \node (a1selfint) at (0.5,0.9) {$-3$};
      \node (a2selfint) at (0.5,-0.9) {$-3$};
    \end{tikzpicture}
  } 
  \caption*{Figure \ref{sec:pf.of.main.th.2}-1-2C}   
\end{figure}
Then, we obtain an explicit expression of 
the canonical divisor of $X$.
\[
  K_{X} \sim 
  2(O)+3mF+A_{0}
  +\pi^{*}(E_{\mathrm{can}}-E_{\mathrm{fin}}^{(2)}),   
\]
where $F$ is a general fiber. 

\subsubsection*{Case \texorpdfstring{$2$}{2} \texorpdfstring{$d=5n+2$ $(n>0)$}{d=5n+2}} 

We have to blow-up ($n+2$) times to obtain the nonsingular curve $\overline{C}$ (Figure \ref{sec:pf.of.main.th.2}-2). 
\begin{figure}[H]  
  \centering
  \begin{tikzpicture}
    \draw (0.5,1)--(1.5,0)node[below]{$E_{n+1}$}; 
    \draw (1,0)--(2,1)node[right]{$\overline{M_{\infty}}$};  
    \draw (-0.5,0) to [out=0,in=225] (0.25,0.25) 
    to [out=45,in=315] (0.2,0.9)node[above]{$\overline{C}$};  
    \draw (1,1)node[above]{$E_{n+2}$}--(0,0); 
    \draw (0.5,0)node[below]{$E_{n}$}--(-0.5,1); 
    \draw (0,1)--(-1,0)node[below]{$E_{n-1}$};   
    \node (tenten) at (-1.5,0.5) {$\cdots$}; 
    \draw (-2,1)--(-3,0)node[below]{$E_{2}$};  
    \draw (-2.5,0)--(-3.5,1)node[above]{$E_{1}$}; 
    \draw (-3,1)--(-4,0)node[below]{$\overline{L_{\infty}}$};  
    \node (Lselfint) at (-4,-0.8) {$\Circled{-1}$}; 
    \node (E1selfint) at (-3.5,1.7) {$\Circled{-2}$}; 
    \node (E2selfint) at (-3,-0.7) {$\Circled{-2}$}; 
    \node (En-1selfint) at (-1,-0.7) {$\Circled{-2}$};
    \node (Enselfint) at (0.5,-0.7) {$\Circled{-3}$}; 
    \node (En1selfint) at (1.5,-0.7) {$\Circled{-2}$}; 
    \node (En2selfint) at (1.3,1.7) {$\Circled{-1}$}; 
    \node (Mselfint) at (2.5,1.5) {$\Circled{-(n+1)}$}; 
    \node (intersection) at (3,0) {$(\overline{C} \cdot E_{n+2})=2$}; 
  \end{tikzpicture} 
  \caption*{Figure \ref{sec:pf.of.main.th.2}-2}   
\end{figure}
Then, we have 
\begin{align*}
  \overline{\sigma}^{*}A 
  &= \overline{\sigma}^{*}C - 5\overline{\sigma}^{*}M_{\infty} 
  -d\overline{\sigma}^{*}L_{\infty} \\
  &= (\overline{C}+5(E_{1}+2E_{2}+ \dots + nE_{n})
  +(5n+2)E_{n+1}+(10n+4)E_{n+2}+E_{\mathrm{fin}}) \\ 
  & -5(\overline{M_{\infty}}+E_{1}+2E_{2}+ \dots +nE_{n}
  +(n+1)E_{n+1}+(2n+1)E_{n+2})  \\ 
  &-d(\overline{L_{\infty}}+E_{1}+E_{2}+ \dots +E_{n}
  +E_{n+1}+2E_{n+2}) \\ 
  &= \overline{C}-5\overline{M_{\infty}}-d\overline{L_{\infty}}
  -d(E_{1}+E_{2}+ \dots +E_{n}) \\
  &-(d+3)E_{n+1}-(2d+1)E_{n+2} 
  + E_{\mathrm{fin}},  
\end{align*}
and 
\begin{align*}
  K_{\overline{\Sigma}} 
  &\sim -2\overline{\sigma}^{*}M_{\infty} 
  -2\overline{\sigma}^{*}L_{\infty} +(E_{1}+2E_{2}+ \dots nE_{n}) \\
  &+(n+1)E_{n+1}+(2n+2)E_{n+2}  
  +E_{\mathrm{can}} \\ 
  &= -2\overline{M_{\infty}}-2\overline{L_{\infty}}
  -(3E_{1}+4E_{2}+ \dots + (n+2)E_{n}) \\
  &-(n+3)E_{n+1}-(2n+4)E_{n+2}+E_{\mathrm{can}}.  
\end{align*} 

\medskip 

\paragraph*{Subcase 2-1 \texorpdfstring{$d=10m+2$ $(n=2m)$}{d=10m+2}}\quad\\ 

We have 
\begin{align*} 
  \overline{B} 
  &= \overline{C}+\overline{M_{\infty}} 
  +E_{2m+1}+E_{2m+2}+E_{\mathrm{fin}}^{(1)}, \\
  \overline{Z} 
  &= 3\overline{M_{\infty}}+(5m+1)\overline{L_{\infty}}
  +(5m+1)(E_{1}+ \dots + E_{2m}) \\
  &+(5m+3)E_{2m+1}+(10m+3)E_{2m+2}-E_{\mathrm{fin}}^{(2)}, \\
  K_{\overline{\Sigma}}+\overline{Z} 
  &\sim \overline{M_{\infty}}+(5m-1)\overline{L_{\infty}} \\
  &+(5m-2)E_{1}+(5m-3)E_{2}+ \dots + (3m-1)E_{2m} \\ 
  &+3mE_{2m+1}+(6m-1)E_{2m+2}+(E_{\mathrm{can}}-E_{\mathrm{fin}}^{(2)}), \\
  p_{a}(\overline{Z}) 
  &= 4m+\mu, \\
  (K_{\overline{\Sigma}}+\overline{Z})^{2} 
  &= 2m-2+\nu. 
\end{align*}
Reflecting the data about the branch locus, 
we have the configuration as in Figure \ref{sec:pf.of.main.th.2}-2-1A.  
\begin{figure}[H] 
  \centering
  \begin{tikzpicture}
    \draw (0.5,1)--(1.5,0)node[below]{$E_{2m+1}$}; 
    \draw (1,0)--(2,1)node[right]{$\overline{M_{\infty}}$};  
    \draw (-0.5,0) to [out=0,in=225] (0.25,0.25) 
    to [out=45,in=315] (0.2,0.9)node[above]{$\overline{C}$};  
    \draw (1,1)node[above]{$E_{2m+2}$}--(0,0); 
    \draw[dashed] (0.5,0)node[below]{$E_{2m}$}--(-0.5,1); 
    \draw[dashed] (0,1)--(-1,0)node[below]{$E_{2m-1}$};   
    \node (tenten) at (-1.5,0.5) {$\cdots$}; 
    \draw[dashed] (-2,1)--(-3,0)node[below]{$E_{2}$};  
    \draw[dashed] (-2.5,0)--(-3.5,1)node[above]{$E_{1}$}; 
    \draw[dashed] (-3,1)--(-4,0)node[below]{$\overline{L_{\infty}}$};  
    \node (Lselfint) at (-4,-0.8) {$\Circled{-1}$}; 
    \node (E1selfint) at (-3.5,1.7) {$\Circled{-2}$}; 
    \node (E2selfint) at (-3,-0.7) {$\Circled{-2}$}; 
    \node (En-1selfint) at (-1,-0.7) {$\Circled{-2}$};
    \node (Enselfint) at (0.5,-0.7) {$\Circled{-3}$}; 
    \node (En1selfint) at (1.5,-0.7) {$\Circled{-2}$}; 
    \node (En2selfint) at (1.3,1.7) {$\Circled{-1}$}; 
    \node (Mselfint) at (2.5,1.5) {$\Circled{-(2m+1)}$}; 
    \node (intersection) at (3,0) {$(\overline{C} \cdot E_{2m+2})=2$}; 
  \end{tikzpicture} 
  \caption*{Figure \ref{sec:pf.of.main.th.2}-2-1A}  
\end{figure}
After some blowing-ups, 
we have the dual graph on $S$ by taking a double covering of $\Sigma$ as in Figure \ref{sec:pf.of.main.th.2}-2-1B, 
\begin{figure}[H]
  \centering
  \scalebox{0.8}{
  \begin{tikzpicture}
    \node[draw,shape=circle,inner sep=2pt] (linf) at (0,1) {};
    \node[draw,shape=circle,inner sep=2pt] (llinf) at (0,-1) {};
    \node[draw,shape=circle,inner sep=2pt] (e1) at (1,1) {};
    \node[draw,shape=circle,inner sep=2pt] (ee1) at (1,-1) {}; 
    \node (ueyoko) at (2,1) {$\cdots$}; 
    \node (shitayoko) at (2,-1) {$\cdots$}; 
    \node[draw,shape=circle,inner sep=2pt] (e2m) at (3,1) {};
    \node[draw,shape=circle,inner sep=2pt] (ee2m) at (3,-1) {};
    \node[draw,shape=circle,inner sep=2pt] (l) at (4,1) {};
    \node[draw,shape=circle,inner sep=2pt] (ll) at (4,-1) {}; 
    \node[draw,shape=circle,inner sep=2pt] (l1) at (5,0) {};
    \node[draw,shape=circle,inner sep=2pt,fill=black] (e2m+2) at (6,0) {}; 
    \node[draw,shape=circle,inner sep=2pt] (l2) at (7,0) {};  
    \node[draw,shape=circle,inner sep=2pt,fill=black] (e2m+1) at (8,0) {}; 
    \node[draw,shape=circle,inner sep=2pt] (l3) at (9,0) {}; 
    \node[draw,shape=circle,inner sep=2pt,fill=black] (m) at (10,0) {};   
    \draw (linf)--(e1)--(1.5,1); 
    \draw (llinf)--(ee1)--(1.5,-1);  
    \draw (2.5,1)--(e2m);  
    \draw (2.5,-1)--(ee2m); 
    \draw (e2m)--(l)--(l1); 
    \draw (ee2m)--(ll)--(l1); 
    \draw (l1)--(e2m+2)--(l2)--(e2m+1)--(l3)--(m); 
    \node (linfname) at (0,0.5) {$\widetilde{L_{\infty}}$}; 
    \node (llinfname) at (0,-1.5) {$\widetilde{L_{\infty}'}$}; 
    \node (e1name) at (1,0.5) {$\widetilde{E_{1}}$};  
    \node (ee1name) at (1,-1.5) {$\widetilde{E_{1}'}$}; 
    \node (e2mname) at (3,0.5) {$\widetilde{E_{2m}}$};  
    \node (ee2mname) at (3,-1.5) {$\widetilde{E_{2m}'}$}; 
    \node (lmname) at (4,0.5) {$\widetilde{L}$}; 
    \node (llmname) at (4,-1.5) {$\widetilde{L'}$};   
    \node (l1name) at (5,-0.5) {$L_{1}$}; 
    \node (e2m+2name) at (6,-0.5) {$\widetilde{E_{2m+2}}$};  
    \node (l2name) at (7,-0.5) {$L_{2}$}; 
    \node (e2m+1) at (8,-0.5) {$\widetilde{E_{2m+1}}$}; 
    \node (l3name) at (9,-0.5) {$L_{3}$}; 
    \node (mname) at (10,-0.5) {$\widetilde{M_{\infty}}$};   
    \node (linfselfint) at (0,1.5) {$-1$}; 
    \node (llinfselfint) at (0,-0.5) {$-1$}; 
    \node (e1selfint) at (1,1.5) {$-2$};  
    \node (ee1selfint) at (1,-0.5) {$-2$}; 
    \node (e2mselfint) at (3,1.5) {$-4$};  
    \node (ee2mselfint) at (3,-0.5) {$-4$}; 
    \node (lselfint) at (4,1.5) {$-2$}; 
    \node (llselfint) at (4,-0.5) {$-2$};   
    \node (l1selfint) at (5,0.5) {$-2$}; 
    \node (e2m+2selfint) at (6,0.5) {$-2$};  
    \node (l2selfint) at (7,0.5) {$-2$}; 
    \node (e2m+1selfint) at (8,0.5) {$-2$}; 
    \node (l3selfint) at (9,0.5) {$-2$}; 
    \node (mselfint) at (10,0.5) {$-(m+1)$};  
  \end{tikzpicture}  
  } 
  \caption*{Figure \ref{sec:pf.of.main.th.2}-2-1B}   
\end{figure} 
\noindent where 
$\pi^{*}(\overline{L_{\infty}})=\widetilde{L_{\infty}}+\widetilde{L_{\infty}'}$, 
$\pi^{*}(E_{i})=\widetilde{E_{i}}+\widetilde{E_{i}'}$ for $1 \le i \le 2m-1$, 
$\pi^{*}(E_{2m})=\widetilde{E_{2m}}+\widetilde{E_{2m}'}+L+L'+L_{1}$, 
$\pi^{*}(E_{2m+2})=L+L'+2L_{1}+2\widetilde{E_{2m+2}}+L_{2}$, 
$\pi^{*}(E_{2m+1})=L_{2}+2\widetilde{E_{2m+1}}+L_{3}$ 
and 
$\pi^{*}(\overline{M_{\infty}})=2\widetilde{M_{\infty}}+L_{3}$.

By contracting $\widetilde{L_{\infty}}$, $\widetilde{E_{1}}$, 
$\dots$, $\widetilde{E_{2m-1}}$ 
and 
$\widetilde{L_{\infty}'}$, $\widetilde{E_{1}'}$, 
$\dots$, $\widetilde{E_{2m-1}'}$ 
and 
other $(-1)$-curves contained in reducible fibers in the affine part,        
we have the relatively minimal model $X$. 
By relabeling suitably, we have the dual graph on $X$ as in Figure \ref{sec:pf.of.main.th.2}-2-1C.    
\begin{figure}[H]
  \centering
  \scalebox{0.8}{
    \begin{tikzpicture} 
      \node[draw,shape=circle] (0) at (5,0) {};
      \node[draw,shape=circle] (a0) at (4,0) {};
      \node[draw,shape=circle] (a1) at (3,0) {};
      \node[draw,shape=circle] (a2) at (2,0) {}; 
      \node[draw,shape=circle] (a3) at (1,0) {};
      \node[draw,shape=circle] (a4) at (0,0) {};
      \node[draw,shape=circle] (a5) at (-1,1) {}; 
      \node[draw,shape=circle] (a6) at (-2,1) {}; 
      \node[draw,shape=circle] (a7) at (-1,-1) {}; 
      \node[draw,shape=circle] (a8) at (-2,-1) {};  
      \draw (0)--(a0)--(a1)--(a2)--(a3)--(a4);
      \draw (a4)--(a5)--(a6); 
      \draw (a4)--(a7)--(a8);  
      \node (0name) at ($(0)+(0,-0.5)$) {$(O)$};
      \node (a0name) at (4,-0.5) {$A_{0}$};
      \node (a1name) at (3,-0.5) {$A_{1}$};
      \node (a2name) at (2,-0.5) {$A_{2}$}; 
      \node (a3name) at (1,-0.5) {$A_{3}$};
      \node (a4name) at (0,-0.5) {$A_{4}$};
      \node (a5name) at (-1,0.5) {$A_{5}$}; 
      \node (a6name) at (-2,0.5) {$A_{6}$};
      \node (a7name) at (-1,-1.5) {$A_{7}$};
      \node (a8name) at (-2,-1.5) {$A_{8}$}; 
      \node (0selfint) at (5,0.5) {$-(m+1)$};
      \node (a6selfint) at (-2,1.5) {$-3$};
      \node (a8selfint) at (-2,-0.5) {$-3$};
    \end{tikzpicture}
  } 
  \caption*{Figure \ref{sec:pf.of.main.th.2}-2-1C}   
\end{figure}
Then, we obtain an explicit expression of 
the canonical divisor of $X$.
\[
  K_{X} \sim 
  2(O)+(3m-1)F+2\sum_{i=0}^{4}A_{i}+(A_{5}+A_{7})
  +\pi^{*}(E_{\mathrm{can}}-E_{\mathrm{fin}}^{(2)}),   
\]
where $F$ is a general fiber. 

\medskip 

\paragraph*{Subcase 2-2 \texorpdfstring{$d=10m+7$ $(n=2m+1)$}{d=10m+7}}\quad\\ 

We have 
\begin{align*}
  \overline{B} 
  &= \overline{C}+\overline{M_{\infty}}+\overline{L_{\infty}}
  +E_{1}+ \dots + E_{2m+1} +E_{2m+3}+ E_{\mathrm{fin}}^{(1)}, \\
  \overline{Z} 
  &= 3\overline{M_{\infty}}+(5m+4)\overline{L_{\infty}}
  +(5m+4)(E_{1}+ \dots + E_{2m+1}) \\
  &+(5m+5)E_{2m+2}+(10m+8)E_{2m+3}-E_{\mathrm{fin}}^{(2)}, \\
  K_{\overline{\Sigma}}+\overline{Z}
  &\sim \overline{M_{\infty}}+(5m+2)\overline{L_{\infty}} \\ 
  &+(5m+1)E_{1}+5mE_{2}+ \dots + (3m+1)E_{2m+1} \\ 
  &+(3m+1)E_{2m+2}+(6m+2)E_{2m+3}+(E_{\mathrm{can}}-E_{\mathrm{fin}}^{(2)}), \\
  p_{a}(\overline{Z}) 
  &= 4m+2+\mu, \\
  (K_{\overline{\Sigma}}+\overline{Z})^{2} 
  &= 2m-1+\nu. 
\end{align*}
Reflecting the data about the branch locus, 
we have the configuration as in Figure \ref{sec:pf.of.main.th.2}-2-2A.  
\begin{figure}[H] 
  \centering
  \begin{tikzpicture}
    \draw[dashed] (0.5,1)--(1.5,0)node[below]{$E_{2m+2}$}; 
    \draw (1,0)--(2,1)node[right]{$\overline{M_{\infty}}$};  
    \draw (-0.5,0) to [out=0,in=225] (0.25,0.25) 
    to [out=45,in=315] (0.2,0.9)node[above]{$\overline{C}$};  
    \draw (1,1)node[above]{$E_{2m+3}$}--(0,0); 
    \draw (0.5,0)node[below]{$E_{2m+1}$}--(-0.5,1); 
    \draw (0,1)--(-1,0)node[below]{$E_{2m}$};   
    \node (tenten) at (-1.5,0.5) {$\cdots$}; 
    \draw (-2,1)--(-3,0)node[below]{$E_{2}$};  
    \draw (-2.5,0)--(-3.5,1)node[above]{$E_{1}$}; 
    \draw (-3,1)--(-4,0)node[below]{$\overline{L_{\infty}}$};  
    \node (Lselfint) at (-4,-0.8) {$\Circled{-1}$}; 
    \node (E1selfint) at (-3.5,1.7) {$\Circled{-2}$}; 
    \node (E2selfint) at (-3,-0.7) {$\Circled{-2}$}; 
    \node (En-1selfint) at (-1,-0.7) {$\Circled{-2}$};
    \node (Enselfint) at (0.5,-0.7) {$\Circled{-3}$}; 
    \node (En1selfint) at (1.5,-0.7) {$\Circled{-2}$}; 
    \node (En2selfint) at (1.3,1.7) {$\Circled{-1}$}; 
    \node (Mselfint) at (2.5,1.5) {$\Circled{-(2m+2)}$}; 
    \node (intersection) at (3,0) {$(\overline{C} \cdot E_{2m+3})=2$}; 
  \end{tikzpicture} 
  \caption*{Figure \ref{sec:pf.of.main.th.2}-2-2A}  
\end{figure}
After some blowing-ups, 
we have the dual graph on $S$ by taking a double covering of $\Sigma$ as in Figure \ref{sec:pf.of.main.th.2}-2-2B, 
\begin{figure}[H] 
  \centering
  \scalebox{0.75}{
  \begin{tikzpicture}
    \node[draw,shape=circle,inner sep=2pt,fill=black] (linf) at (-4,0) {};
    \node[draw,shape=circle,inner sep=2pt] (llinf) at (-3,0) {};
    \node[draw,shape=circle,inner sep=2pt,fill=black] (e1) at (-2,0) {};  
    \node[draw,shape=circle,inner sep=2pt] (ee1) at (-1,0) {}; 
    \node (ueyoko) at (0,0) {$\cdots$}; 
    \node[draw,shape=circle,inner sep=2pt] (e2m) at (1,0) {};
    \node[draw,shape=circle,inner sep=2pt,fill=black] (e2m+1) at (2,0) {};
    \node[draw,shape=circle,inner sep=2pt] (l1) at (3,0) {};
    \node[draw,shape=circle,inner sep=2pt,fill=black] (l2) at (4,0) {}; 
    \node[draw,shape=circle,inner sep=2pt] (l3) at (5,0) {};
    \node[draw,shape=circle,inner sep=2pt,fill=black] (l4) at (6,0) {}; 
    \node[draw,shape=circle,inner sep=2pt] (l5) at (6,-1) {};
    \node[draw,shape=circle,inner sep=2pt] (l6) at (7,0) {}; 
    \node[draw,shape=circle,inner sep=2pt,fill=black] (e2m+3) at (8,0) {};
    \node[draw,shape=circle,inner sep=2pt] (e2m+2) at (9,0) {};  
    \node[draw,shape=circle,inner sep=2pt,fill=black] (m) at (10,0) {};     
    \draw (linf)--(llinf)--(e1)--(ee1);     
    \draw (ee1)--($(ee1)+(0.5,0)$); 
    \draw ($(e2m)+(-0.5,0)$)--(e2m);    
    \draw (e2m)--(e2m+1)--(l1)--(l2)--(l3)--(l4)--(l6)--(e2m+3)--(e2m+2)--(m);  
    \draw (l4)--(l5);   
    \node (linfname) at ($(linf)+(0,-0.5)$) {$\widetilde{L_{\infty}}$}; 
    \node (llinfname) at ($(llinf)+(0,-0.5)$) {$\widetilde{L_{\infty}'}$}; 
    \node (e1name) at ($(e1)+(0,-0.5)$) {$\widetilde{E_{1}}$};   
    \node (ee1name) at ($(ee1)+(0,-0.5)$) {$\widetilde{E_{1}'}$}; 
    \node (e2mname) at ($(e2m)+(0,-0.5)$) {$\widetilde{E_{2m}'}$};  
    \node (e2m+1name) at ($(e2m+1)+(0,-0.5)$) {$\widetilde{E_{2m+1}}$}; 
    \node (l1name) at ($(l1)+(0,-0.5)$) {$L_{1}$}; 
    \node (l2name) at ($(l2)+(0,-0.5)$) {$L_{2}$};  
    \node (l3name) at ($(l3)+(0,-0.5)$) {$L_{3}$}; 
    \node (l4name) at ($(l4)+(0.3,-0.5)$) {$L_{4}$}; 
    \node (l5name) at ($(l5)+(0.5,-0.3)$) {$L_{5}$}; 
    \node (l6name) at ($(l6)+(0,-0.5)$) {$L_{6}$};  
    \node (e2m+3name) at ($(e2m+3)+(0,-0.5)$) {$\widetilde{E_{2m+3}}$}; 
    \node (e2m+2name) at ($(e2m+2)+(0,-0.5)$) {$\widetilde{E_{2m+2}}$};  
    \node (mname) at ($(m)+(0,-0.5)$) {$\widetilde{M_{\infty}}$};  
    \node (linfselfint) at ($(linf)+(0,0.5)$) {$-1$}; 
    \node (llinfselfint) at ($(llinf)+(0,0.5)$) {$-2$};   
    \node (e1selfint) at ($(e1)+(0,0.5)$) {$-2$}; 
    \node (ee1selfint) at ($(ee1)+(0,0.5)$) {$-2$}; 
    \node (e2mselfint) at ($(e2m)+(0,0.5)$) {$-2$};   
    \node (e2m+1selfint) at ($(e2m+1)+(0,0.5)$) {$-3$}; 
    \node (l1selfint) at ($(l1)+(0,0.5)$) {$-2$}; 
    \node (l2selfint) at ($(l2)+(0,0.5)$) {$-2$};  
    \node (l3selfint) at ($(l3)+(0,0.5)$) {$-2$}; 
    \node (l4selfint) at ($(l4)+(0,0.5)$) {$-2$}; 
    \node (l5selfint) at ($(l5)+(-0.5,-0.3)$) {$-2$};  
    \node (l6selfint) at ($(l6)+(0,0.5)$) {$-2$};
    \node (e2m+3selfint) at ($(e2m+3)+(0,0.5)$) {$-2$}; 
    \node (e2m+2selfint) at ($(e2m+2)+(0,0.5)$) {$-4$};
    \node (mselfint) at ($(m)+(0,0.5)$) {$-(m+1)$};  
  \end{tikzpicture}  
  } 
  \caption*{Figure \ref{sec:pf.of.main.th.2}-2-2B}  
\end{figure} 
\noindent where 
$\pi^{*}(\overline{L_{\infty}})=2\widetilde{L_{\infty}}+\widetilde{L_{\infty}'}$, 
$\pi^{*}(E_{1})=\widetilde{L_{\infty}'}+2\widetilde{E_{1}}+\widetilde{E_{1}'}$, 
$\pi^{*}(E_{i})=\widetilde{E_{i-1}'}+2\widetilde{E_{i}}+\widetilde{E_{i}'}$ for $2 \le i \le 2m$, 
$\pi^{*}(E_{2m+1})=\widetilde{E_{2m}'}+2\widetilde{E_{2m+1}}
+2(L_{1}+L_{2}+L_{3}+L_{4})+L_{5}+L_{6}$,
$\pi^{*}(E_{2m+3})=2\widetilde{E_{2m+3}}+L_{1}+2L_{2}+3L_{3}+4L_{4}+2L_{5}+3L_{6}$, 
$\pi^{*}(E_{2m+2})=2\widetilde{E_{2m+2}}$, 
and 
$\pi^{*}(\overline{M_{\infty}})=2\widetilde{M_{\infty}}$.

By contracting $\widetilde{L_{\infty}}$, $\widetilde{E_{1}}$, $\dots$, 
$\widetilde{E_{2m-1}}$, $\widetilde{E_{2m}}$    
and 
other $(-1)$-curves contained in reducible fibers in the affine part,        
we have the relatively minimal model $X$. 
By relabeling suitably, we have the dual graph on $X$ as in Figure \ref{sec:pf.of.main.th.2}-2-2C.     
\begin{figure}[H]
  \centering
  \scalebox{0.8}{
    \begin{tikzpicture} 
      \node[draw,shape=circle] (0) at (1,0) {};
      \node[draw,shape=circle] (a0) at (0,0) {};
      \node[draw,shape=circle] (a1) at (-1,0) {};
      \node[draw,shape=circle] (a2) at (-2,0) {}; 
      \node[draw,shape=circle] (a3) at (-3,0) {};
      \node[draw,shape=circle] (a4) at (-3,-1) {}; 
      \node[draw,shape=circle] (a5) at (-4,0) {};
      \node[draw,shape=circle] (a6) at (-5,0) {};
      \node[draw,shape=circle] (a7) at (-6,0) {}; 
      \node[draw,shape=circle] (a8) at (-7,0) {};
      \draw (0)--(a0)--(a1)--(a2)--(a3)--(a5)--(a6)--(a7)--(a8);
      \draw (a3)--(a4); 
      \node (0name) at (1,-0.5) {$(O)$};
      \node (a0name) at (0,-0.5) {$A_{0}$};
      \node (a1name) at (-1,-0.5) {$A_{1}$};
      \node (a2name) at (-2,-0.5) {$A_{2}$};
      \node (a3name) at (-2.7,-0.5) {$A_{3}$}; 
      \node (a4name) at (-2.7,-1.5) {$A_{4}$}; 
      \node (a5name) at (-4,-0.5) {$A_{5}$}; 
      \node (a6name) at (-5,-0.5) {$A_{6}$}; 
      \node (a7name) at (-6,-0.5) {$A_{7}$}; 
      \node (a8name) at (-7,-0.5) {$A_{8}$}; 
      \node (0selfint) at (1,0.5) {$-(m+1)$};
      \node (a0selfint) at (0,0.5) {$-4$}; 
    \end{tikzpicture}
  } 
  \caption*{Figure \ref{sec:pf.of.main.th.2}-2-2C}   
\end{figure}
Then, we obtain an explicit expression of 
the canonical divisor of $X$.
\[
  K_{X} \sim 
  2(O)+(3m+1)F
  +\pi^{*}(E_{\mathrm{can}}-E_{\mathrm{fin}}^{(2)}),   
\]
where $F$ is a general fiber. 

\subsubsection*{Case \texorpdfstring{$3$}{3} \texorpdfstring{$d=5n+3$ $(n>0)$}{d=5n+3}}

We have to blow-up $(n+2)$ times to obtain the nonsigular curve $\overline{C}$ (Figure \ref{sec:pf.of.main.th.2}-3). 
\begin{figure}[H]
  \centering
  \begin{tikzpicture}
    \draw (0.4,0.8) to [out=225,in=135] (0.25,0.25) 
    to [out=315,in=225] (0.7,0.2)node[right]{$\overline{C}$}; 
    \draw (1,1)node[above]{$E_{n+1}$}--(0,0);  
    \draw (0.5,1)--(1.5,0)node[right]{$\overline{M_{\infty}}$};
    \draw (0.5,0)node[below]{$E_{n+2}$}--(-0.5,1); 
    \draw (0,1)--(-1,0)node[below]{$E_{n}$};   
    \node (tenten) at (-1.5,0.5) {$\cdots$}; 
    \draw (-2,1)--(-3,0)node[below]{$E_{2}$};  
    \draw (-2.5,0)--(-3.5,1)node[above]{$E_{1}$}; 
    \draw (-3,1)--(-4,0)node[below]{$\overline{L_{\infty}}$};  
    \node (Lselfint) at (-4,-0.8) {$\Circled{-1}$}; 
    \node (E1selfint) at (-3.5,1.7) {$\Circled{-2}$}; 
    \node (E2selfint) at (-3,-0.7) {$\Circled{-2}$}; 
    \node (Enselfint) at (-1,-0.7) {$\Circled{-3}$}; 
    \node (En+2selfint) at (0.5,-0.7) {$\Circled{-1}$}; 
    \node (En+1selfint) at (1,1.7) {$\Circled{-2}$}; 
    \node (Mselfint) at (2,-0.6) {$\Circled{-(n+1)}$}; 
    \node (intersection) at (4,0) {$(\overline{C} \cdot E_{n+2})=2$}; 
  \end{tikzpicture} 
  \caption*{Figure \ref{sec:pf.of.main.th.2}-3}  
\end{figure} 
Then, we have 
\begin{align*}
  \overline{\sigma}^{*}A 
  &= \overline{\sigma}^{*}C - 5\overline{\sigma}^{*}M_{\infty} 
  -d\overline{\sigma}^{*}L_{\infty} \\
  &= (\overline{C}+5(E_{1}+2E_{2}+ \dots + nE_{n})
  +(5n+3)E_{n+1}+(10n+5)E_{n+2}+E_{\mathrm{fin}}) \\ 
  & -5(\overline{M_{\infty}}+E_{1}+2E_{2}+ \dots +nE_{n}
  +(n+1)E_{n+1}+(2n+1)E_{n+2})  \\ 
  &-d(\overline{L_{\infty}}+E_{1}+E_{2}+ \dots +E_{n}
  +E_{n+1}+2E_{n+2}) \\ 
  &= \overline{C}-5\overline{M_{\infty}}-d\overline{L_{\infty}}
  -d(E_{1}+E_{2}+ \dots +E_{n}) \\
  &-(d+2)E_{n+1}-2dE_{n+2}  
  + E_{\mathrm{fin}},  
\end{align*}
and 
\begin{align*}
  K_{\overline{\Sigma}} 
  &\sim -2\overline{\sigma}^{*}M_{\infty} 
  -2\overline{\sigma}^{*}L_{\infty} +E_{1}+2E_{2}+ \dots +nE_{n} \\
  &+(n+1)E_{n+1}+(2n+2)E_{n+2}  
  +E_{\mathrm{can}} \\ 
  &= -2\overline{M_{\infty}}-2\overline{L_{\infty}}
  -(3E_{1}+4E_{2}+ \dots + (n+2)E_{n}) \\
  &-(n+3)E_{n+1}-(2n+4)E_{n+2}+E_{\mathrm{can}}.  
\end{align*} 

\medskip 

\paragraph*{Subcase 3-1 \texorpdfstring{$d=10m+3$ $(n=2m)$}{d=10m+3}}\quad\\ 

We have 
\begin{align*}
  \overline{B} 
  &= \overline{C}+\overline{M_{\infty}}+\overline{L_{\infty}}
  +E_{1}+ \dots + E_{2m}+E_{2m+1}+E_{\mathrm{fin}}^{(1)}, \\
  \overline{Z} 
  &= 3\overline{M_{\infty}}+(5m+2)\overline{L_{\infty}}
  +(5m+2)(E_{1}+ \dots + E_{2m}) \\
  &+(5m+3)E_{2m+1}+(10m+3)E_{2m+2}-E_{\mathrm{fin}}^{(2)}, \\
  K_{\overline{\Sigma}}+\overline{Z}
  &\sim \overline{M_{\infty}}+5m\overline{L_{\infty}} \\
  &+(5m-1)E_{1}+(5m-2)E_{2}+ \dots + 3mE_{2m} \\
  &+3mE_{2m+1}+(6m-1)E_{2m+2} 
  +(E_{\mathrm{can}}-E_{\mathrm{fin}}^{(2)}), \\
  p_{a}(\overline{Z}) 
  &= 4m+\mu, \\
  (K_{\overline{\Sigma}}+\overline{Z})^{2} 
  &= 2m-2+\nu. 
\end{align*}
Reflecting the data about the branch locus, 
we have the configuration as in Figure \ref{sec:pf.of.main.th.2}-3-1A.  
\begin{figure}[H]
  \centering
  \begin{tikzpicture}
    \draw (0.4,0.8) to [out=225,in=135] (0.25,0.25) 
    to [out=315,in=225] (0.7,0.2)node[right]{$\overline{C}$}; 
    \draw (1,1)node[above]{$E_{2m+1}$}--(0,0);  
    \draw (0.5,1)--(1.5,0)node[right]{$\overline{M_{\infty}}$};
    \draw[dashed] (0.5,0)node[below]{$E_{2m+2}$}--(-0.5,1); 
    \draw (0,1)--(-1,0)node[below]{$E_{2m}$};   
    \node (tenten) at (-1.5,0.5) {$\cdots$}; 
    \draw (-2,1)--(-3,0)node[below]{$E_{2}$};  
    \draw (-2.5,0)--(-3.5,1)node[above]{$E_{1}$}; 
    \draw (-3,1)--(-4,0)node[below]{$\overline{L_{\infty}}$};  
    \node (Lselfint) at (-4,-0.8) {$\Circled{-1}$}; 
    \node (E1selfint) at (-3.5,1.7) {$\Circled{-2}$}; 
    \node (E2selfint) at (-3,-0.7) {$\Circled{-2}$}; 
    \node (Enselfint) at (-1,-0.7) {$\Circled{-3}$}; 
    \node (En+2selfint) at (0.5,-0.7) {$\Circled{-1}$}; 
    \node (En+1selfint) at (1,1.7) {$\Circled{-2}$}; 
    \node (Mselfint) at (2,-0.6) {$\Circled{-(2m+1)}$}; 
    \node (intersection) at (4,0) {$(\overline{C} \cdot E_{2m+2})=2$}; 
  \end{tikzpicture} 
  \caption*{Figure \ref{sec:pf.of.main.th.2}-3-1A}  
\end{figure} 
After some blowing-ups, 
we have the dual graph on $S$ by taking a double covering of $\Sigma$ as in Figure \ref{sec:pf.of.main.th.2}-3-1B,   
\begin{figure}[H] 
  \centering
  \scalebox{0.8}{
  \begin{tikzpicture}
    \node[draw,shape=circle,inner sep=2pt,fill=black] (linf) at (-4,0) {};
    \node[draw,shape=circle,inner sep=2pt] (llinf) at (-3,0) {};
    \node[draw,shape=circle,inner sep=2pt,fill=black] (e1) at (-2,0) {};  
    \node[draw,shape=circle,inner sep=2pt] (ee1) at (-1,0) {}; 
    \node (ueyoko) at (0,0) {$\cdots$}; 
    \node[draw,shape=circle,inner sep=2pt] (e2m-1) at (1,0) {};
    \node[draw,shape=circle,inner sep=2pt,fill=black] (e2m) at (2,0) {};
    \node[draw,shape=circle,inner sep=2pt] (e2m+2) at (3,0) {};
    \node[draw,shape=circle,inner sep=2pt] (l) at (4,0) {}; 
    \node[draw,shape=circle,inner sep=2pt,fill=black] (e2m+1) at (5,0) {};
    \node[draw,shape=circle,inner sep=2pt] (ll) at (6,0) {}; 
    \node[draw,shape=circle,inner sep=2pt,fill=black] (m) at (7,0) {};     
    \draw (linf)--(llinf)--(e1)--(ee1);     
    \draw (ee1)--(-0.5,0);
    \draw (0.5,0)--(e2m-1);    
    \draw (e2m-1)--(e2m)--(e2m+2);
    \draw[double distance=2pt,nfold=2] (e2m+2)--(l);  
    \draw (l)--(e2m+1)--(ll)--(m);  
    \node (linfname) at (-4,-0.5) {$\widetilde{L_{\infty}}$}; 
    \node (llinfname) at (-3,-0.5) {$\widetilde{L_{\infty}'}$}; 
    \node (e1name) at (-2,-0.5) {$\widetilde{E_{1}}$};   
    \node (ee1name) at (-1,-0.5) {$\widetilde{E_{1}'}$}; 
    \node (e2m-2name) at (1,-0.5) {$\widetilde{E_{2m-1}}$};  
    \node (ee2mname) at (2,-0.5) {$\widetilde{E_{2m}}$}; 
    \node (e2m+2mname) at (3,-0.5) {$\widetilde{E_{2m+2}}$}; 
    \node (lmname) at (4,-0.5) {$L$};   
    \node (e2m+1name) at (5,-0.5) {$\widetilde{E_{2m+1}}$}; 
    \node (llname) at (6,-0.5) {$L'$};  
    \node (mname) at (7,-0.5) {$\widetilde{M_{\infty}}$};  
    \node (linfselfint) at (-4,0.5) {$-1$}; 
    \node (llinfselfint) at (-3,0.5) {$-2$};   
    \node (e1selfint) at (-2,0.5) {$-2$}; 
    \node (ee1selfint) at (-1,0.5) {$-2$}; 
    \node (e2m-2selfint) at (1,0.5) {$-2$};   
    \node (ee2selfint) at (2,0.5) {$-2$}; 
    \node (e2m+2selfint) at (3,0.5) {$-4$};  
    \node (lselfint) at (4,0.5) {$-2$}; 
    \node (e2m+1selfint) at (5,0.5) {$-2$}; 
    \node (llselfint) at (6,0.5) {$-2$}; 
    \node (mselfint) at (7,0.5) {$-(m+1)$};  
  \end{tikzpicture}  
  } 
  \caption*{Figure \ref{sec:pf.of.main.th.2}-3-1B}   
\end{figure} 
\noindent where 
$\pi^{*}(\overline{L_{\infty}})=2\widetilde{L_{\infty}}+\widetilde{L_{\infty}'}$, 
$\pi^{*}(E_{1})=\widetilde{L_{\infty}'}+2\widetilde{E_{1}}+\widetilde{E_{1}'}$, 
$\pi^{*}(E_{i})=\widetilde{E_{i-1}'}+2\widetilde{E_{i}}+\widetilde{E_{i}'}$ for $2 \le i \le 2m-1$, 
$\pi^{*}(E_{2m})=\widetilde{E_{2m-1}'}+2\widetilde{E_{2m}}$, 
$\pi^{*}(E_{2m+2})=\widetilde{E_{2m+2}}+L$, 
$\pi^{*}(E_{2m+1})=L+2\widetilde{E_{2m+1}}+L'$
and 
$\pi^{*}(\overline{M_{\infty}})=2\widetilde{M_{\infty}}+L'$.

By contracting $\widetilde{L_{\infty}}$, $\widetilde{E_{1}}$, $\dots$, 
$\widetilde{E_{2m-1}}$, $\widetilde{E_{2m}}$ 
and 
other $(-1)$-curves contained in reducible fibers in the affine part,        
we have the relatively minimal model $X$. 
By relabeling suitably, we have the dual graph on $X$ as in Figure \ref{sec:pf.of.main.th.2}-3-1C.   
\begin{figure}[H]
  \centering
  \scalebox{0.8}{
    \begin{tikzpicture} 
      \node[draw,shape=circle] (0) at (1,0) {};
      \node[draw,shape=circle] (a0) at (0,0) {};
      \node[draw,shape=circle] (a1) at (-1,0) {};
      \node[draw,shape=circle] (a2) at (-2,0) {}; 
      \node[draw,shape=circle] (a3) at (-3,0) {}; 
      \draw (0)--(a0)--(a1)--(a2);
      \draw[double distance=4pt] (a2)--(a3); 
      \node (0name) at (1,-0.5) {$(O)$};
      \node (a0name) at (0,-0.5) {$A_{0}$};
      \node (a1name) at (-1,-0.5) {$A_{1}$};
      \node (a2name) at (-2,-0.5) {$A_{2}$};
      \node (a3name) at (-3,-0.5) {$A_{3}$};
      \node (0selfint) at (1,0.5) {$-(m+1)$};
      \node (a3selfint) at (-3,0.5) {$-3$};
    \end{tikzpicture}
  } 
  \caption*{Figure \ref{sec:pf.of.main.th.2}-3-1C}   
\end{figure}
Then, we obtain an explicit expression of 
the canonical divisor of $X$. 
\[
  K_{X} \sim 
  2(O)+(3m-1)F+2(A_{0}+A_{1}+A_{2})+A_{3}
  +\pi^{*}(E_{\mathrm{can}}-E_{\mathrm{fin}}^{(2)}),   
\]
where $F$ is a general fiber. 

\medskip 

\paragraph*{Subcase 3-2 \texorpdfstring{$d=10m+8$ $(n=2m+1)$}{d=10m+8}}\quad\\ 

We have 
\begin{align*}
  \overline{B} 
  &= \overline{C}+\overline{M_{\infty}}+E_{\mathrm{fin}}^{(1)}, \\
  \overline{Z} 
  &= 3\overline{M_{\infty}}+(5m+4)\overline{L_{\infty}}
  +(5m+4)(E_{1}+ \dots + E_{2m+1}) \\
  &+(5m+5)E_{2m+2}+(10m+8)E_{2m+3}-E_{\mathrm{fin}}^{(2)}, \\
  K_{\overline{\Sigma}}+\overline{Z}
  &\sim \overline{M_{\infty}}+(5m+2)\overline{L_{\infty}} \\
  &+(5m+1)E_{1}+5mE_{2}+ \dots + (3m+2)E_{2m+1} \\
  &+(3m+1)E_{2m+3}+(6m+2)E_{2m+3}  
  +(E_{\mathrm{can}}-E_{\mathrm{fin}}^{(2)}), \\
  p_{a}(\overline{Z}) 
  &= 4m+2+\mu, \\
  (K_{\overline{\Sigma}}+\overline{Z})^{2} 
  &= 2m-1+\nu. 
\end{align*}
Reflecting the data about the branch locus, 
we have the configuration as in Figure \ref{sec:pf.of.main.th.2}-3-2A.   
\begin{figure}[H]
  \centering
  \begin{tikzpicture}
    \draw (0.4,0.8) to [out=225,in=135] (0.25,0.25) 
    to [out=315,in=225] (0.7,0.2)node[right]{$\overline{C}$}; 
    \draw[dashed] (1,1)node[above]{$E_{2m+2}$}--(0,0);  
    \draw (0.5,1)--(1.5,0)node[right]{$\overline{M_{\infty}}$};
    \draw[dashed] (0.5,0)node[below]{$E_{2m+3}$}--(-0.5,1); 
    \draw[dashed] (0,1)--(-1,0)node[below]{$E_{2m+1}$};   
    \node (tenten) at (-1.5,0.5) {$\cdots$}; 
    \draw[dashed] (-2,1)--(-3,0)node[below]{$E_{2}$};  
    \draw[dashed] (-2.5,0)--(-3.5,1)node[above]{$E_{1}$}; 
    \draw[dashed] (-3,1)--(-4,0)node[below]{$\overline{L_{\infty}}$};  
    \node (Lselfint) at (-4,-0.8) {$\Circled{-1}$}; 
    \node (E1selfint) at (-3.5,1.7) {$\Circled{-2}$}; 
    \node (E2selfint) at (-3,-0.7) {$\Circled{-2}$}; 
    \node (Enselfint) at (-1,-0.7) {$\Circled{-3}$}; 
    \node (En+2selfint) at (0.5,-0.7) {$\Circled{-1}$}; 
    \node (En+1selfint) at (1,1.7) {$\Circled{-2}$}; 
    \node (Mselfint) at (2,-0.6) {$\Circled{-(2m+2)}$}; 
    \node (intersection) at (4,0) {$(\overline{C} \cdot E_{2m+3})=2$}; 
  \end{tikzpicture} 
  \caption*{Figure \ref{sec:pf.of.main.th.2}-3-2A}  
\end{figure} 
We have the dual graph on $S$ by taking a double covering of $\Sigma$ as in Figure \ref{sec:pf.of.main.th.2}-3-2B,    
\begin{figure}[H] 
  \centering
  \scalebox{0.8}{
  \begin{tikzpicture}
    \node[draw,shape=circle,inner sep=2pt] (linf) at (2,1) {};
    \node[draw,shape=circle,inner sep=2pt] (llinf) at (2,-1) {};
    \node[draw,shape=circle,inner sep=2pt] (e1) at (3,1) {};
    \node[draw,shape=circle,inner sep=2pt] (ee1) at (3,-1) {}; 
    \node (ueyoko) at (4,1) {$\cdots$}; 
    \node (shitayoko) at (4,-1) {$\cdots$}; 
    \node[draw,shape=circle,inner sep=2pt] (e2m) at (5,1) {};
    \node[draw,shape=circle,inner sep=2pt] (ee2m) at (5,-1) {};
    \node[draw,shape=circle,inner sep=2pt] (e2m+1) at (6,1) {};
    \node[draw,shape=circle,inner sep=2pt] (ee2m+1) at (6,-1) {};
    \node[draw,shape=circle,inner sep=2pt] (e2m+3) at (7,1) {};
    \node[draw,shape=circle,inner sep=2pt] (ee2m+3) at (7,-1) {};
    \node[draw,shape=circle,inner sep=2pt] (e2m+2) at (8,0) {};
    \node[draw,shape=circle,inner sep=2pt,fill=black] (m) at (9,0) {};   
    \draw (linf)--(e1)--($(e1)+(0.5,0)$); 
    \draw (llinf)--(ee1)--($(ee1)+(0.5,0)$);  
    \draw ($(e2m)+(-0.5,0)$)--(e2m);  
    \draw ($(ee2m)+(-0.5,0)$)--(ee2m);   
    \draw (e2m)--(e2m+1)--(e2m+3)--(e2m+2); 
    \draw (ee2m)--(ee2m+1)--(ee2m+3)--(e2m+2);   
    \draw (e2m+2)--(m);   
    \draw (e2m+3)--(ee2m+3);  
    \node (linfname) at ($(linf)+(0,-0.5)$) {$\widetilde{L_{\infty}}$}; 
    \node (llinfname) at ($(llinf)+(0,-0.5)$) {$\widetilde{L_{\infty}'}$}; 
    \node (e1name) at ($(e1)+(0,-0.5)$) {$\widetilde{E_{1}}$};  
    \node (ee1name) at ($(ee1)+(0,-0.5)$) {$\widetilde{E_{1}'}$}; 
    \node (e2mname) at ($(e2m)+(0,-0.5)$) {$\widetilde{E_{2m}}$};  
    \node (ee2mname) at ($(ee2m)+(0,-0.5)$) {$\widetilde{E_{2m}'}$}; 
    \node (e2m+1name) at ($(e2m+1)+(0,-0.5)$) {$\widetilde{E_{2m+1}}$};  
    \node (ee2m+1name) at ($(ee2m+1)+(0,-0.5)$) {$\widetilde{E_{2m+1}'}$};
    \node (e2m+3name) at ($(e2m+3)+(1,0.5)$) {$\widetilde{E_{2m+3}}$};  
    \node (ee2m+3name) at ($(ee2m+3)+(1,-0.5)$) {$\widetilde{E_{2m+3}'}$};   
    \node (e2m+2name) at ($(e2m+2)+(0.5,-0.5)$) {$\widetilde{E_{2m+2}}$}; 
    \node (mname) at ($(m)+(0.5,-0.5)$) {$\widetilde{M_{\infty}}$};   
    \node (linfselfint) at ($(linf)+(0,0.5)$) {$-1$}; 
    \node (llinfselfint) at ($(llinf)+(0,0.5)$) {$-1$}; 
    \node (e1selfint) at ($(e1)+(0,0.5)$) {$-2$};  
    \node (ee1selfint) at ($(ee1)+(0,0.5)$) {$-2$};  
    \node (e2mselfint) at ($(e2m)+(0,0.5)$) {$-2$};  
    \node (ee2mselfint) at ($(ee2m)+(0,0.5)$) {$-2$};  
    \node (e2m+1selfint) at ($(e2m+1)+(0,0.5)$) {$-3$};  
    \node (ee2m+1selfint) at ($(ee2m+1)+(0,0.5)$) {$-3$}; 
    \node (e2m+3selfint) at ($(e2m+3)+(0,0.5)$) {$-2$};  
    \node (ee2m+3selfint) at ($(ee2m+3)+(0,-0.5)$) {$-2$};  
    \node (e2m+2selfint) at ($(e2m+2)+(0,0.5)$) {$-4$}; 
    \node (mselfint) at ($(m)+(0,0.5)$) {$-(m+1)$};   
  \end{tikzpicture}  
  } 
  \caption*{Figure \ref{sec:pf.of.main.th.2}-3-2B}   
\end{figure} 
\noindent where 
$\pi^{*}(\overline{L_{\infty}})=\widetilde{L_{\infty}}+\widetilde{L_{\infty}'}$, 
$\pi^{*}(E_{i})=\widetilde{E_{i}}+\widetilde{E_{i}'}$ for $1 \le i \le 2m+1$ or $i=2m+3$, 
$\pi^{E_{2m+2}}=\widetilde{E_{2m+2}}$ 
and 
$\pi^{*}(\overline{M_{\infty}})=2\widetilde{M_{\infty}}$. 
Note that three components $\widetilde{E_{2m+3}}, \widetilde{E_{2m+3}'}$ and 
$\widetilde{E_{2m+2}}$ meet each other in one point.

By contracting 
$\widetilde{L_{\infty}}$, $\widetilde{E_{1}}$, $\dots$, $\widetilde{E_{2m}}$ 
and 
$\widetilde{L_{\infty}'}$, $\widetilde{E_{1}'}$, $\dots$, $\widetilde{E_{2m}'}$, 
and 
other $(-1)$-curves contained in reducible fibers in the affine part,        
we have the relatively minimal model $X$. 
By relabeling suitably, we have the dual graph on $X$ as in Figure \ref{sec:pf.of.main.th.2}-3-2C.      
\begin{figure}[H]
  \centering
  \scalebox{0.8}{
    \begin{tikzpicture} 
      \node[draw,shape=circle] (0) at (1,0) {};
      \node[draw,shape=circle] (a0) at (0,0) {};
      \node[draw,shape=circle] (a1) at (-1,1) {};
      \node[draw,shape=circle] (a2) at (-2,1) {}; 
      \node[draw,shape=circle] (a3) at (-1,-1) {};
      \node[draw,shape=circle] (a4) at (-2,-1) {}; 
      \draw (0)--(a0)--(a1)--(a2);
      \draw (a0)--(a3)--(a4);  
      \draw (a1)--(a3); 
      \node (0name) at (1,-0.5) {$(O)$};
      \node (a0name) at (0,-0.5) {$A_{0}$};
      \node (a1name) at (-1.3,0.5) {$A_{1}$};
      \node (a2name) at (-2,0.5) {$A_{2}$};
      \node (a3name) at (-1.3,-1.5) {$A_{3}$};
      \node (a4name) at (-2,-1.5) {$A_{4}$};
      \node (0selfint) at (1,0.5) {$-(m+1)$};
      \node (a0selfint) at (0,0.5) {$-4$};
    \end{tikzpicture}
  } 
  \caption*{Figure \ref{sec:pf.of.main.th.2}-3-2C}   
\end{figure}
Then, we obtain an explicit expression of 
the canonical divisor of $X$.
\[
  K_{X} \sim 
  2(O)+(3m+1)F
  +\pi^{*}(E_{\mathrm{can}}-E_{\mathrm{fin}}^{(2)}),   
\]
where $F$ is a general fiber. 

\subsubsection*{Case \texorpdfstring{$4$}{4} \texorpdfstring{$d=5n+4$ $(n>0)$}{d=5n+4}}  

We have to blow-up $(n+1)$ times to obtain the nonsingular curve $\overline{C}$ (Figure \ref{sec:pf.of.main.th.2}-4). 
\begin{figure}[H]
  \centering
  \begin{tikzpicture}
    \draw (-0.5,0.3)node[above]{$\overline{C}$} 
    to [out=315,in=225] (0.25,0.25) 
    to [out=45,in=315] (0.3,0.8); 
    \node (intersection) at (4,0) {$(\overline{C} \cdot E_{n+1})=4$}; 
    \draw (1,1)node[above]{$E_{n+1}$}--(0,0);  
    \draw (0.5,1)--(1.5,0)node[right]{$\overline{M_{\infty}}$}; 
    \draw (-0.5,1)--(0.5,0)node[below]{$E_{n}$};   
    \node (tenten) at (-1,0.5) {$\cdots$}; 
    \draw (-1.5,1)--(-2.5,0)node[below]{$E_{2}$};  
    \draw (-2,0)--(-3,1)node[above]{$E_{1}$}; 
    \draw (-2.5,1)--(-3.5,0)node[below]{$\overline{L_{\infty}}$};  
    \node (Lselfint) at (-3.5,-0.8) {$\Circled{-1}$}; 
    \node (E1selfint) at (-3,1.7) {$\Circled{-2}$}; 
    \node (E2selfint) at (-2.5,-0.7) {$\Circled{-2}$}; 
    \node (Enselfint) at (0.5,-0.7) {$\Circled{-2}$}; 
    \node (En+1selfint) at (1,1.7) {$\Circled{-1}$}; 
    \node (Mselfint) at (2,-0.6) {$\Circled{-(n+1)}$}; 
  \end{tikzpicture} 
  \caption*{Figure \ref{sec:pf.of.main.th.2}-4}  
\end{figure} 
Then, we have 
\begin{align*} 
  \overline{\sigma}^{*}A 
  &= \overline{\sigma}^{*}C - 5\overline{\sigma}^{*}M_{\infty} 
  -d\overline{\sigma}^{*}L_{\infty} \\
  &= (\overline{C}+5(E_{1}+2E_{2}+ \dots + nE_{n})
  +(5n+4)E_{n+1}+E_{\mathrm{fin}}) \\ 
  & -5(\overline{M_{\infty}}+E_{1}+2E_{2}+ \dots +nE_{n}
  +(n+1)E_{n+1})  \\ 
  &-d(\overline{L_{\infty}}+E_{1}+E_{2}+ \dots +E_{n}+E_{n+1}) \\ 
  &= \overline{C}-5\overline{M_{\infty}}-d\overline{L_{\infty}}
  -d(E_{1}+E_{2}+ \dots E_{n})-(d+1)E_{n+1}+ E_{\mathrm{fin}},  
\end{align*} 
and 
\begin{align*}
  K_{\overline{\Sigma}} 
  &\sim -2\overline{\sigma}^{*}M_{\infty} 
  -2\overline{\sigma}^{*}L_{\infty} +(E_{1}+2E_{2}+ \dots nE_{n})
  +(n+1)E_{n+1}+E_{\mathrm{can}} \\ 
  &= -2\overline{M_{\infty}}-2\overline{L_{\infty}}
  -(3E_{1}+4E_{2}+ \dots + (n+2)E_{n})-(n+3)E_{n+1}+E_{\mathrm{can}}.  
\end{align*} 

\medskip

\paragraph*{Subcase 4-1 \texorpdfstring{$d=10m+4$ $(n=2m)$}{d=10m+4}}\quad\\ 

We have 
\begin{align*}
  \overline{B} 
  &= \overline{C}+\overline{M_{\infty}}
  +E_{2m+1}+E_{\mathrm{fin}}^{(1)}, \\
  \overline{Z} 
  &= 3\overline{M_{\infty}}+(5m+2)\overline{L_{\infty}}
  +(5m+2)(E_{1}+ \dots + E_{2m}) \\
  &+(5m+3)E_{2m+1}-E_{\mathrm{fin}}^{(2)}, \\
  K_{\overline{\Sigma}}+\overline{Z}
  &\sim \overline{M_{\infty}}+5m\overline{L_{\infty}} \\
  &+(5m-1)E_{1}+(5m-2)E_{2}+ \dots +3mE_{2m} \\ 
  &+3mE_{2m+1}+(E_{\mathrm{can}}-E_{\mathrm{fin}}^{(2)}), \\
  p_{a}(\overline{Z}) 
  &= 4m+1+\mu, \\
  (K_{\overline{\Sigma}}+\overline{Z})^{2} 
  &= 2m-1+\nu. 
\end{align*}
Reflecting the data about the branch locus, 
we have the configuration as in Figure \ref{sec:pf.of.main.th.2}-4-1A.   
\begin{figure}[H]
  \centering
  \begin{tikzpicture}
    \draw (-0.5,0.3)node[above]{$\overline{C}$} 
    to [out=315,in=225] (0.25,0.25) 
    to [out=45,in=315] (0.3,0.8); 
    \node (intersection) at (4,0) {$(\overline{C} \cdot E_{2m+1})=4$}; 
    \draw (1,1)node[above]{$E_{2m+1}$}--(0,0);  
    \draw (0.5,1)--(1.5,0)node[right]{$\overline{M_{\infty}}$}; 
    \draw[dashed] (-0.5,1)--(0.5,0)node[below]{$E_{2m}$};   
    \node (tenten) at (-1,0.5) {$\cdots$}; 
    \draw[dashed] (-1.5,1)--(-2.5,0)node[below]{$E_{2}$};  
    \draw[dashed] (-2,0)--(-3,1)node[above]{$E_{1}$}; 
    \draw[dashed] (-2.5,1)--(-3.5,0)node[below]{$\overline{L_{\infty}}$};  
    \node (Lselfint) at (-3.5,-0.8) {$\Circled{-1}$}; 
    \node (E1selfint) at (-3,1.7) {$\Circled{-2}$}; 
    \node (E2selfint) at (-2.5,-0.7) {$\Circled{-2}$}; 
    \node (Enselfint) at (0.5,-0.7) {$\Circled{-2}$}; 
    \node (En+1selfint) at (1,1.7) {$\Circled{-1}$}; 
    \node (Mselfint) at (2,-0.6) {$\Circled{-(2m+1)}$}; 
  \end{tikzpicture} 
  \caption*{Figure \ref{sec:pf.of.main.th.2}-4-1A}  
\end{figure} 
After some blowing-ups, 
we have the dual graph on $S$ by taking a double covering of $\Sigma$ as in Figure \ref{sec:pf.of.main.th.2}-4-1B,  
\begin{figure}[H] 
  \centering
  \scalebox{0.8}{
  \begin{tikzpicture}
    \node[draw,shape=circle,inner sep=2pt] (linf) at (0,1) {};
    \node[draw,shape=circle,inner sep=2pt] (llinf) at (0,-1) {};
    \node[draw,shape=circle,inner sep=2pt] (e1) at (1,1) {};
    \node[draw,shape=circle,inner sep=2pt] (ee1) at (1,-1) {}; 
    \node (ueyoko) at (2,1) {$\cdots$}; 
    \node (shitayoko) at (2,-1) {$\cdots$}; 
    \node[draw,shape=circle,inner sep=2pt] (e2m-1) at (3,1) {};
    \node[draw,shape=circle,inner sep=2pt] (ee2m-1) at (3,-1) {};
    \node[draw,shape=circle,inner sep=2pt] (e2m) at (4,1) {};
    \node[draw,shape=circle,inner sep=2pt] (ee2m) at (4,-1) {}; 
    \node[draw,shape=circle,inner sep=2pt] (l1) at (5,1) {};
    \node[draw,shape=circle,inner sep=2pt] (ll1) at (5,-1) {};
    \node[draw,shape=circle,inner sep=2pt] (l2) at (6,1) {};
    \node[draw,shape=circle,inner sep=2pt] (ll2) at (6,-1) {};
    \node[draw,shape=circle,inner sep=2pt] (l3) at (7,1) {};
    \node[draw,shape=circle,inner sep=2pt] (ll3) at (7,-1) {};
    \node[draw,shape=circle,inner sep=2pt] (l4) at (8,0) {};
    \node[draw,shape=circle,inner sep=2pt,fill=black] (e2m+1) at (9,0) {}; 
    \node[draw,shape=circle,inner sep=2pt] (l5) at (10,0) {}; 
    \node[draw,shape=circle,inner sep=2pt,fill=black] (m) at (11,0) {};   
    \draw (linf)--(e1)--(1.5,1); 
    \draw (llinf)--(ee1)--(1.5,-1);  
    \draw (2.5,1)--(e2m-1);  
    \draw (2.5,-1)--(ee2m-1);  
    \draw (e2m-1)--(e2m)--(l1)--(l2)--(l3)--(l4); 
    \draw (ee2m-1)--(ee2m)--(ll1)--(ll2)--(ll3)--(l4);  
    \draw (l4)--(e2m+1)--(l5)--(m);  
    \node (linfname) at ($(linf)+(0,-0.5)$) {$\widetilde{L_{\infty}}$}; 
    \node (llinfname) at ($(llinf)+(0,-0.5)$) {$\widetilde{L_{\infty}'}$}; 
    \node (e1name) at ($(e1)+(0,-0.5)$) {$\widetilde{E_{1}}$};  
    \node (ee1name) at ($(ee1)+(0,-0.5)$) {$\widetilde{E_{1}'}$};
    \node (e2m-1name) at ($(e2m-1)+(0,-0.5)$) {$\widetilde{E_{2m}}$};  
    \node (ee2m-1name) at ($(ee2m-1)+(0,-0.5)$) {$\widetilde{E_{2m}'}$};  
    \node (e2mname) at ($(e2m)+(0,-0.5)$) {$\widetilde{E_{2m}}$};  
    \node (ee2mname) at ($(ee2m)+(0,-0.5)$) {$\widetilde{E_{2m}'}$}; 
    \node (l1name) at ($(l1)+(0,-0.5)$) {$L_{1}$};
    \node (ll1name) at ($(ll1)+(0,-0.5)$) {$L_{1}'$}; 
    \node (l2name) at ($(l2)+(0,-0.5)$) {$L_{2}$};  
    \node (ll2name) at ($(ll2)+(0,-0.5)$) {$L_{2}'$};  
    \node (l3name) at ($(l3)+(0,-0.5)$){$L_{3}$}; 
    \node (ll3name) at ($(ll3)+(0,-0.5)$){$L_{3}'$}; 
    \node (l4name) at ($(l4)+(0,-0.5)$) {$L_{4}$}; 
    \node (e2m+1name) at ($(e2m+1)+(0,-0.5)$) {$\widetilde{E_{2m+1}}$}; 
    \node (l5name) at ($(l5)+(0,-0.5)$) {$L_{5}$};
    \node (mname) at ($(m)+(0,-0.5)$) {$\widetilde{M_{\infty}}$};   
    \node (linfselfint) at ($(linf)+(0,0.5)$) {$-1$}; 
    \node (llinfselfint) at ($(llinf)+(0,0.5)$) {$-1$}; 
    \node (e1selfint) at ($(e1)+(0,0.5)$) {$-2$};  
    \node (ee1selfint) at ($(ee1)+(0,0.5)$) {$-2$}; 
    \node (e2m-1selfint) at ($(e2m-1)+(0,0.5)$) {$-2$};  
    \node (ee2m-1selfint) at ($(ee2m-1)+(0,0.5)$) {$-2$}; 
    \node (e2mselfint) at ($(e2m)+(0,0.5)$) {$-3$};  
    \node (ee2mselfint) at ($(ee2m)+(0,0.5)$) {$-3$}; 
    \node (l1selfint) at ($(l1)+(0,0.5)$) {$-2$}; 
    \node (ll1selfint) at ($(ll1)+(0,0.5)$) {$-2$};  
    \node (l2selfint) at ($(l2)+(0,0.5)$) {$-2$}; 
    \node (ll2selfint) at ($(ll2)+(0,0.5)$) {$-2$}; 
    \node (l3selfint) at ($(l3)+(0,0.5)$) {$-2$}; 
    \node (ll3selfint) at ($(ll3)+(0,0.5)$) {$-2$}; 
    \node (l4selfint) at ($(l4)+(0,0.5)$) {$-2$}; 
    \node (e2m+1selfint) at ($(e2m+1)+(0,0.5)$) {$-3$};  
    \node (l5selfint) at ($(l5)+(0,0.5)$) {$-2$}; 
    \node (mselfint) at ($(m)+(0,0.5)$) {$-(m+1)$};   
  \end{tikzpicture}  
  } 
  \caption*{Figure \ref{sec:pf.of.main.th.2}-4-1B}  
\end{figure} 
\noindent where 
$\pi^{*}(\overline{L_{\infty}})=\widetilde{L_{\infty}}+\widetilde{L_{\infty}'}$,  
$\pi^{*}(E_{i})=\widetilde{E_{i}}+\widetilde{E_{i}'}$ for $1 \le i \le 2m-1$, 
$\pi^{*}(E_{2m})=\widetilde{E_{2m}}+\widetilde{E_{2m}'}+\sum_{i=1}^{3}(L_{i}+L_{i}')+L_{4}$, 
$\pi^{*}(E_{2m+1})=\sum_{i=1}^{3}i(L_{i}+L_{i}')+4L_{4}+2\widetilde{E_{2m+1}}+L_{5}$ 
and 
$\pi^{*}(\overline{M_{\infty}})=2\widetilde{M_{\infty}}+L_{5}$.

By contracting 
$\widetilde{L_{\infty}}$, $\widetilde{E_{1}}$, $\dots$, $\widetilde{E_{2m}}$ 
and 
$\widetilde{L_{\infty}'}$, $\widetilde{E_{1}'}$, $\dots$, $\widetilde{E_{2m}'}$ 
and 
other $(-1)$-curves contained in reducible fibers in the affine part,        
we have the relatively minimal model $X$.   
By relabeling suitably, we have the dual graph on $X$ as in Figure \ref{sec:pf.of.main.th.2}-4-1C.     
\begin{figure}[H]
  \centering
  \scalebox{0.8}{
    \begin{tikzpicture} 
      \node[draw,shape=circle] (0) at (1,0) {};
      \node[draw,shape=circle] (a0) at (0,0) {};
      \node[draw,shape=circle] (a1) at (-1,0) {};
      \node[draw,shape=circle] (a2) at (-2,0) {}; 
      \node[draw,shape=circle] (a3) at (-3,1) {};
      \node[draw,shape=circle] (a4) at (-4,1) {};
      \node[draw,shape=circle] (a5) at (-5,1) {}; 
      \node[draw,shape=circle] (a6) at (-6,1) {}; 
      \node[draw,shape=circle] (a7) at (-3,-1) {};
      \node[draw,shape=circle] (a8) at (-4,-1) {};
      \node[draw,shape=circle] (a9) at (-5,-1) {}; 
      \node[draw,shape=circle] (a10) at (-6,-1) {}; 
      \draw (0)--(a0)--(a1)--(a2); 
      \draw (a2)--(a3)--(a4)--(a5)--(a6); 
      \draw (a2)--(a7)--(a8)--(a9)--(a10);  
      \node (0name) at (1,-0.5) {$(O)$}; 
      \node (a0name) at (0,-0.5) {$A_{0}$};
      \node (a1name) at (-1,-0.5) {$A_{1}$};
      \node (a2name) at (-2,-0.5) {$A_{2}$}; 
      \node (a3name) at (-3,0.5) {$A_{3}$}; 
      \node (a4name) at (-4,0.5) {$A_{4}$}; 
      \node (a5name) at (-5,0.5) {$A_{5}$}; 
      \node (a6name) at (-6,0.5) {$A_{6}$};  
      \node (a7name) at (-3,-1.5) {$A_{7}$}; 
      \node (a8name) at (-4,-1.5) {$A_{8}$}; 
      \node (a9name) at (-5,-1.5) {$A_{9}$}; 
      \node (a10name) at (-6,-1.5) {$A_{10}$};  
      \node (0selfint) at (1,0.5) {$-(m+1)$};
      \node (a1selfint) at (-1,0.5) {$-3$}; 
    \end{tikzpicture}
  } 
  \caption*{Figure \ref{sec:pf.of.main.th.2}-4-1C}   
\end{figure}
Then, we obtain an explicit expression of 
the canonical divisor of $X$.
\[
  K_{X} \sim 
  2(O)+3mF+A_{0}
  +\pi^{*}(E_{\mathrm{can}}-E_{\mathrm{fin}}^{(2)}),   
\]
where $F$ is a general fiber. 

\medskip 

\paragraph*{Subcase 4-2 \texorpdfstring{$d=10m+9$ $(n=2m+1)$}{d=10m+9}}\quad\\ 

We have 
\begin{align*}
  \overline{B} 
  &= \overline{C}+\overline{M_{\infty}}+\overline{L_{\infty}}
  +E_{1}+ \dots + E_{2m+1}+E_{\mathrm{fin}}^{(1)}, \\
  \overline{Z} 
  &= 3\overline{M_{\infty}}+(5m+5)\overline{L_{\infty}}
  +(5m+5)(E_{1}+ \dots + E_{2m+1}) \\ 
  &+(5m+5)E_{2m+2}-E_{\mathrm{fin}}^{(2)}, \\
  K_{\overline{\Sigma}}+\overline{Z}
  &\sim \overline{M_{\infty}}+(5m+3)\overline{L_{\infty}} \\
  &+(5m-2)E_{1}+(5m-1)E_{2}+ \dots +(3m+2)E_{2m+1} \\ 
  &+(3m+1)E_{2m+2}+(E_{\mathrm{can}}-E_{\mathrm{fin}}^{(2)}), \\
  p_{a}(\overline{Z}) 
  &= 4m+2+\mu, \\
  (K_{\overline{\Sigma}}+\overline{Z})^{2} 
  &= 2m-2+\nu. 
\end{align*}
Reflecting the data about the branch locus, 
we have the configuration as in Figure \ref{sec:pf.of.main.th.2}-4-2A. 
\begin{figure}[H]
  \centering
  \begin{tikzpicture}
    \draw (-0.5,0.3)node[above]{$\overline{C}$} 
    to [out=315,in=225] (0.25,0.25) 
    to [out=45,in=315] (0.3,0.8); 
    \node (intersection) at (4,0) {$(\overline{C} \cdot E_{2m+2})=4$}; 
    \draw[dashed] (1,1)node[above]{$E_{2m+2}$}--(0,0);  
    \draw (0.5,1)--(1.5,0)node[right]{$\overline{M_{\infty}}$}; 
    \draw (-0.5,1)--(0.5,0)node[below]{$E_{2m+1}$};   
    \node (tenten) at (-1,0.5) {$\cdots$}; 
    \draw (-1.5,1)--(-2.5,0)node[below]{$E_{2}$};  
    \draw (-2,0)--(-3,1)node[above]{$E_{1}$}; 
    \draw (-2.5,1)--(-3.5,0)node[below]{$\overline{L_{\infty}}$};  
    \node (Lselfint) at (-3.5,-0.8) {$\Circled{-1}$}; 
    \node (E1selfint) at (-3,1.7) {$\Circled{-2}$}; 
    \node (E2selfint) at (-2.5,-0.7) {$\Circled{-2}$}; 
    \node (Enselfint) at (0.5,-0.7) {$\Circled{-2}$}; 
    \node (En+1selfint) at (1,1.7) {$\Circled{-1}$}; 
    \node (Mselfint) at (2,-0.6) {$\Circled{-(2m+2)}$}; 
  \end{tikzpicture} 
  \caption*{Figure \ref{sec:pf.of.main.th.2}-4-2A}  
\end{figure} 
After some blowing-ups, 
we have the dual graph on $S$ by taking a double covering of $\Sigma$ as in Figure \ref{sec:pf.of.main.th.2}-4-2B,  
\begin{figure}[H] 
  \centering
  \scalebox{0.8}{
  \begin{tikzpicture}
    \node[draw,shape=circle,inner sep=2pt,fill=black] (linf) at (-3,0) {};
    \node[draw,shape=circle,inner sep=2pt] (llinf) at (-2,0) {};
    \node[draw,shape=circle,inner sep=2pt,fill=black] (e1) at (-1,0) {};  
    \node[draw,shape=circle,inner sep=2pt] (ee1) at (0,0) {}; 
    \node (ueyoko) at (1,0) {$\cdots$}; 
    \node[draw,shape=circle,inner sep=2pt,fill=black] (e2m) at (2,0) {};
    \node[draw,shape=circle,inner sep=2pt] (ee2m) at (3,0) {};
    \node[draw,shape=circle,inner sep=2pt,fill=black] (e2m+1) at (4,0) {}; 
    \node[draw,shape=circle,inner sep=2pt] (l) at (5,0) {};
    \node (e2m+2) at (6,0) {$\triangle$}; 
    \node[draw,shape=circle,inner sep=2pt,fill=black] (m) at (7,0) {};     
    \draw (linf)--(llinf)--(e1)--(ee1);     
    \draw (ee1)--($(ee1)+(0.5,0)$);
    \draw ($(e2m)+(-0.5,0)$)--(e2m);    
    \draw (e2m)--(ee2m)--(e2m+1)--(l);
    \draw[double distance=2pt,nfold=2] (l)--($(e2m+2)+(-0.1,0)$);  
    \draw ($(e2m+2)+(0.1,0)$)--(m);  
    \node (linfname) at ($(linf)+(0,-0.5)$) {$\widetilde{L_{\infty}}$}; 
    \node (llinfname) at($(llinf)+(0,-0.5)$)  {$\widetilde{L_{\infty}'}$}; 
    \node (e1name) at ($(e1)+(0,-0.5)$) {$\widetilde{E_{1}}$};   
    \node (ee1name) at ($(ee1)+(0,-0.5)$) {$\widetilde{E_{1}'}$}; 
    \node (e2mname) at ($(e2m)+(0,-0.5)$) {$\widetilde{E_{2m}}$};  
    \node (ee2mname) at ($(ee2m)+(0,-0.5)$) {$\widetilde{E_{2m}'}$}; 
    \node (e2m+1name) at ($(e2m+1)+(0,-0.5)$) {$\widetilde{E_{2m+1}}$};
    \node (lname) at ($(l)+(0,-0.5)$) {$L$}; 
    \node (e2m+2name) at ($(e2m+2)+(0,-0.5)$) {$\widetilde{E_{2m+2}}$}; 
    \node (mname) at ($(m)+(0,-0.5)$) {$\widetilde{M_{\infty}}$};  
    \node (linfselfint) at ($(linf)+(0,0.5)$) {$-1$}; 
    \node (llinfselfint) at ($(llinf)+(0,0.5)$) {$-2$};   
    \node (e1selfint) at ($(e1)+(0,0.5)$) {$-2$}; 
    \node (ee1selfint) at ($(ee1)+(0,0.5)$) {$-2$}; 
    \node (e2mselfint) at ($(e2m)+(0,0.5)$) {$-2$};   
    \node (ee2mselfint) at ($(ee2m)+(0,0.5)$) {$-2$}; 
    \node (e2m+1selfint) at ($(e2m+1)+(0,0.5)$) {$-2$};  
    \node (lselfint) at ($(l)+(0,0.5)$) {$-2$}; 
    \node (e2m+2selfint) at ($(e2m+2)+(0,0.5)$) {$-4$}; 
    \node (mselfint) at ($(m)+(0,0.5)$) {$-(m+1)$};  
  \end{tikzpicture}  
  } 
  \caption*{Figure \ref{sec:pf.of.main.th.2}-4-2B}   
\end{figure} 
\noindent where  
$\pi^{*}(\overline{L_{\infty}})=2\widetilde{L_{\infty}}+\widetilde{L_{\infty}'}$, 
$\pi^{*}(E_{1})=\widetilde{L_{\infty}'}+2\widetilde{E_{1}}+\widetilde{E_{1}'}$, 
$\pi^{*}(E_{i})=\widetilde{E_{i-1}'}+2\widetilde{E_{i}}+\widetilde{E_{i}'}$ for $2 \le i \le 2m$, 
$\pi^{*}(E_{2m+1})=\widetilde{E_{2m}'}+2\widetilde{E_{2m+1}}+L$, 
$\pi^{*}(E_{2m+2})=L+2\widetilde{E_{2m+2}}$,
$\pi^{*}(\overline{M_{\infty}})=2\widetilde{M_{\infty}}$, 
$\triangle$ stands for an irreducible rational curve $\widetilde{E_{2m+2}}$ with 
an ordinary cusp of multiplicity $2$, 
$L$ intersects $\widetilde{E_{2m+2}}$ at 
the cusp point with $(L \cdot \widetilde{E_{2m+2}})=2$, 
and $\widetilde{M_{\infty}}$ intersects $\widetilde{E_{2m+2}}$ transversally at a nonsingular point.

By contracting 
$\widetilde{L_{\infty}}$, $\widetilde{L_{\infty}'}$, $\dots$, $\widetilde{E_{2m+1}}$, $L$ 
and 
other $(-1)$-curves contained in reducible fibers in the affine part,        
we have the relatively minimal model $X$. 
By relabeling suitably, we have the dual graph on $X$ as in Figure \ref{sec:pf.of.main.th.2}-4-2C.     
\begin{figure}[H]
  \centering
  \scalebox{0.8}{
    \begin{tikzpicture} 
      \node[draw,shape=circle] (0) at (1,0) {};
      \node (a0) at (0,0) {$\triangle$}; 
      \draw (0)--($(a0)+(0.1,0)$);  
      \node (0name) at (1,-0.5) {$(O)$};
      \node (a0name) at (0,-0.5) {$A_{0}$};
      \node (0selfint) at (1,0.5) {$-(m+1)$};
      \node (a0selfint) at (0,0.5) {$0$}; 
    \end{tikzpicture}
  } 
  \caption*{Figure \ref{sec:pf.of.main.th.2}-4-2C}   
\end{figure}
Then, we obtain an explicit expression of 
the canonical divisor of $X$.
\[
  K_{X} \sim 
  2(O)+(3m+1)F
  +\pi^{*}(E_{\mathrm{can}}-E_{\mathrm{fin}}^{(2)}),   
\]
where $F$ is a general fiber. 

\medskip 
   
In global calculations, we note that $\varphi'(t)$ 
may have a root $\alpha \neq 0$ such that $e_{\alpha} = 5$. If $\varphi'(t)$ has such a root $\alpha$, 
we need to compute the local contribution after the coordinate change mapping $\alpha$ to $0$.  
By Lemma \ref{lem:Artin}, we can compute $p_{a}(X)$ and $(K_{X}^{2})$ from the data on branch loci.   
Therefore, we obtain the desired formulas by 
combining calculations in this section and the previous secction (see also the definition of $\mu$ anf $\nu$ at the end of Section \ref{sec:pf.of.main.th.1}). 
This completes the proof of Theorem \ref{main.th.2}. $\qed$ 

\medskip 

As we saw in the proof of Theorem \ref{main.th.2}, 
it always holds that 
$\overline{C}, \overline{M_{\infty}} \subseteq \Supp(\overline{B})$ 
independently of the degree $d$.      
Thus, we have $\sigma'(\overline{C}), \sigma'(\overline{M_{\infty}}) \subseteq \Supp(B)$. 
By Lemma \ref{double covering},  
we see that $\rho^{*}(\sigma'(\overline{C}))=2\Gamma$ and 
$\rho^{*}(\sigma'(\overline{M_{\infty}}))=2D$, 
where $\Gamma$ and $D$ are nonsingular curves.  

\begin{prop.}
As the above notations, 
$\Gamma$ is a moving cusp and $D$ is a section of $f$.  
\end{prop.}

\begin{proof}
We can easily check them by using Lemma \ref{double covering}. 
\end{proof}  

Finally, we give an explicit formula for the self-intersection number of the moving cusp $\Gamma$. 
We denote $n_{i}$ by the number of reducible fibers of corresponding type in the affine part as below:   
\begin{table}[H] 
  \scalebox{0.9}{ 
  \centering 
  \begin{tabular}{|c|c|c|c|c|c|c|c|} 
  \hline
  type &
  $C(5,1)$ & $C(9,2)$ & $C(3,3)$ & 
  $C(11,5)$ & $C(4,6)$ & $C(9,7)$ & 
  $C(13,8)$
  \\ \hline   
  the number & $n_{1}$ & $n_{2}$ & $n_{3}$ & $n_{5}$ & $n_{6}$ & $n_{7}$ & $n_{8}$  
  \\ \hline
  \end{tabular}
  }
\end{table} 
\noindent Note that $N_{i}$ in Theorem \ref{main.th.2} represents the number of reducible fibers of type $C(-,i)$ including the fiber at infinity, 
but $n_{i}$ does not count the fiber at infinity. 
With these notations, we can give a formula for $(\Gamma^{2})$. 

\begin{cor.} \label{selfint of the moving cusp} 
  For a relatively minimal unirational quasi-hyperelliptic fibration $f : X \to \P^{1}$ 
  of genus $2$ in characteristic $5$ with a section, we have     
\[
  (\Gamma^{2})=4n_{3}+5n_{5}+7n_{6}+8n_{7}+9n_{8}-(15m+\delta(d)), 
\]
where 
\[
  \delta(d) \coloneqq 
  \begin{cases}
    2 & \text{if} \ d=10m+1,10m+4, \\
    -1 & \text{if} \ d=10m+2, \\
    0 & \text{if} \ d=10m+3, \\
    3 & \text{if} \ d=10m+6, \\
    7 & \text{if} \ d=10m+7,10m+8,10m+9.  
    \end{cases}
\] 
\end{cor.} 

\begin{proof}
  By the adjunction formula, we have 
  \[
    (\Gamma^{2}) = -(\Gamma \cdot K_{X}) +2p_{a}(\Gamma)-2 = -(\Gamma \cdot K_{X}) -2. 
  \]
  Using the expression of the canonical divisor obtained in the proof of Theorem \ref{main.th.2}, 
  we get 
  {\small 
  \[ 
  -(\Gamma \cdot K_{X}) = 
  (\Gamma \cdot \pi^{*}(E_{\mathrm{can}}-E_{\mathrm{fin}}^{(2)})) 
  + 
  \begin{cases}
    -15m & \text{if} \ d=10m+1,10m+4, \\
    -15m+3 & \text{if} \ d=10m+2, \\
    -15m+2 & \text{if} \ d=10m+3, \\
    -15m-1 & \text{if} \ d=10m+6, \\
    -15m-5 & \text{if} \ d=10m+7,10m+8,10m+9.  
  \end{cases} 
  \]
  }  
  Since we already computed $\pi^{*}(E_{\mathrm{can}}-E_{\mathrm{fin}}^{(2)})$ in the proof of Theorem \ref{main.th.1}, we get the desired formula.  
\end{proof}

\section{Rational quasi-hyperelliptic surfaces} \label{sec:rat QE} 

In this section, we classify rational quasi-hyperelliptic surfaces in characteristic $5$.  

\begin{th.} \label{main.th.3}
  Let $f : X \to \P^{1}$ be a rational quasi-hyperelliptic fibration of genus $2$ in characteristic $5$ which is birationally equivalent to an affine hypersurface defined by $y^{2} = x^{5} + \varphi(t)$ satisfying Condition \ref{assumption}. 
  Then, $X$ is classified into the $10$ cases as in Table \ref{table:rational},   
  where $\rho$ is the Picard number of $X$ and 
  $r$ is the torsion-rank of the Mordell-Weil group 
  (cf.\ Theorem \ref{th:MWG of rat q-HE} and Corollary \ref{th:str of MWG}). 
  By using automorphisms of $\P^{1}$, we can assume that 
  three reducible fibers are located over $\infty, 0$ and $1$ of $\P^{1}$.     
  Under this assumption, we give defining equations of affine hypersurfaces birationally equivalent to each $X$.     
  \begin{landscape}  
  \begin{table}[H]  
  \centering 
  \scalebox{0.85}{
    \begin{tabular}{|c|l|c|c|l|l|}  
    \hline 
    No. & singular fibers & $\rho$ & $r$ & the defining equation &
    sections other than $(O)$ 
    \\ \hline 
    1 & $C(13,8) \times 1$ & $14$ & $0$ & $y^2 = x^{5}+t$ & 
    not exist      
    \\ \hline
    2 & $C(9,7) \times 1$, $C(5,1) \times 1$
    & $14$ & $1$ & 
    $y^{2} = x^{5}+t^{2}$ & $(x,y) = (0,\pm t)$
    \\ \hline 
    3 & $C(4,6) \times 1$, $C(9,2) \times 1$ & $13$ & $0$ &
    $y^{2} = x^{5}+t^{3}$ & 
    not exist 
    \\ \hline 
    4 & $C(4,6) \times 1$, $C(5,1) \times 2$ & $13$ & $1$ &
    $y^{2} = x^{5}+t^{3}+t^{2}$ & not exist  
    \\ \hline 
    5 & $C(11,5) \times 1$, $C(3,3) \times 1$ & $14$ & $1$ & 
    $y^{2} = x^{5}+t^{4}$ & $(x,y) = (0, \pm t^{2})$
    \\ \hline 
    6 & $C(3,3) \times 2$, $C(5,1) \times 2$ & $14$ & $2$ & 
    $y^{2} = x^{5}+t^{6}+t^{4}$ & 
      \begin{tabular}{r}
        $(x,y) = 
        (2t , \pm (t^{3} + t^{2}))$, \\  
        $(-2t , \pm (t^{3} - t^{2}))$   
      \end{tabular}   
    \\ \hline  
    7 & $C(3,3) \times 6$ &  $14$ & $3$ & 
    \begin{tabular}{l}  
    $y^{2} = x^{5} + \varphi(t)$  
    \ \text{with} \ \\  
    $\varphi'(t) = t^{3}(t-1)^{3}(t-\alpha)^{3}(t-\beta)^{3}(t-\gamma)^{3}$, \\ 
    \ \text{where} \ $\alpha, \beta, \gamma$ 
    \ \text{are distinct and neither 0 nor 1}. \\  
    \end{tabular}   
    &  
    \hspace{3pt} When $\{ \alpha,\beta,\gamma \} \neq \{ 2,3,4 \}$, not exist. 
    \\ \cline{6-6} 
    & & & & & 
    \begin{tabular}{l} 
    When $\{ \alpha, \beta, \gamma \} = \{ 2,3,4 \}$),  
    i.e., $\varphi(t) = t^{16}+t^{12}+t^{8}+t^{4}$, \\  
    $(x,y) = 
    (t^{4}, \pm (t^{5}-t)^{2})$, \\     
    $\left( 4it^{3}+4i^{2}t^{2}+4i^{3}t+4i^{4}+1 , \pm \left(\frac{t^{5}-t}{t-i}\right)^{2} \right)$ \\ for $i \in \mathbb{F}_{5} = \{ 0,1,2,3,4 \}$. \\  
    \end{tabular}   
    \\ \hline   
    8 & $C(3,3) \times 4$, $C(4,6) \times 1$ & $13$ & $2$ & 
    \begin{tabular}{l}  
    $y^{2} = x^{5} + \varphi(t)$  
    \ \text{with} \ \\  
    $\varphi'(t) = t^{3}(t-1)^{3}(t-\alpha)^{3}(t-\beta)^{3}$, \\ 
    \ \text{where} \ $\alpha, \beta$ 
    \ \text{are distinct and neither 0 nor 1}.   
    \end{tabular}   
    &  not exist
    \\ \hline   
    9 & $C(3,3) \times 2$, $C(4,6) \times 2$ & $12$ & $1$ & 
    $y^{2} = x^{5}+t^{13}+t^{11}+t^{9}+t^{7}$ 
    &  not exist
    \\ \hline   
    10 & $C(4,6) \times 3$ & $11$ & $0$ & 
    $y^{2} = x^{5}+t^{13}+t^{12}+4t^{8}+4t^{7}$  
    & not exist
    \\ \hline   
    \end{tabular}  
  } 
  \caption{Rational quasi-hyperelliptic surfaces} \label{table:rational}
  \end{table} 
  \end{landscape}  
\end{th.} 

\subsection{The bound of the Picard number}  

In this subsection, we give a bound of the Picard number of a rational quasi-hyperelliptic surface. 
The same bound is obtained in \cite{SS94} for 
a fibration of genus $\ge 1$ in characteristic $0$. 
In fact, we do not need the bound to prove 
the classification theorem \ref{main.th.3}, 
but we believe that the bound will play an important role  
to treat a rational quasi-hyperelliptic fibration in characteristic $p$.  

It is known that the following inequality so called Xiao's inequality holds. 

\begin{th.}[{\cite[Theorem 1.1]{GSZ23}}]   
  Let $X$ (resp.\ $C$) be a nonsingular projective algebraic surface (resp.\ curve) 
  defined over an algebraically closed field $k$ of characteristic $p>0$. 
  Let $f : X \to C$ be a relatively minimal fibration of arithmetic genus $g \ge 2$.  
  We assume that one of the following conditions; 
  \begin{enumerate}[label=(\roman*)]
    \item The generic fiber of $f$ is hyperelliptic or quasi-hyperelliptic. 
    \item The generic fiber of $f$ is nonsingular. 
    \item The genus of $C$ is less than or equal to 1.   
  \end{enumerate} 
  Then, the following inequality holds;  
  \[
    (K_{X/{C}}^{2}) \ge \frac{4g-4}{g}\chi_{f},    
  \]
  where $K_{X/{C}} = K_{X}-f^{*}(K_{C})$, $\chi_{f}=\deg(f_{*}K_{X/{C}})$. 
\end{th.}

By Xiao's inequality, we can obtain a bound of the Picard number of 
a rational quasi-hyperelliptic surface in characteristic $p \ge 5$. 

\begin{lem.} \label{bound of the Picard number for rational surfaces}
  Let $f : X \to \P^{1}$ be a relatively minimal 
  quasi-hyperelliptic fibration of genus $g=(p-1)/2$ in characteristic $p \ge 5$.    
  If $X$ is rational, 
  then the Picard number $\rho \coloneqq \rho(X)$ of $X$ satisfies 
  $\rho \le 4g+6 = 2p+4$.     
\end{lem.} 

\begin{proof}  
  Note that the following relations holds for a genus $g$ fibration 
  $f : X \to C$ by the Riemann-Roch theorem;   
  \begin{align*} 
    (K_{X/{C}}^{2}) &= (K_{X}^{2})-8(g-1)(g(C)-1), \\
    \chi_{f} &= \chi(\O_{X})-(g-1)(g(C)-1) +\dim_{k}(R^{1}f_{*}\O_{X})_{\mathrm{tor}}, 
  \end{align*} 
  \noindent where $g(C)$ is the genus of $C$.    
  Since $X$ is rational, we have $\chi(\O_{X})=1$, $b_{1}=0$ and  $b_{2}=\rho$, 
  where each $b_{i}$ is the $i$-th Betti number for $i=1,2$.     
  By Noether's formula, we have $12 = (K_{X}^{2})+e(X)$, 
  where $e(X)$ is the topological Euler number of $X$. 
  Since $e(X) = 2+\rho$, we get $\rho = 10 - (K_{X}^{2})$. 
  Combining the above remarks and Xiao's inequality, 
  we obtain 
  \[
  (K_{X}^{2}) \ge \frac{4g-4}{g}(\dim_{k}(R^{1}f_{*}\O_{X})_{\mathrm{tor}}-g) \ge -(4g-4).   
  \]
  Thus, we have   
  \[ 
    \rho = 10 - (K_{X}^{2}) \le 10+(4g-4) \le 4g+6.  
  \]
\end{proof}  

\subsection{The combinatorial possibilities} 

Let $f : X \to \P^{1}$ be a unirational quasi-hyperelliptic fibration 
of genus 2 in characteristic $5$ with a section, and  
assume that $X$ is birationally equivalent to 
the affine hypersurface defined by $y^{2} = x^{5} + \varphi(t)$ satisfying Condition \ref{assumption}. 
Put $d = \deg \varphi(t)$ and $m = [d/{10}]$.   
In this subsection, 
we determine combinations of reducible fibers 
which a rational quasi-hyperelliptic fibration can have.  

\begin{lem.} \label{rational candidates} 
The value $m$ and the combinations of reducible fibers for which X can be rational are as follows.
 \begin{enumerate}[label=(\roman*)]
    \item $m=0$ ; $C(13,8) \times 1$. 
    \item $m=0$ ; $C(9,7) \times 1$, $C(5,1) \times 1$. 
    \item $m=0$ ; $C(4,6) \times 1$, $C(9,2) \times 1$. 
    \item $m=0$ ; $C(4,6) \times 1$, $C(5,1) \times 2$. 
    \item $m=0$ ; $C(11,5) \times 1$, $C(3,3) \times 1$. 
    \item $m=0$ ; $C(3,3) \times 2$, $C(9,2) \times 1$. 
    \item  $m=0$ ; $C(3,3) \times 2$, $C(5,1) \times 2$.  
    \item  $m=1$ ; $C(3,3) \times 6$. 
    \item  $m=1$ ; $C(3,3) \times 4$, $C(4,6) \times 1$. 
    \item  $m=1$ ; $C(3,3) \times 2$, $C(4,6) \times 2$. 
    \item  $m=1$ ; $C(4,6) \times 3$. 
  \end{enumerate}
\end{lem.} 

\begin{proof} 
First, we can assume that 
the fiber at $t=\infty$ is 
of type $C(0,1)$ by 
using an automorphism of $\P^{1}$. 
Then we have $d = 10m+9$ ($m \in \Z_{\ge 0}$). 
We denote $n_{i}$ by  
the numbers of reducible fibers of corresponding type in the affine part 
as below:   
\begin{table}[H] 
  \scalebox{0.9}{ 
  \centering 
  \begin{tabular}{|c|c|c|c|c|c|c|c|} 
  \hline
  type &
  $C(5,1)$ & $C(9,2)$ & $C(3,3)$ & 
  $C(11,5)$ & $C(4,6)$ & $C(9,7)$ & 
  $C(13,8)$
  \\ \hline   
  the number & $n_{1}$ & $n_{2}$ & $n_{3}$ & $n_{5}$ & $n_{6}$ & $n_{7}$ & $n_{8}$  
  \\ \hline
  \end{tabular}
  }
\end{table} 
Focusing on the multiplicities of roots of 
$\varphi'(t)=0$, we have the equation 
\begin{equation} \label{eq:multiplicity in affine part} 
  n_{1}+2n_{2}+3n_{3}+5n_{5}+6n_{6}+7n_{7}+8n_{8}=d-1
  =10m+8. 
\end{equation}   
Since we have $p_{g}(X)=0$ according to Castelnuovo's criterion of rationality, and $q(X)=0$ for a unirational quasi-hyperelliptic surface $X$ by Lemma \ref{irregularity}, we must have $p_{a}(X)=0$.   
By Theorem \ref{main.th.2}, we have 
\begin{equation} \label{eq:arithmetic genus zero} 
n_{3} + n_{5} + 2n_{6} + 2n_{7} + 2n_{8} = 4m+2. 
\end{equation} 
Combining equations (\ref{eq:multiplicity in affine part}) and (\ref{eq:arithmetic genus zero}), we have 
\[
0 \le n_{1}+2n_{2}+2n_{5}+n_{7}+2n_{8} = -2m+2.  
\]
In particular, we have $m \le 1$. 
The remaining combinatorial computation 
is left to the reader. 
\end{proof} 

\subsection{Existence}  

In this subsection, we determine the rational ones 
from the candidates (i), $\dots$, (xi) which we listed in the previous section. 

\begin{lem.} 
The surfaces listed in Lemma \ref{rational candidates} 
exist except for (vi). 
Furthermore, the defining equations 
are uniquely determined up to the automorphisms of $\P^{1}$
except for (vi), (viii) and (xi).   
\end{lem.}  

\begin{proof}
  Note that up to three, we can assume that reducible fibers emerge over $\infty, 0$ and $1$ of $\P^{1}$ 
  by using automorphisms of $\P^{1}$. 
  Under this assumption, we will check existence and non-existence. 
  
  First, we show that (vi) does not exist.     
  Suppose that a polynomial $\varphi(t)$ realizing (vi) exists. Then $\varphi'(t)$ has exactly two roots $0$ and $1$ whose multiplicities are $3$ and $2$, respectively, 
  and $\varphi'(t)$ has a term of degree $4$, but 
  this is impossible because of Condition \ref{assumption} 
  
  Next, we check the existence for the remaining cases. 
  For (i), (ii), (iii), (iv), and (v), 
  we can assume that the degree of $\varphi'(t)$ is less than $4$ by choosing the fiber at infinity appropriately. Therefore, we get the polynomial $\varphi(t)$ realizing each surface.  
  For (xi), since $\varphi'(t) = t^{6}(t-1)^{6}$ and it does not have any terms of degree $4$ modulo $5$, we can obtain $\varphi(t)$ realizing (xi). 
  For (v), we start with $\varphi'(t) = t^{3}(t-1)(t-\alpha)$, where $\alpha \in k$ is neither $0$ nor $1$ 
  (the fiber at infinity is of type $C(3,3)$). 
  The term of degree $4$ of $\varphi'(t)$ is $-(\alpha +1)t^{4}$, 
  so we must have $\alpha = 4$. 
  By the same computation, we obtain the defining equation of (x). 
  In summary, the surfaces such that the number of reducible fibers is less than $5$ except for (vi) exist and the defining equations are rigid. 

  For (viii), we start with 
  \[ 
  \varphi'(t) = t^{3} (t-1)^{3} (t-\alpha)^{3} (t-\beta)^{3} (t-\gamma)^{3}, 
  \]
  where $\alpha, \beta$ and $\gamma$ are distinct from each othe and neither $0$ nor $1$.  
  Choosing $\alpha, \beta$ and $\gamma$ such that 
  all coefficients of the terms of degree $4,9$ and $14$ 
  of $\varphi'(t)$ vanish 
  (take $\alpha = 2, \beta =3, \gamma=4$ for example), 
  we get $\varphi(t)$ realizing (viii). 
  In the same fashion, we can check the existence of (ix). 
\end{proof} 

The following lemma is a well-known result 
due to Castelnuovo's criterion for rationality.   

\begin{lem.} \label{rationality criterion} 
  Let $X$ be a nonsingular surface with $q(X)=0$.  
  If there exists a rational curve $C$ of genus zero on $X$ with nonnegative self-intersection number, then $X$ is rational.   
\end{lem.}

\begin{lem.}
The surfaces listed in 
Lemma \ref{rational candidates} are rational, except for (vi).
\end{lem.} 

\begin{proof} 
Recall that the irregularity of each surface is zero 
because of the unirationality (Lemma \ref{irregularity}). 
By the previous Lemma \ref{rationality criterion}, 
we enough to find 
a rational curve of genus zero with nonnegative self-intersection number. 
We can take the moving cusp $\Gamma$ as a desired curve.  
Indeed, we can check $(\Gamma^{2})$ is nonnegative 
by Corollary \ref{selfint of the moving cusp} except for (xi). 
For (xi), we can find such a curve after contracting $(-1)$-curves three times.     
\end{proof} 

This completes the proof of Theorem \ref{main.th.3} and 
we can conclude that rational quasi-hyperelliptic surfaces 
are classified into 10 cases.   

\subsection{Sections}   

In this subsection, we observe 
the configurations of sections and reducible fibers.    
For the cases where the torsion-rank of the Mordell-Weil group is zero, we have nothing to do.  
Therefore, we treat No.\ 2,4,5,6,7,8 and 9 in Table \ref{table:rational}.   

\begin{prop.} \label{prop:non-existence of rat pt}
For No.\ 4, 8 and 9, 
$f : X \to \P^{1}$ has no nonzero sections.     
\end{prop.} 

\begin{proof} 
We only prove for No.\ 4 because proofs for other two cases are same. 
Suppose that there exists a nontrivial 
$K$-rational point on $X_{\eta} \setminus \{ P_{\infty} \}$, say $P$. 
Using the explicit formula for the height pairing (cf.\ the equation (\ref{explcit height pairing}) in Section \ref{sec:MWG}), 
the adjunction formula and 
the expression of the canonical divisor (cf.\ Section \ref{sec:pf.of.main.th.2}), we have 
\begin{align*} 
  \langle \iota(P), \iota(P) \rangle 
  &= -(P^{2}) +2(P \cdot O) -(O^{2}) - \sum_{v \in R}\mathrm{contr}_{v}(P) \\ 
  &= (P \cdot K_{X}) + 2 + 2(P \cdot O) +1 - \sum_{v \in R}\mathrm{contr}_{v}(P) \\ 
  &= 4(P \cdot O) + 8 - \sum_{v \in R} \mathrm{contr}_{v}(P).   
\end{align*} 
No.\ 4 has one reducible fiber of type $C(4,6)$ and two reducible fibers of type $C(5,1)$, 
so the sum of local contributions cannot exceed $8$ (cf.\ Table \ref{table:contribution}). 
Hence we have $\langle \iota(P),\iota(P) \rangle > 0$, 
but this contradicts the fact that $P$ is a torsion element.
\end{proof} 

For No.\ 2,5 and 6, we can easily find nontrivial $K$-rational points as in Table \ref{table:rational} 
and we can check that those are all of the $K$-rational points of the generic fiber by the same calculation as in \cite[Section 5]{Ito92}. 

\begin{prop.} 
For No.\ 7, that is, 
a rational quasi-hyperelliptic fibration $f : X \to \P^{1}$
which has exactly six reducible fibers of type $C(3,3)$, 
it has nonzero sections if and only if 
it is birationally equivalent to 
the affine hypersurface defined by 
$y^{2} = x^{5} + t^{16}+t^{12}+t^{8}+t^{4}$. 
Then nonzero sections are as follows.   
\begin{itemize}
    \item $(x,y) = (t^{4}, \pm (t^{5}-t)^{2})$,   
    \item $(x,y) = 
    \left( 4it^{3}+4i^{2}t^{2}+4i^{3}t+4i^{4}+1, \pm \left(\frac{t^{5}-t}{t-i}\right)^{2} \right)$ for $i \in \mathbb{F}_{5} = \{ 0,1,2,3,4 \}$. 
\end{itemize} 
\end{prop.} 

\begin{proof}
Let $(x,y)$ be a nonzero section.  
By the same way as Proposition \ref{prop:non-existence of rat pt}, we can see that a section is integral if it exists. 
Thus we consider the equation 
\[
  2y(t)y'(t) = \varphi'(t),  
\]
for $y(t) \in k[t]$. 
Put $\varphi'(t) = g(t)^{3}$, where 
$g(t) = \prod_{i=1}^{5}(t-\alpha_{i})$  
is a polynomial whose roots $\alpha_{i} \in k$ 
are distinct each other.  
Since $\ord_{t=\alpha_{i}}(y')$ must be 
equal to $\ord_{t=\alpha_{i}}(y)-1$ 
if $\ord_{t=\alpha_{i}}(y) \neq 0$, 
we have 
$(\ord_{t=\alpha_{i}}(y),\ord_{t=\alpha_{i}}(y')) = (2,1)$ or $(0,3)$ for each $i$. 
It holds for the order at $t=\infty$ too. 
This implies that 
$(\deg(y),\deg(y')) = (8,7)$ or $(10,5)$ because the order of $y$ (resp.\ $y'$) at the infinity is equal to $10-\deg(y)$ (resp.\ $8-\deg(y')$).  
For five out of the six fibers over $t= \alpha_{1}, \dots , \alpha_{5}$ and $\infty$, 
the orders of $y$ and $y'$ are $2$ and $1$ respectively.  
For the remaining one, those are $0$ and $3$ respectively. 
We first consider the case where 
$(\ord_{t=\infty}(y),\ord_{t=\infty}(y'))
=(0,3)$ and 
$(\ord_{t=\alpha_{i}}(y),\ord_{t=\alpha_{i}}(y')) = (2,1)$ for $i=1, \dots , 5$. 
In this case, we have 
\[
y = cg(t)^{2}, \ \ y' = 2cg(t)g'(t) \ \ (c \in k^{\times}),  
\]
and $2yy' = 4c^{2}g^{3}g'$. 
Since $2yy' = \varphi'$, 
we must have $g'(t) = -c^{-2}$, 
that is, $g(t) = t^{5} -c^{-2}t + d$ ($d \in k$). 
By suitable coordinate changes,  
we obtain $g(t) = t^{5} - t$. 
For the case where 
$(\ord_{t=\alpha_{i}}(y),\ord_{t=\alpha_{i}}(y')) = (0,3)$ for some $i$, 
we can reduce this case to the former case by the coordinate change mapping $\alpha_{i}$ to $\infty$. 

To conclude the above discussion, 
we know that $g(t)=t^{5}-t$ if $f$ has a nonzero section, and there is no nonzero sections of $f$ except for the twelve sections in the statement.   
\end{proof} 

\begin{rem.} 
By the above observations, 
we see that all of the sections of rational quasi-hyperelliptic fibration
are integral. 
\end{rem.} 

\begin{rem.}
For No.\ 2,5,6 and 7, each Mordell-Weil group is generated by the image of nontrivial $K$-rational points of 
the generic fiber under the closed immersion $\iota$ 
while it is not true for the others. 
\end{rem.}

Finally, we describe the configuration for No.\ 5 as an example (Figure 9-1).         

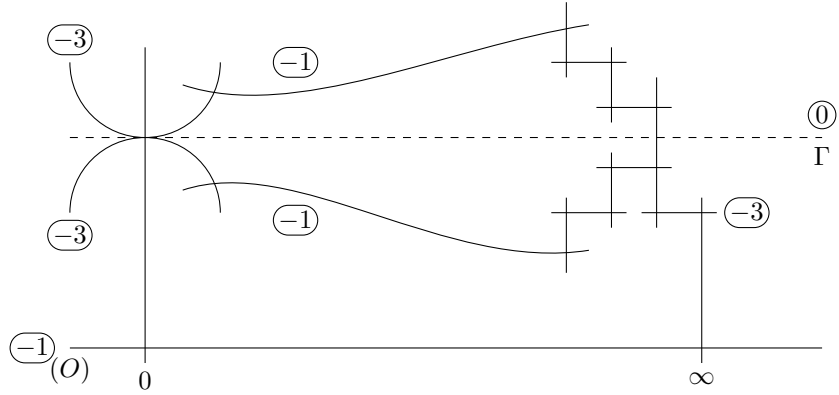
\begin{figure}[H] 
  \centering  
  \begin{tikzpicture}
    \draw (-4,-1) node[below] {$(O)$} --(6,-1); 
    \draw (-3,-1.2) node[below] {$0$} --(-3,3); 
    \draw (4.4,-1.2)node[below] {$\infty$} --(4.4,1); 
    \draw (4.6,0.8)--(3.6,0.8); 
    \draw (3.8,0.6)--(3.8,2.6); 
    \draw (4,1.4)--(3,1.4); 
    \draw (3.2,1.6)--(3.2,0.6); 
    \draw (3.4,0.8)--(2.4,0.8); 
    \draw (2.6,1)--(2.6,0);
    \draw (4,2.2)--(3,2.2); 
    \draw (3.2,2)--(3.2,3); 
    \draw (3.4,2.8)--(2.4,2.8); 
    \draw (2.6,2.6)--(2.6,3.6); 
    \draw[dashed] (-4,1.8)--(6,1.8) node[below] {$\Gamma$}; 
    \draw (-4,2.8) to [out=270,in=180] (-3,1.8) to [out=0,in=270] (-2,2.8); 
    \draw (-4,0.8) to [out=90,in=180] (-3,1.8) to [out=0,in=90] (-2,0.8); 
    \node (a) at (-4,3.1) {$\Circled{-3}$};  
    \node (b) at (-4,0.5) {$\Circled{-3}$};  
    \node (e) at (-4.5,-1) {$\Circled{-1}$}; 
    \node (f) at (6,2.1) {$\Circled{0}$}; 
    \node (g) at (5,0.8) {$\Circled{-3}$};
    \draw (-2.5,2.5) .. controls (-1,2) and (1,3) .. (2.9,3.3); 
    \node (aa) at (-1,2.8) {$\Circled{-1}$}; 
    \draw (-2.5,1.1) .. controls (-1,1.6) and (1,0) .. (2.9,0.3); 
    \node (bb) at (-1,0.7) {$\Circled{-1}$};  
  \end{tikzpicture} 
  \caption*{Figure 9-1: Configuration for No.\ 5} 
\end{figure} 

In Figure 9-1, 
the fiber at $0$ is of type $C(3,3)$ and the fiber at $\infty$ is of type $C(11,5)$. 
There are two sections $(x,y) = (0,t^{2}), (0,-t^{2})$ and 
both self-intersection numbers are equal to $-1$, and    
they are disjoint (we can check these by using the height pairing).  
Each section intersects with the two reducible fibers at each non-identity component 
and does not pass the same non-identity component of each reducible fiber. 

\section*{Acknowledgements} 
The authors thank to Professor Yuya Matsumoto 
for valuable comments on the earlier version of this preprint, which improve our paper.    
The authors also thank to Professor Fuetaro Yobuko for valuable comments. 


\end{document}